\documentclass{article}%
\usepackage{amsmath}
\usepackage{amsfonts}
\usepackage{amssymb}
\usepackage{graphicx}%
\providecommand{\U}[1]{\protect\rule{.1in}{.1in}}
\newtheorem{theorem}{Theorem}

\newtheorem{corollary}[theorem]{Corollary}

\newtheorem{definition}[theorem]{Definition}
\newtheorem{example}[theorem]{Example}

\newtheorem{lemma}[theorem]{Lemma}
\newtheorem{notation}[theorem]{Notation}

\newtheorem{proposition}[theorem]{Proposition}
\newtheorem{remark}[theorem]{Remark}

\newenvironment{proof}[1][Proof]{\noindent\textbf{#1.} }{\ \rule{0.5em}{0.5em}}
\evensidemargin=\oddsidemargin
\newdimen\dummy
\dummy=\oddsidemargin
\begin{document}

\title{A Family of Crouzeix-Raviart Finite Elements in 3D}
\author{Patrick Ciarlet, Jr.\thanks{(patrick.ciarlet@ensta-paristech.fr), POEMS
(CNRS/ENSTA ParisTech/INRIA), 828, Boulevard des Mar{\'{e}}chaux, 91762
Palaiseau Cedex, France}
\and Charles F. Dunkl\thanks{(cfd5z@virginia.edu), Dept. of Mathematics, University
of Virginia, Charlottesville, Virginia 22904-4137, USA}
\and Stefan A. Sauter\thanks{(stas@math.uzh.ch), Institut f\"{u}r Mathematik,
Universit\"{a}t Z\"{u}rich, Winterthurerstr 190, CH-8057 Z\"{u}rich,
Switzerland}}
\maketitle

\begin{abstract}
In this paper we will develop a family of non-conforming \textquotedblleft
Crouzeix-Raviart\textquotedblright\ type finite elements in three dimensions.
They consist of local polynomials of maximal degree $p\in\mathbb{N}$ on
simplicial finite element meshes while certain jump conditions are imposed
across adjacent simplices. We will prove optimal a priori estimates for these
finite elements.

The characterization of this space via jump conditions is implicit and the
derivation of a local basis requires some deeper theoretical tools from
orthogonal polynomials on triangles and their representation. We will derive
these tools for this purpose. These results allow us to give explicit
representations of the local basis functions. Finally we will analyze the
linear independence of these sets of functions and discuss the question
whether they span the whole non-conforming space.

\end{abstract}

\noindent\textbf{AMS-Classification}: Primary 33C45, 33C50, 65N12, 65N30;
Secondary 33C80.

\vspace{0.5cm}


\vspace{0.5cm}

\section{Introduction}

For the numerical solution of partial differential equations, Galerkin finite
element methods are among the most popular discretization methods. In the last
decades, \textit{non-conforming} Galerkin discretizations have become very
attractive where the test and trial spaces are not subspaces of the natural
energy spaces and/or the variational formulation is modified on the discrete
level. These methods have nice properties, e.g. in different parts of the
domain different discretizations can be easily used and glued together or, for
certain classes of problems (Stokes problems, highly indefinite Helmholtz and
Maxwell problems, problems with \textquotedblleft locking\textquotedblright,
etc.), the non-conforming discretization enjoys a better stability behavior
compared to the conforming one. One of the first non-conforming finite element
space was the Crouzeix-Raviart element (\cite{Crouzeix_Raviart}, see
\cite{Brenner_Crouzeix} for a survey). It is piecewise affine with respect to
a triangulation of the domain while interelement continuity is required only
at the barycenters of the edges/facets (2D/3D).

In \cite{ccss_2012}, a family of high order non-conforming (intrinsic) finite
elements have been introduced which corresponds to a family of high-order
Crouzeix-Raviart elements in two dimensions. For Poisson's equation, this
family includes the non-conforming Crouzeix-Raviart element
\cite{Crouzeix_Raviart}, the Fortin-Soulie element \cite{Fortin_Soulie}, the
Crouzeix-Falk element \cite{Crouzeix_Falk}, and the Gauss-Legendre elements
\cite{BaranGL}, \cite{BaranCVD} as well as the standard conforming $hp$-finite elements.

In our paper we will characterize a family of high-order Crouzeix-Raviart type
finite elements in \textit{three} dimensions, first implicitly by imposing
certain jump conditions at the interelement facets. Then we derive a local
basis for these finite elements. These new finite element spaces are
non-conforming but the (broken version of the) continuous bilinear form can
still be used. Thus, our results also give insights on how far one can go in
the non-conforming direction while keeping the original forms.

The explicit construction of a basis for these new finite element spaces
require some deeper theoretical tools in the field of orthogonal polynomials
on triangles and their representations which we develop here for this purpose.

As a simple model problem for the introduction of our method, we consider
Poisson's equation but emphasize that this method is applicable also for much
more general (systems of) elliptic equations.

There is a vast literature on various conforming and non-conforming, primal,
dual, mixed formulations of elliptic differential equations and conforming as
well as non-conforming discretization. Our main focus is the\ characterization
and construction of non-conforming Crouzeix-Raviart type finite elements from
theoretical principles. For this reason, we do not provide an extensive list
of references on the analysis of specific families of finite elements spaces
but refer to the classical monographs \cite{Ciarlet}, \cite{Schwab}, and
\cite{BoffiBrezziFortin} and the references therein.

The paper is organized as follows.

In Section \ref{SecModelProblem} we introduce our model problem, Poisson's
equation, the relevant function spaces and standard conditions on its well-posedness.

In Section \ref{sectconfbasis} we briefly recall classical, conforming
$hp$-finite element spaces and their Lagrange basis.

The new non-conforming finite element spaces are introduced in Section
\ref{section_nonconforme}. We introduce an appropriate compatibility condition
at the interfaces between elements of the mesh so that the non-conforming
perturbation of the original bilinear form is consistent with the local error
estimates. We will see that this compatibility condition can be inferred from
the proof of the second Strang lemma applied to our setting. The weak
compatibility condition allows to characterize the non-conforming family of
high-order Crouzeix-Raviart type elements in an \textit{implicit} way. In this
section, we will also present explicit representations of non-conforming basis
functions of general degree $p$ while their derivation and analysis is the
topic of the following sections.

Section \ref{SecExplConstr} is devoted to the explicit construction of a basis
for these new non-conforming finite elements. It requires deeper theoretical
tools from orthogonal polynomials on triangles and their representation which
we will derive for this purpose in this section.

It is by no means obvious whether the constructed set of functions is linearly
independent and span the non-conforming space which was defined implicitly in
Section \ref{section_nonconforme}. These questions will be treated in Section
\ref{SecPropNC}.

Finally, in Section \ref{SecConclusion} we summarize the main results and give
some comparison with the two-dimensional case which was developed in
\cite{ccss_2012}.

\section{Model Problem\label{SecModelProblem}}

As a model problem we consider the Poisson equation in a bounded Lipschitz
domain $\Omega\subset\mathbb{R}^{d}$ with boundary $\Gamma:=\partial\Omega$.
First, we introduce some spaces and sets of functions for the coefficient
functions and solution spaces.

The Euclidean scalar product in $\mathbb{R}^{d}$ is denoted for $\mathbf{a}%
,\mathbf{b}\in\mathbb{R}^{d}$ by $\mathbf{a}\cdot\mathbf{b}$. For $s\geq0$,
$1\leq p\leq\infty$, let $W^{s,p}\left(  \Omega\right)  $ denote the classical
(real-valued) Sobolev spaces with norm $\left\Vert \cdot\right\Vert
_{W^{s,p}\left(  \Omega\right)  }$. The space $W_{0}^{s,p}\left(
\Omega\right)  $ is the closure with respect to the $\left\Vert \cdot
\right\Vert _{W^{s,p}\left(  \Omega\right)  }$ of all $C^{\infty}\left(
\Omega\right)  $ functions with compact support. As usual we write
$L^{p}\left(  \Omega\right)  $ short for $W^{0,p}\left(  \Omega\right)  $. The
scalar product and norm in $L^{2}\left(  \Omega\right)  $ are denoted by
$\left(  u,v\right)  :=\int_{\Omega}uv$ and $\left\Vert \cdot\right\Vert
:=\left(  \cdot,\cdot\right)  ^{1/2}$. For $p=2$, we use $H^{s}\left(
\Omega\right)  $, $H_{0}^{s}\left(  \Omega\right)  $ as shorthands for
$W^{s,2}\left(  \Omega\right)  $, $W_{0}^{s,2}\left(  \Omega\right)  $. The
dual space of $H_{0}^{s}\left(  \Omega\right)  $ is denoted by $H^{-s}\left(
\Omega\right)  $. We recall that, for positive integers $s$, the seminorm
$\left\vert \cdot\right\vert _{H^{s}\left(  \Omega\right)  }$ in $H^{s}\left(
\Omega\right)  $ which contains only the derivatives of order $s$ is a norm in
$H_{0}^{s}\left(  \Omega\right)  $.

We consider the Poisson problem in weak form:%
\begin{equation}
\text{Given }f\in L^{2}\left(  \Omega\right)  \text{ find }u\in H_{0}%
^{1}\left(  \Omega\right)  \text{\quad}a\left(  u,v\right)  :=\left(
\mathbb{A}\nabla u,\nabla v\right)  =\left(  f,v\right)  \text{\quad}\forall
v\in H_{0}^{1}\left(  \Omega\right)  . \label{varform}%
\end{equation}
Throughout the paper we assume that the diffusion matrix $\mathbb{A}\in
L^{\infty}\left(  \Omega,\mathbb{R}_{\operatorname*{sym}}^{d\times d}\right)
$ is symmetric and satisfies%
\begin{equation}
0<a_{\min}:=\underset{\mathbf{x}\in\Omega}{\operatorname*{ess}\inf}%
\inf_{\mathbf{v\in}\mathbb{R}^{d}\backslash\left\{  0\right\}  }\frac{\left(
\mathbb{A}\left(  \mathbf{x}\right)  \mathbf{v}\right)  \cdot\mathbf{v}%
}{\mathbf{v}\cdot\mathbf{v}}\,\leq\underset{\mathbf{x}\in\Omega}%
{\operatorname*{ess}\sup}\sup_{\mathbf{v}\in\mathbb{R}^{d}\backslash\left\{
0\right\}  }\frac{\left(  \mathbb{A}\left(  \mathbf{x}\right)  \mathbf{v}%
\right)  \cdot\mathbf{v}}{\mathbf{v}\cdot\mathbf{v}}=:a_{\max}<\infty
\label{aeps}%
\end{equation}
and{ that there exists a partition $\mathcal{P}:=(\Omega_{j})_{j=1}^{J}$ of
$\Omega$ into $J$ (possibly curved) polygons (polyhedra for }$d=3$) {such
that, for some appropriate }$r\in\mathbb{N}$, it holds%
\begin{equation}
\left\Vert \mathbb{A}\right\Vert _{{PW^{{r},\infty}\left(  \Omega\right)  }%
}:=\max_{1\leq j\leq J}\left\Vert \left.  \mathbb{A}\right\vert _{\Omega_{j}%
}\right\Vert _{W^{r,\infty}\left(  \Omega_{j}\right)  }<\infty. \label{aeps2}%
\end{equation}
Assumption (\ref{aeps}) implies the well-posedness of problem (\ref{varform})
via the Lax-Milgram lemma.

\section{Conforming hp-Finite Element Galerkin
Discretization\label{sectconfbasis}}

In this paper we restrict our studies to bounded, polygonal ($d=2$) or
polyhedral ($d=3$) Lipschitz domains $\Omega\subset\mathbb{R}^{d}$ and regular
finite element meshes $\mathcal{G}$ (in the sense of \cite{Ciarlet})
consisting of (closed) simplices $K$, where hanging nodes are not allowed. The
local and global mesh width is denoted by $h_{K}:=\operatorname*{diam}K$ and
$h:=\max_{K\in\mathcal{G}}h_{K}$. The boundary of a simplex $K$ can be split
into $\left(  d-1\right)  $-dimensional simplices (facets for $d=3$ and
triangle edges for $d=2$) which are denoted by $T$. The set of all facets in
$\mathcal{G}$ is called $\mathcal{F}$; the set of facets lying on
$\partial\Omega$ is denoted by $\mathcal{F}_{\partial\Omega}$ and defines a
triangulation of the surface $\partial\Omega$. The set of facets in $\Omega$
is denoted by $\mathcal{F}_{\Omega}$. As a convention we assume that simplices
and facets are closed sets. The interior of a simplex $K$ is denoted by
$\overset{\circ}{K}$ and we write $\overset{\circ}{T}$ to denote the
(relative) interior of a facet $T$. The set of all simplex vertices in the
mesh $\mathcal{G}$ is denoted by $\mathcal{V}$, those lying on $\partial
\Omega$ by $\mathcal{V}_{\partial\Omega}$, and those lying in $\Omega$ by
$\mathcal{V}_{\Omega}$. Similar the set of simplex edges in $\mathcal{G}$ is
denoted by $\mathcal{E}$, those lying on $\partial\Omega$ by $\mathcal{E}%
_{\partial\Omega}$, and those lying in $\Omega$ by $\mathcal{E}_{\Omega}$.

We recall the definition of conforming $hp$-finite element spaces (see, e.g.,
\cite{Schwab}). For $p\in\mathbb{N}_{0}:=\left\{  0,1,\ldots\right\}  $, let
$\mathbb{P}_{p}^{d}$ denote the space of $d$-variate polynomials of total
degree $\leq p$. For a connected subset $\omega\subset\Omega$, we write
$\mathbb{P}_{d}^{p}\left(  \omega\right)  $ for polynomials of degree $\leq p$
defined on $\omega$. For a connected $m$-dimensional manifold $\omega
\subset\mathbb{R}^{d}$, for which there exists a subset $\hat{\omega}%
\in\mathbb{R}^{m}$ along an affine bijection $\chi_{\omega}:\hat{\omega
}\rightarrow\omega$, we set $\mathbb{P}_{p}^{m}\left(  \omega\right)
:=\left\{  v\circ\chi_{\omega}^{-1}:v\in\mathbb{P}_{p}^{m}\left(  \hat{\omega
}\right)  \right\}  $. If the dimension $m$ is clear from the context, we
write $\mathbb{P}_{p}\left(  \omega\right)  $ short for $\mathbb{P}_{p}%
^{m}\left(  \omega\right)  $.

The conforming $hp$-finite element space is given by%
\begin{equation}
S_{\mathcal{G},\operatorname*{c}}^{p}:=\left\{  u\in C^{0}\left(
\overline{\Omega}\right)  \mid\forall K\in\mathcal{G}\quad\left.  u\right\vert
_{K}\in\mathbb{P}_{p}\left(  K\right)  \right\}  \cap H_{0}^{1}\left(
\Omega\right)  . \label{hpfinele}%
\end{equation}
A Lagrange basis for $S_{\mathcal{G},\operatorname*{c}}^{p}$ can be defined as
follows. Let
\begin{equation}
\widehat{\mathcal{N}}^{p}:=\left\{  \frac{\mathbf{i}}{p}:\mathbf{i}%
\in\mathbb{N}_{0}^{d}\text{ with }i_{1}+\ldots+i_{d}\leq p\right\}
\label{defnodalpointsref}%
\end{equation}
denote the equispaced unisolvent set of nodal points on the $d$-dimensional
unit simplex
\begin{equation}
\widehat{K}:=\left\{  \mathbf{x}\in\mathbb{R}_{\geq0}^{d}\mid x_{1}%
+\ldots+x_{d}\leq1\right\}  . \label{defrefelement}%
\end{equation}
For a simplex $K\in\mathcal{G}$, let $\chi_{K}:\widehat{K}\rightarrow K$
denote an affine mapping. The set of nodal points is given by%
\begin{equation}%
\begin{array}
[c]{l}%
\mathcal{N}^{p}:=\left\{  \chi_{K}\left(  \mathbf{\hat{N}}\right)
\mid\mathbf{\hat{N}}\in\widehat{\mathcal{N}}^{p},K\in\mathcal{G}\right\}  ,\\
\mathcal{N}_{\Omega}^{p}:=\mathcal{N}^{p}\cap\Omega,\qquad\mathcal{N}%
_{\partial\Omega}^{p}:=\mathcal{N}^{p}\cap\partial\Omega.
\end{array}
\label{defNodalpoints}%
\end{equation}
The Lagrange basis for $S_{\mathcal{G},\operatorname*{c}}^{p}$ can be indexed
by the nodal points $\mathbf{N}\in\mathcal{N}_{\Omega}^{p}$ and is
characterized by
\begin{equation}
B_{p,\mathbf{N}}^{\mathcal{G}}\in S_{\mathcal{G},\operatorname*{c}}^{p}%
\quad\text{and\quad}\forall\mathbf{N}^{\prime}\in\mathcal{N}_{\Omega}^{p}\quad
B_{p,\mathbf{N}}^{\mathcal{G}}\left(  \mathbf{N}^{\prime}\right)
=\delta_{\mathbf{N},\mathbf{N}^{\prime}}, \label{basisfunctions}%
\end{equation}
where $\delta_{\mathbf{N},\mathbf{N}^{\prime}}$ is the Kronecker delta.

\begin{definition}
\label{Deflocbasis}For all $K\in\mathcal{G}$, $T\in\mathcal{F}_{\Omega}$,
$E\in\mathcal{E}_{\Omega}$, $\mathbf{V}\in\mathcal{V}_{\Omega}$, the
conforming spaces $S_{K,\operatorname*{c}}^{p}$, $S_{T,\operatorname*{c}}^{p}%
$, $S_{E,\operatorname*{c}}^{p}$, $S_{\mathbf{V},\operatorname*{c}}^{p}$ are
given as the spans of the following basis functions%
\begin{align*}
S_{K,\operatorname*{c}}^{p}  &  :=\operatorname*{span}\left\{  B_{p,\mathbf{N}%
}^{\mathcal{G}}\mid\mathbf{N}\in\overset{\circ}{K}\cap\mathcal{N}_{\Omega}%
^{p}\right\}  ,\\
S_{T,\operatorname*{c}}^{p}  &  :=\operatorname*{span}\left\{  B_{p,\mathbf{N}%
}^{\mathcal{G}}\mid\mathbf{N}\in\overset{\circ}{T}\cap\mathcal{N}_{\Omega}%
^{p}\right\}  ,\\
S_{E,\operatorname*{c}}^{p}  &  :=\operatorname*{span}\left\{  B_{p,\mathbf{N}%
}^{\mathcal{G}}\mid\mathbf{N}\in\overset{\circ}{E}\cap\mathcal{N}_{\Omega}%
^{p}\right\}  ,\\
S_{\mathbf{V},\operatorname*{c}}^{p}  &  :=\operatorname*{span}\left\{
B_{p,\mathbf{V}}^{\mathcal{G}}\right\}  \text{.}%
\end{align*}

\end{definition}

The following proposition shows that these spaces give rise to a direct sum
decomposition and that these spaces are locally defined. To be more specific
we first have to introduce some notation.

For any facet $T\in\mathcal{F}_{\Omega}$, vertex $\mathbf{V}\in\mathcal{V}%
_{\Omega}$, and $E\in\mathcal{E}_{\Omega}$ we define the sets%
\begin{equation}%
\begin{array}
[c]{ll}%
\mathcal{G}_{T}:=\left\{  K\in\mathcal{G}:T\subset\partial K\right\}  , &
\omega_{T}:=%
{\displaystyle\bigcup\limits_{K\in\mathcal{G}_{T}}}
K,\\
\mathcal{G}_{\mathbf{V}}:=\left\{  K\in\mathcal{G}:\mathbf{V}\in\partial
K\right\}  , & \omega_{\mathbf{V}}:=%
{\displaystyle\bigcup\limits_{K\in\mathcal{G}_{\mathbf{V}}}}
K,\\
\mathcal{G}_{E}:=\left\{  K\in\mathcal{G}:E\subset\partial K\right\}  , &
\omega_{E}:=%
{\displaystyle\bigcup\limits_{K\in\mathcal{G}_{E}}}
K.
\end{array}
\label{defoftrianglesubsets}%
\end{equation}

\begin{proposition}
\label{PropSpaceDecomp}Let $S_{K,\operatorname*{c}}^{p}$,
$S_{T,\operatorname*{c}}^{p}$, $S_{E,\operatorname*{c}}^{p}$, $S_{\mathbf{V}%
,\operatorname*{c}}^{p}$ be as in Definition \ref{Deflocbasis}. Then the
direct sum decomposition holds%
\begin{equation}
S_{\mathcal{G},\operatorname*{c}}^{p}=\left(
{\displaystyle\bigoplus\limits_{\mathbf{V}\in\mathcal{V}_{\Omega}}}
S_{\mathbf{V},\operatorname*{c}}^{p}\right)  \oplus\left(
{\displaystyle\bigoplus\limits_{E\in\mathcal{E}_{\Omega}}}
S_{E,\operatorname*{c}}^{p}\right)  \oplus\left(
{\displaystyle\bigoplus\limits_{T\in\mathcal{F}_{\Omega}}}
S_{T,\operatorname*{c}}^{p}\right)  \oplus\left(
{\displaystyle\bigoplus\limits_{K\in\mathcal{G}}}
S_{K,\operatorname*{c}}^{p}\right)  . \label{spacedecomposition}%
\end{equation}

\end{proposition}

\section{Galerkin Discretization with Non-Conforming Crouzeix-Raviart Finite
Elements\label{section_nonconforme}}

\subsection{Non-Conforming Finite Elements with Weak Compatibility Conditions}

In this section, we will characterize a class of non-conforming finite element
spaces implicitly by a weak compatibility condition across the facets. For
each facet $T\in\mathcal{F}$, we fix a unit vector $\mathbf{n}_{T}$ which is
orthogonal to $T$. The orientation for the inner facets is arbitrary but fixed
while the orientation for the boundary facets is such that $\mathbf{n}_{T}$
points toward the exterior of $\Omega$. Our non-conforming finite element
spaces will be a subspace of%
\[
C_{\mathcal{G}}^{0}\left(  \Omega\right)  :=\left\{  u\in L^{\infty}\left(
\Omega\right)  \mid\forall K\in\mathcal{G\quad}\left.  u\right\vert
_{\overset{\circ}{K}}\in C^{0}\left(  \overset{\circ}{K}\right)  \right\}
\]
and we consider the skeleton $%
{\displaystyle\bigcup\limits_{T\in\mathcal{F}}}
T$ as a set of measure zero.

For $K\in\mathcal{G}$, we define the restriction operator $\gamma
_{K}:C_{\mathcal{G}}^{0}\left(  \Omega\right)  \rightarrow C^{0}\left(
K\right)  $ by%
\[
\left(  \gamma_{K}w\right)  \left(  \mathbf{x}\right)  =w\left(
\mathbf{x}\right)  \quad\forall\mathbf{x}\in\overset{\circ}{K}%
\]
and on the boundary $\partial K$ by continuous extension. For the inner facets
$T\in\mathcal{F}$, let $K_{T}^{1},K_{T}^{2}$ be the two simplices which share
$T$ as a common facet with the convention that $\mathbf{n}_{T}$ points into
$K_{2}$. We set $\omega_{T}:=K_{T}^{1}\cup K_{T}^{2}$. The jump $\left[
\cdot\right]  _{T}:C_{\mathcal{G}}^{0}\left(  \Omega\right)  \rightarrow
C^{0}\left(  T\right)  $ across $T$ is defined by%
\begin{equation}
\left[  w\right]  _{T}=\left.  \left(  \gamma_{K_{2}}w\right)  \right\vert
_{T}-\left.  \left(  \gamma_{K_{1}}w\right)  \right\vert _{T}.
\label{defjumps}%
\end{equation}
For vector-valued functions, the jump is defined component-wise. The
definition of the non-conforming finite elements involves orthogonal
polynomials on triangles which we introduce first.

Let $\widehat{T}$ denote the (closed) unit simplex in $\mathbb{R}^{d-1}$, with
vertices $\mathbf{0}$, $\left(  1,0,\ldots,0\right)  ^{\intercal}$, $\left(
0,1,0,\ldots,0\right)  ^{\intercal}$, $\left(  0,\ldots,0,1\right)
^{\intercal}$. For $n\in\mathbb{N}_{0}$, the set of orthogonal polynomials on
$\widehat{T}$ is given by%
\begin{equation}
\mathbb{P}_{n,n-1}^{\perp}\left(  \widehat{T}\right)  :=\left\{
\begin{array}
[c]{ll}%
\mathbb{P}_{0}\left(  \widehat{T}\right)  & n=0,\\
\left\{  u\in\mathbb{P}_{n}\left(  \widehat{T}\right)  \mid\int_{\widehat{T}%
}uv=0\quad\forall v\in\mathbb{P}_{n-1}\left(  \widehat{T}\right)  \right\}  &
n\geq1.
\end{array}
\right.  \label{defPnn-1senk}%
\end{equation}
We lift this space to a facet $T\in\mathcal{F}$ by employing an affine
transform $\chi_{T}:\widehat{T}\rightarrow T$
\[
\mathbb{P}_{n,n-1}^{\perp}\left(  T\right)  :=\left\{  v\circ\chi_{T}%
^{-1}:v\in\mathbb{P}_{n,n-1}^{\perp}\left(  T\right)  \right\}  .
\]
The orthogonal polynomials on triangles allows us to formulate the
\textit{weak compatibility condition} which is employed for the definition of
non-conforming finite element spaces:%
\begin{equation}
\left[  u\right]  _{T}\in\mathbb{P}_{p,p-1}^{\perp}\left(  T\right)
,\ \forall T\in\mathcal{F}_{\Omega}\quad\text{and\quad}\left.  u\right\vert
_{T}\in\mathbb{P}_{p,p-1}^{\perp}\left(  T\right)  ,\ \forall T\in
\mathcal{F}_{{\partial\Omega}}. \label{poinprep}%
\end{equation}
We have collected all ingredients for the (implicit) characterization of the
non-conforming Crouzeix-Raviart finite element space.

\begin{definition}
\label{DefinitionE}The non-conforming finite element space $S_{\mathcal{G}%
}^{p}$ with weak compatibility conditions across facets is given by%
\begin{equation}
S_{\mathcal{G}}^{p}:=\left\{  u\in L^{\infty}\left(  \Omega\right)
\mid\forall K\in\mathcal{G}\quad\gamma_{K}u\in\mathbb{P}_{p}\left(  K\right)
\text{ and }u\text{ satisfies\ (\ref{poinprep})}\right\}  . \label{DefSGp}%
\end{equation}

\end{definition}

The non-conforming Galerkin discretization of (\ref{varform}) for a given
finite element space $S$ which satisfies $S_{\mathcal{G},\operatorname*{nc}%
}^{p}\subset S\subset S_{\mathcal{G}}^{p}$ reads:%
\begin{equation}
\text{Given }f\in L^{2}\left(  \Omega\right)  \text{ find }u_{S}\in
S\text{\quad}a_{\mathcal{G}}\left(  u_{S},v\right)  :=\left(  \mathbb{A}%
\nabla_{\mathcal{G}}u_{S},\nabla_{\mathcal{G}}v\right)  =\left(  f,v\right)
\mathbf{\qquad\forall}v\in S \label{ncGalDisc}%
\end{equation}
where
\[
\nabla_{\mathcal{G}}u\left(  \mathbf{x}\right)  :=\nabla u\left(
\mathbf{x}\right)  \qquad\forall\mathbf{x}\in\Omega\backslash\left(
{\displaystyle\bigcup\limits_{T\in\mathcal{F}}}
\partial T\right)  .
\]

\subsection{Non-Conforming Finite Elements of Crouzeix-Raviart Type in
3D\label{SecNCFECR}}

The definition of the non-conforming space $S_{\mathcal{G}}^{p}$ in
(\ref{DefSGp}) is implicit via the weak compatibility condition. In this
section, we will present explicit representations of non-conforming basis
functions of Crouzeix-Raviart type for general polynomial order $p$. These
functions together with the conforming basis functions span a space
$S_{\mathcal{G},\operatorname*{nc}}^{p}$ which satisfies the inclusions
$S_{\mathcal{G},\operatorname*{c}}^{p}\subsetneq S_{\mathcal{G}%
,\operatorname*{nc}}^{p}\subseteq S_{\mathcal{G}}^{p}$ (cf. Theorem
\ref{thm-convergence-piecewise}). The derivation of the formula and their
algebraic properties will be the topic of the following sections.

We will introduce two types of non-conforming basis functions: those whose
support is one tetrahedron and those whose support consists of two adjacent
tetrahedrons, that is tetrahedrons which have a common facet. For details and
their derivation we refer to Section \ref{SecExplConstr} while here we focus
on the representation formulae.

\subsubsection{Non-Conforming Basis Functions Supported on One Tetrahedron}

The construction starts by defining \textit{symmetric orthogonal polynomials}
$b_{p,k}^{\operatorname*{sym}}$, $0\leq k\leq d_{\operatorname*{triv}}\left(
p\right)  -1$ on the reference triangle $\widehat{T}$ with vertices $\left(
0,0\right)  ^{\intercal}$, $\left(  1,0\right)  ^{\intercal}$, $\left(
0,1\right)  ^{\intercal}$, where%
\begin{equation}
d_{\operatorname*{triv}}\left(  p\right)  :=\left\lfloor \frac{p}%
{2}\right\rfloor -\left\lfloor \frac{p-1}{3}\right\rfloor . \label{dtrivfirst}%
\end{equation}
We define the coefficients%
\[
M_{i,j}^{\left(  p\right)  }=\left(  -1\right)  ^{p}~_{4}F_{3}\left(
\genfrac{}{}{0pt}{}{-j,j+1,-i,i+1}{-p,p+2,1}%
;1\right)  \frac{2i+1}{p+1}\qquad0\leq i,j\leq p,
\]
where $_{p}F_{q}$ denotes the generalized hypergeometric function (cf.
\cite[Chap. 16]{NIST:DLMF}). The $_{4}F_{3}$-sum is understood to terminate at
$i$ to avoid the $0/0$ ambiguities in the formal $_{4}F_{3}$-series. These
coefficients allow to define the polynomials%
\[
r_{p,2k}\left(  x_{1},x_{2}\right)  :=2\sum\limits_{0\leq j\leq p/2}%
M_{2j,2k}^{\left(  n\right)  }b_{p,2j}+b_{p,2k}\qquad0\leq k\leq p/2,
\]
where $b_{p,k}$, $0\leq k\leq p$, are the basis for the orthogonal polynomials
of degree $p$ on $\widehat{T}$ as defined afterwards in (\ref{defbnk}). Then,
a basis for the symmetric orthogonal polynomials is given by%
\begin{equation}
b_{p,k}^{\operatorname*{sym}}:=\left\{
\begin{array}
[c]{ll}%
r_{p,p-2k} & \text{if }p\text{ is even,}\\
r_{p,p-1-2k} & \text{if }p\text{ is odd,}%
\end{array}
\right.  \qquad k=0,1,\ldots,d_{\operatorname*{triv}}\left(  p\right)  -1.
\label{defbpksym}%
\end{equation}

The non-conforming Crouzeix-Raviart basis function $B_{p,k}^{\widehat
{K},\operatorname*{nc}}\in\mathbb{P}_{p}\left(  \widehat{K}\right)  $ on the
unit tetrahedron $\widehat{K}$ is characterized\ by its values at the nodal
points in $\widehat{\mathcal{N}}^{p}$ (cf. (\ref{defnodalpointsref})). For a
facet $T\subset\partial\widehat{K}$, let $\chi_{T}:\widehat{T}\rightarrow T$
denote an affine pullback to the reference triangle. Then $B_{p,k}%
^{\widehat{K},\operatorname*{nc}}\in\mathbb{P}_{p}\left(  \widehat{K}\right)
$ is uniquely defined by%
\begin{equation}
B_{p,k}^{\widehat{K},\operatorname*{nc}}\left(  \mathbf{N}\right)  :=\left\{
\begin{array}
[c]{ll}%
b_{p,k}^{\operatorname*{sym}}\circ\chi_{T}^{-1}\left(  \mathbf{N}\right)  &
\forall\mathbf{N}\in\widehat{\mathcal{N}}^{p}\text{ s.t. }\mathbf{N}\in
T\text{ for some facet }T\subset\partial\widehat{K},\\
0 & \forall\mathbf{N}\in\widehat{\mathcal{N}}^{p}\backslash\partial\widehat{K}%
\end{array}
\right.  \qquad k=0,1,\ldots,d_{\operatorname*{triv}}\left(  p\right)  -1.
\label{DefBpkKhut}%
\end{equation}

\begin{remark}
\label{RemDef}In Sec. \ref{SecFullSymPoly}, we will prove that the polynomials
$b_{p,k}^{\operatorname*{sym}}$ are totally symmetric, i.e., invariant under
affine bijections $\chi:\widehat{K}\rightarrow\widehat{K}$. Thus, any of these
functions can be lifted to the facets of a tetrahedron via affine pullbacks
and the resulting function on the surface is continuous. As a consequence, the
value $B_{p,k}^{\widehat{K},\operatorname*{nc}}\left(  \mathbf{N}\right)  $ in
definition (\ref{DefBpkKhut}) is independent of the choice of $T$ also for
nodal points $\mathbf{N}$ which belong to different facets.

It will turn out that the value $0$ at the inner nodes could be replaced by
other values without changing the arising non-conforming space. Other choices
could be preferable in the context of inverse inequalities and the condition
number of the stiffness matrix. However, we recommend to choose these values
such that the symmetries of $B_{p,k}^{\widehat{K},\operatorname*{nc}}$ are preserved.
\end{remark}

\begin{definition}
\label{DefSymSpaceK}The non-conforming tetrahedron-supported basis functions
on the reference element are given by%
\begin{equation}
B_{p,k}^{\widehat{K},\operatorname*{nc}}=\sum_{\mathbf{N}\in\widehat
{\mathcal{N}}^{p}\cap\partial\widehat{K}}B_{p,k}^{\widehat{K}%
,\operatorname*{nc}}\left(  \mathbf{N}\right)  B_{p,\mathbf{N}}^{\mathcal{G}%
}\qquad k=0,1,\ldots,d_{\operatorname*{triv}}\left(  p\right)  -1 \label{18b}%
\end{equation}
with values $B_{p,k}^{\widehat{K},\operatorname*{nc}}\left(  \mathbf{N}%
\right)  $ as in (\ref{DefBpkKhut}). For a simplex $K\in\mathcal{G}$ the
corresponding non-conforming basis functions $B_{p,k}^{K,\operatorname*{nc}}$
are given by lifting $B_{p,k}^{\widehat{K},\operatorname*{nc}}$ via an affine
pullback $\chi_{K}$ from $\widehat{K}$ to $K\in\mathcal{G}$:%
\[
\left.  B_{p,k}^{K,\operatorname*{nc}}\right\vert _{\overset{\circ}{K^{\prime
}}}:=\left\{
\begin{array}
[c]{ll}%
B_{p,k}^{\widehat{K},\operatorname*{nc}}\circ\chi_{K}^{-1} & K=K^{\prime},\\
0 & K\neq K^{\prime}.
\end{array}
\right.
\]
and span the space%
\begin{equation}
S_{K,\operatorname*{nc}}^{p}:=\operatorname*{span}\left\{  B_{p,k}%
^{K,\operatorname*{nc}}:k=0,1,\ldots,d_{\operatorname*{triv}}\left(  p\right)
-1\right\}  . \label{DefSkpsym}%
\end{equation}

\end{definition}

\begin{example}
The lowest order of $p$ such that $d_{\operatorname*{triv}}\left(  p\right)
\geq1$ is $p=2$. In this case, we get $d_{\operatorname*{triv}}\left(
p\right)  =1$. In Figure \ref{FigSym} the function $b_{p,k}%
^{\operatorname*{sym}}$ and corresponding basis functions $B_{p,k}%
^{K,\operatorname*{nc}}$ are depicted for $\left(  p,k\right)  \in\left\{
\left(  2,0\right)  ,\left(  3,0\right)  ,\left(  6,0\right)  ,\left(
6,1\right)  \right\}  $.%
\begin{table}[tbp] \centering
\begin{tabular}
[c]{|l|l|l|l|}\hline%
{\includegraphics[
height=1.1779in,
width=1.4503in
]%
{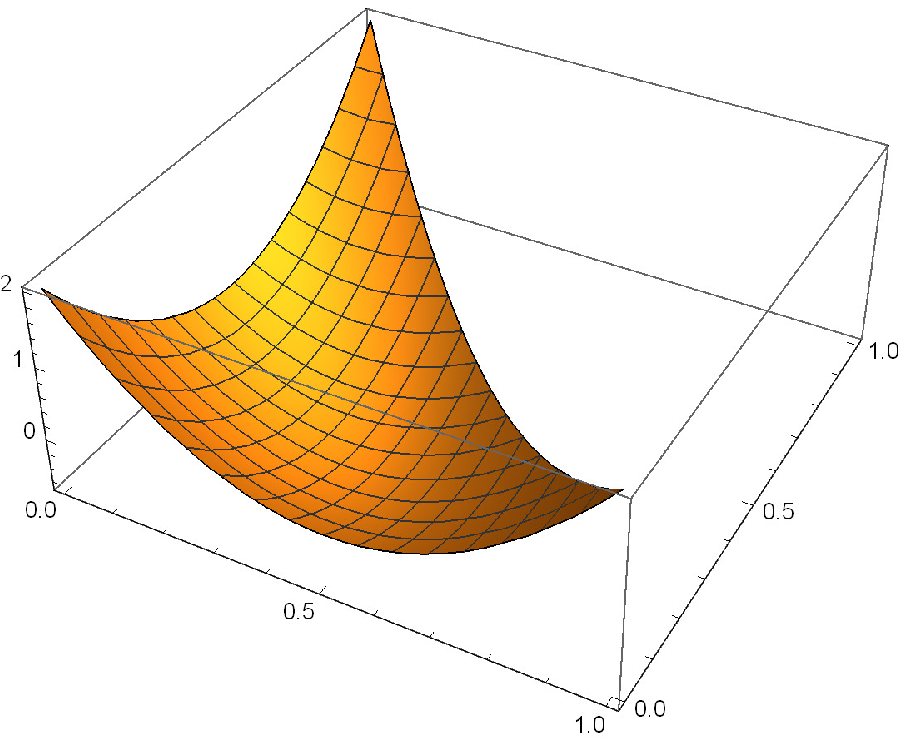}%
}%
&
{\includegraphics[
height=1.177in,
width=1.4434in
]%
{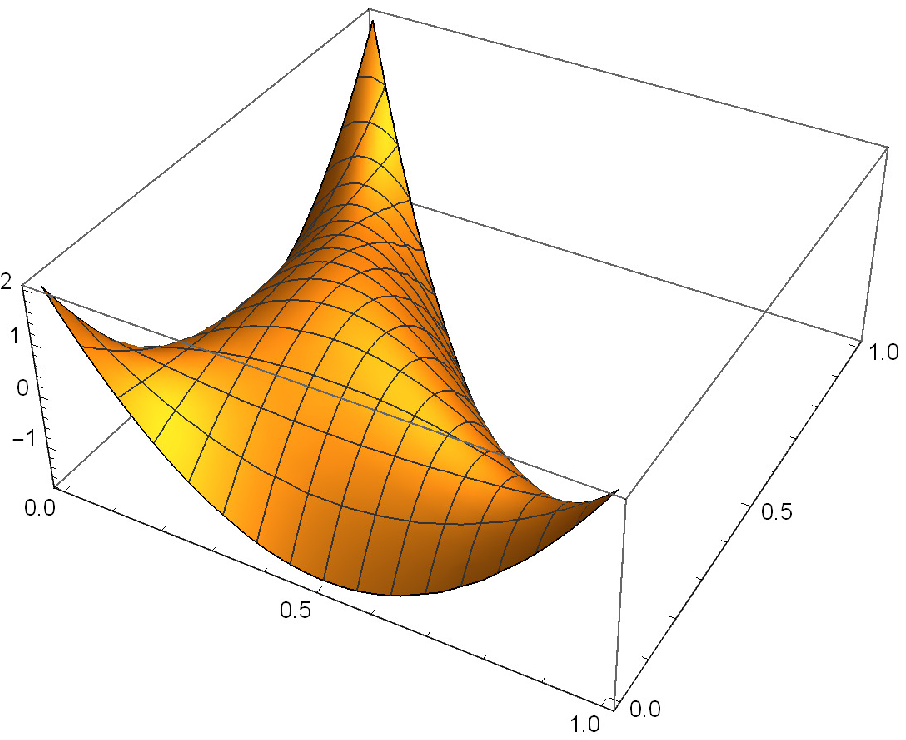}%
}%
&
{\includegraphics[
height=1.1718in,
width=1.4434in
]%
{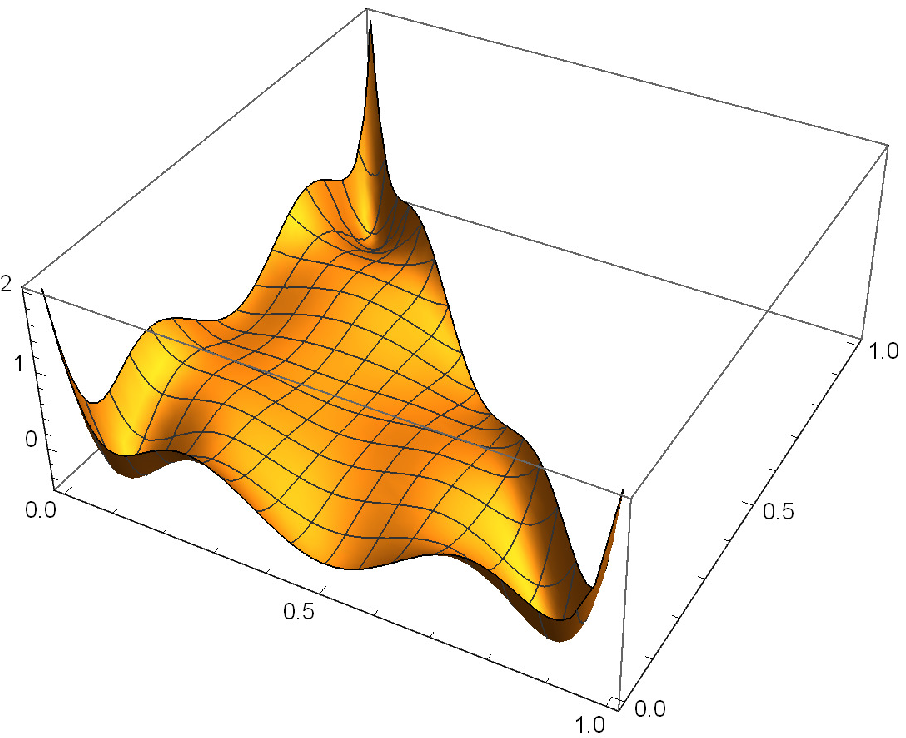}%
}%
&
{\includegraphics[
height=1.1666in,
width=1.4364in
]%
{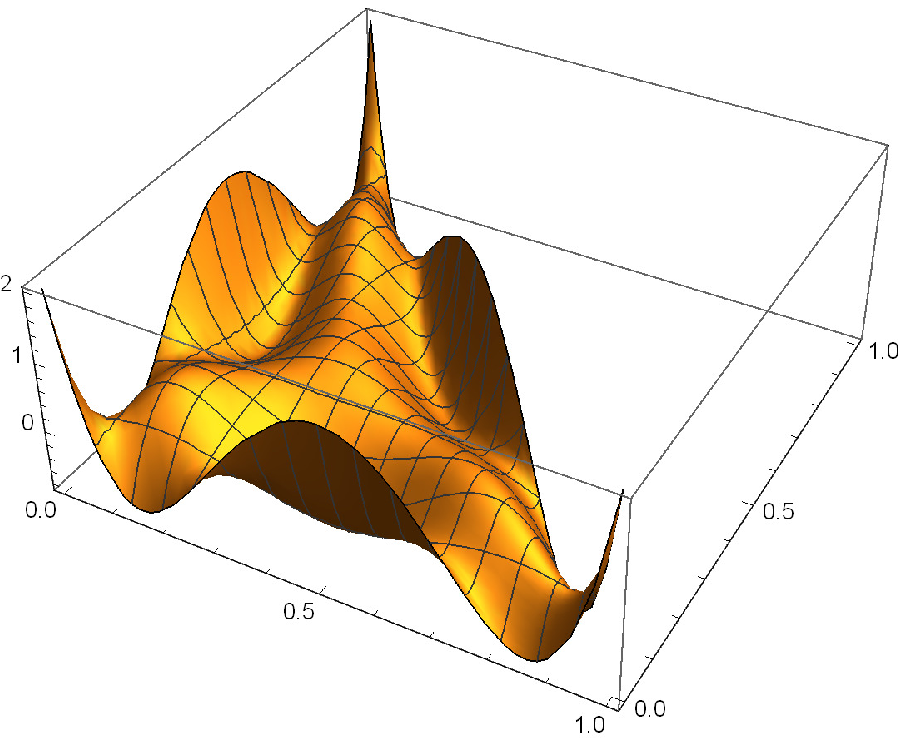}%
}%
\\%
{\includegraphics[
height=1.1018in,
width=1.4676in
]%
{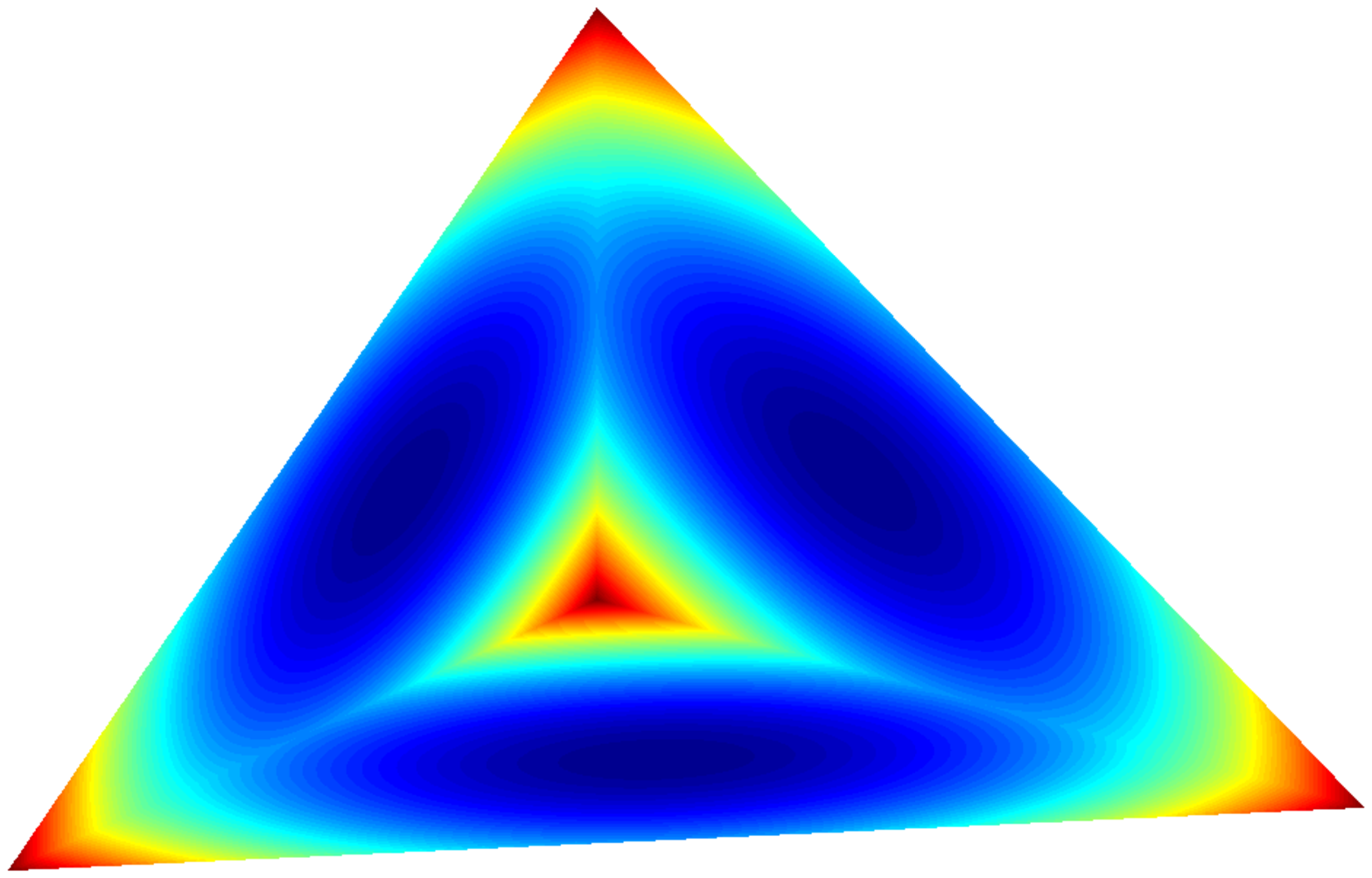}%
}%
&
{\includegraphics[
height=1.0914in,
width=1.4538in
]%
{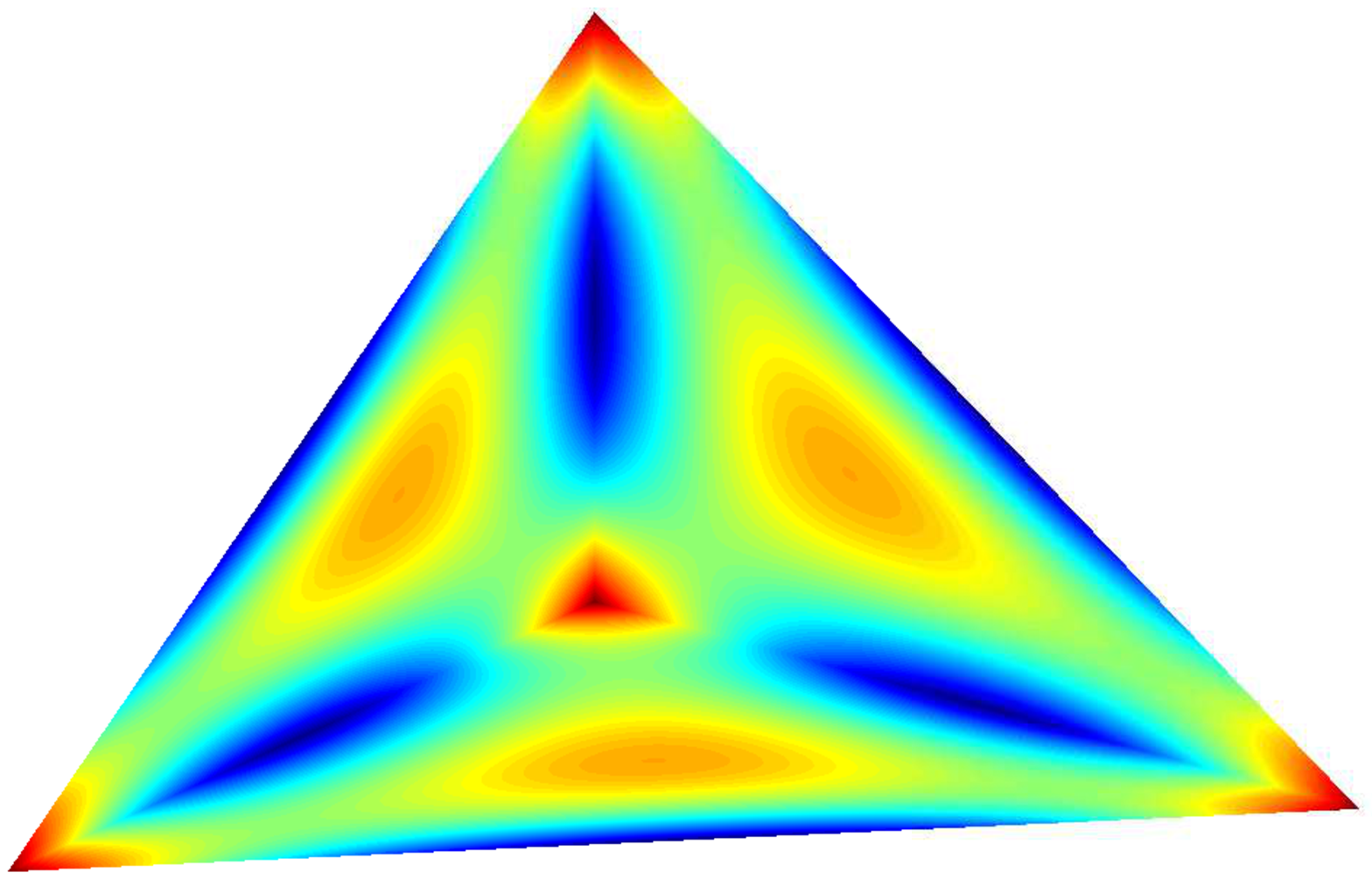}%
}%
&
{\includegraphics[
height=1.0862in,
width=1.4468in
]%
{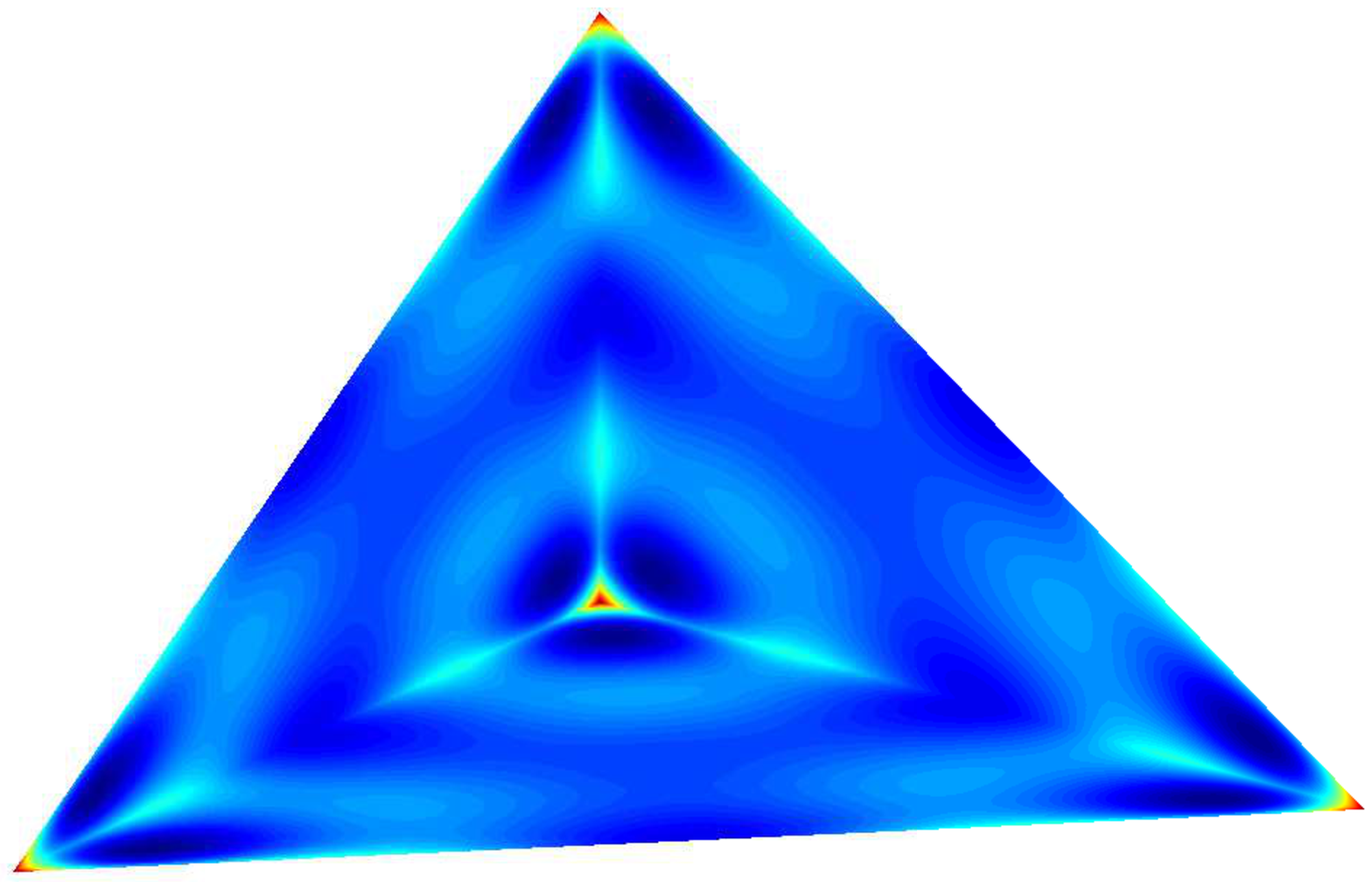}%
}%
&
{\includegraphics[
height=1.081in,
width=1.4399in
]%
{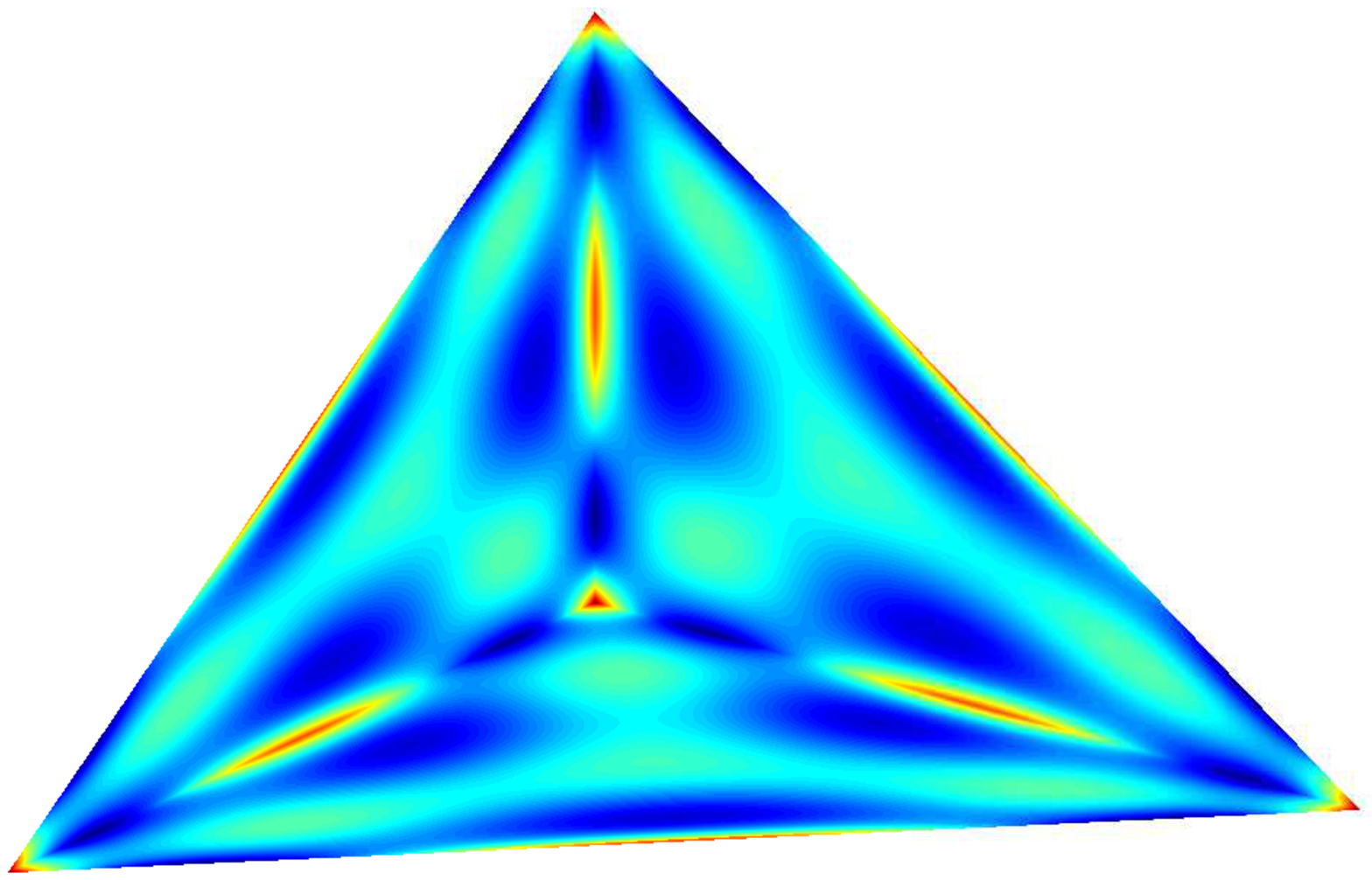}%
}%
\\
\multicolumn{1}{|c|}{$b_{2,0}^{\operatorname*{sym}},B_{2,0}%
^{K,\operatorname*{nc}}$} & \multicolumn{1}{|c|}{$b_{3,0}^{\operatorname*{sym}%
},B_{3,0}^{K,\operatorname*{nc}}$} & \multicolumn{1}{|c|}{$b_{6,0}%
^{\operatorname*{sym}},B_{6,0}^{K,\operatorname*{nc}}$} &
\multicolumn{1}{|c|}{$b_{6,1}^{\operatorname*{sym}},B_{6,1}%
^{K,\operatorname*{nc}}$}\\\hline
\end{tabular}
\caption{Symmetric orthogonal polynomials on the reference triangle and corresponding tetrahedron-supported non-conforming basis functions.\label{FigSym}}%
\end{table}%

\end{example}

\subsubsection{Non-Conforming Basis Functions Supported on Two Adjacent
Tetrahedrons}

The starting point is to define orthogonal polynomials $b_{p,k}%
^{\operatorname*{refl}}$ on the reference triangle $\widehat{T}$ which are
mirror symmetric\footnote{The superscript \textquotedblleft%
$\operatorname*{refl}$\textquotedblright\ is a shorthand for \textquotedblleft
reflection\textquotedblright\ and explained in Section \ref{SubSecIrr}.} with
respect to the angular bisector in $\widehat{T}$ through $\mathbf{0}$ and
linear independent from the fully symmetric functions $b_{p,k}%
^{\operatorname*{sym}}$. We set%
\begin{equation}
b_{p,k}^{\operatorname*{refl}}:=\frac{1}{3}\left(  2b_{p,2k}\left(
x_{1},x_{2}\right)  -b_{p,2k}\left(  x_{2},1-x_{1}-x_{2}\right)
-b_{p,2k}\left(  1-x_{1}-x_{2},x_{1}\right)  \right)  \qquad0\leq k\leq
d_{\operatorname*{refl}}\left(  p\right)  -1, \label{brefl1stdef}%
\end{equation}
where%
\begin{equation}
d_{\operatorname*{refl}}\left(  p\right)  :=\left\lfloor \frac{p+2}%
{3}\right\rfloor . \label{defdreflp}%
\end{equation}

Let $K_{1}$, $K_{2}$ denote two tetrahedrons which share a common facet, say
$T$. The vertex of $K_{i}$ which is opposite to $T$ is denoted by
$\mathbf{V}_{i}$. The procedure of lifting the nodal values to the facets of
$\omega_{T}:=K_{1}\cup K_{2}$ is analogous as for the basis functions
$B_{n,k}^{K,\operatorname*{nc}}$. However, it is necessary to choose the
pullback $\chi_{i,\tilde{T}}:\widehat{T}\rightarrow\tilde{T}$ of a facet
$\tilde{T}\subset\partial K_{i}\backslash\overset{\circ}{T}$ such that the
origin is mapped to $\mathbf{V}_{i}$.
\begin{equation}
B_{p,k}^{T,\operatorname*{nc}}\left(  \mathbf{N}\right)  :=\left\{
\begin{array}
[c]{ll}%
b_{p,k}^{\operatorname*{refl}}\circ\chi_{i,\tilde{T}}^{-1}\left(
\mathbf{N}\right)  & \forall\mathbf{N}\in\mathcal{N}^{p}\text{ s.t.
}\mathbf{N}\in\tilde{T}\text{ for some facet }\tilde{T}\subset\partial
K\backslash\overset{\circ}{T}_{i},\\
0 & \forall\mathbf{N}\in\mathcal{N}^{p}\cap\overset{\circ}{\omega_{T}}%
\end{array}
\right.  \qquad k=0,1,\ldots,d_{\operatorname*{refl}}\left(  p\right)  -1.
\label{defvaluesreflt}%
\end{equation}
Again, the value $0$ at the inner nodes of $\omega_{T}$ could be replaced by
other values without changing the arising non-conforming space.

\begin{definition}
\label{PropBpmTncbasis}The non-conforming facet-oriented basis functions are
given by%
\begin{equation}
B_{p,k}^{T,\operatorname*{nc}}=\sum_{\mathbf{N}\in\mathcal{N}^{p}\cap
\partial\omega_{T}}B_{p,k}^{T,\operatorname*{nc}}\left(  \mathbf{N}\right)
\left.  B_{p,\mathbf{N}}^{\mathcal{G}}\right\vert _{\omega_{T}}\qquad\forall
T\in\mathcal{F}_{\Omega},\quad k=0,1,\ldots,d_{\operatorname*{refl}}\left(
p\right)  -1 \label{defedgesupp}%
\end{equation}
with values $B_{p,k}^{T,\operatorname*{nc}}\left(  \mathbf{N}\right)  $ as in
(\ref{defvaluesreflt}) and span the space%
\begin{equation}
S_{T,\operatorname*{nc}}^{p}:=\operatorname*{span}\left\{  B_{p,k}%
^{T,\operatorname*{nc}}:k=0,1,\ldots,d_{\operatorname*{refl}}\left(  p\right)
-1\right\}  . \label{STpncdef}%
\end{equation}

The non-conforming finite element space of Crouzeix-Raviart type is given by%
\begin{equation}
S_{\mathcal{G},\operatorname*{nc}}^{p}:=\left(
{\displaystyle\bigoplus\limits_{E\in\mathcal{E}_{\Omega}}}
S_{E,\operatorname*{c}}^{p}\right)  \oplus\left(
{\displaystyle\bigoplus\limits_{T\in\mathcal{F}_{\Omega}}}
S_{T,\operatorname*{c}}^{p}\right)  \oplus\left(
{\displaystyle\bigoplus\limits_{K\in\mathcal{G}}}
S_{K,\operatorname*{c}}^{p}\right)  \oplus\left(
{\displaystyle\bigoplus\limits_{K\in\mathcal{G}}}
S_{K,\operatorname*{nc}}^{p}\right)  \oplus\left(
{\displaystyle\bigoplus\limits_{T\in\mathcal{F}_{\Omega}}}
\operatorname*{span}\left\{  B_{p,0}^{T,\operatorname*{nc}}\right\}  \right)
. \label{DefSCR}%
\end{equation}

\end{definition}

\begin{remark}
In Sec. \ref{Sectaureflcomp}, we will show that the polynomials $b_{p,k}%
^{\operatorname*{refl}}$ are mirror symmetric with respect to the angular
bisector in $\widehat{T}$ through $\mathbf{0}$. Thus, any of these functions
can be lifted to the outer facets of two adjacent tetrahedrons via (oriented)
affine pullbacks as employed in (\ref{defvaluesreflt}) and the resulting
function on the surface is continuous. As a consequence, the value
$B_{p,k}^{T,\operatorname*{nc}}\left(  \mathbf{N}\right)  $ in definition
(\ref{defvaluesreflt}) is independent of the choice of $T$ also for nodal
points $\mathbf{N}$ which belong to different facets.

In Theorem \ref{Theorem33}, we will prove that (\ref{DefSCR}), in fact, is a
direct sum and a basis is given by the functions%
\[
B_{p,\mathbf{N}}^{\mathcal{G}}\quad\forall\mathbf{N}\in\mathcal{N}_{\Omega
}\backslash\mathcal{V},\quad B_{p,k}^{K,\operatorname*{nc}}\quad\forall
K\in\mathcal{G},0\leq k\leq d_{\operatorname*{triv}}\left(  p\right)  -1,\quad
B_{p,0}^{T,\operatorname*{nc}}\quad\forall T\in\mathcal{F}_{\Omega}.
\]
Also we will prove that $S_{\mathcal{G},\operatorname*{c}}^{p}\subsetneq
S_{\mathcal{G},\operatorname*{nc}}^{p}\subseteq S_{\mathcal{G}}^{p}$. This
condition implies that the convergence estimates as in Theorem
\ref{thm-convergence-piecewise} are valid for this space. We restricted the
reflection-type non-conforming basis functions to the lowest order $k=0$ in
order to keep the functions linearly independent.
\end{remark}

\begin{example}
The lowest order of $p$ such that $d_{\operatorname*{refl}}\left(  p\right)
\geq1$ is $p=1$. In this case, we get $d_{\operatorname*{refl}}\left(
p\right)  =1$. In Figure \ref{FigRefl} the function $b_{p,k}%
^{\operatorname*{refl}}$ and corresponding basis functions $B_{p,k}%
^{T,\operatorname*{nc}}$ are depicted for $\left(  p,k\right)  \in\left\{
\left(  1,0\right)  ,\left(  2,0\right)  ,\left(  4,0\right)  ,\left(
4,1\right)  \right\}  $.%
\begin{table}[tbp] \centering
\begin{tabular}
[c]{|c|c|c|c|}\hline%
{\includegraphics[
height=1.222in,
width=1.5065in
]%
{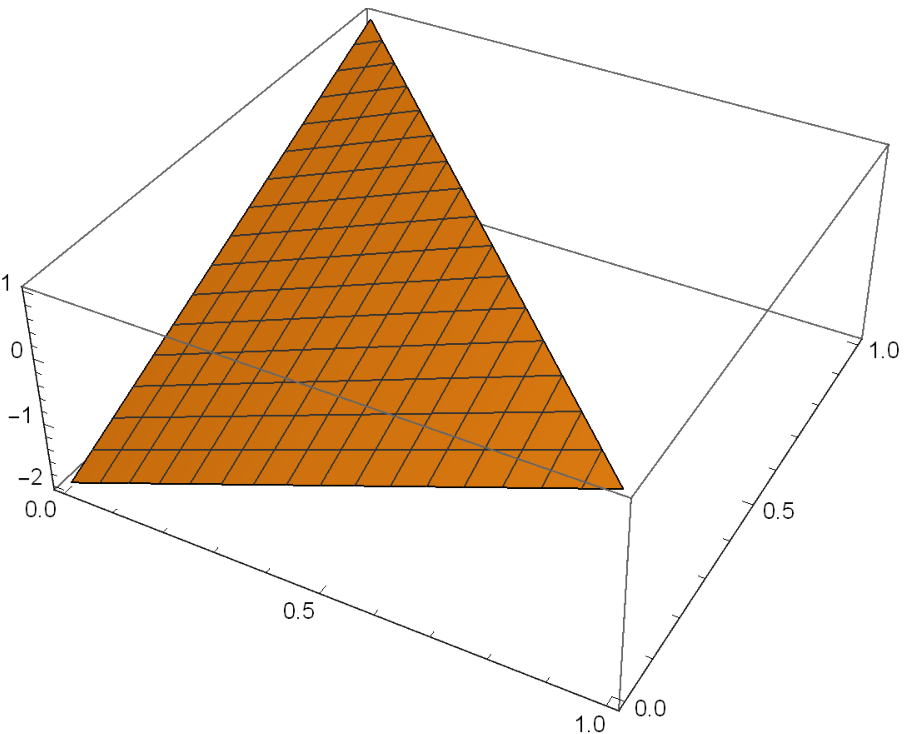}%
}%
&
{\includegraphics[
height=1.2168in,
width=1.5273in
]%
{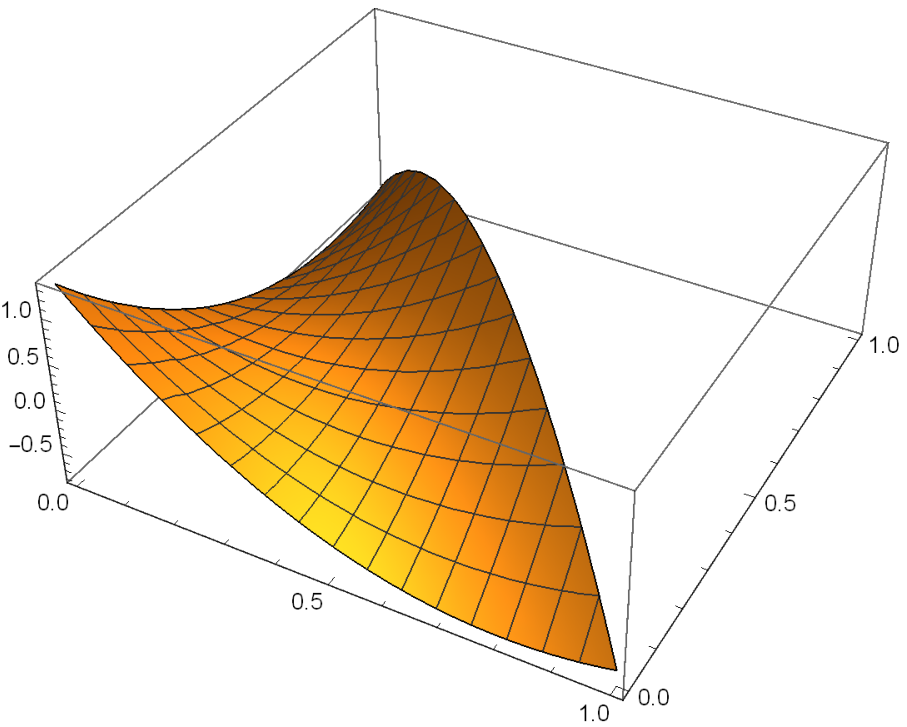}%
}%
&
{\includegraphics[
height=1.222in,
width=1.4987in
]%
{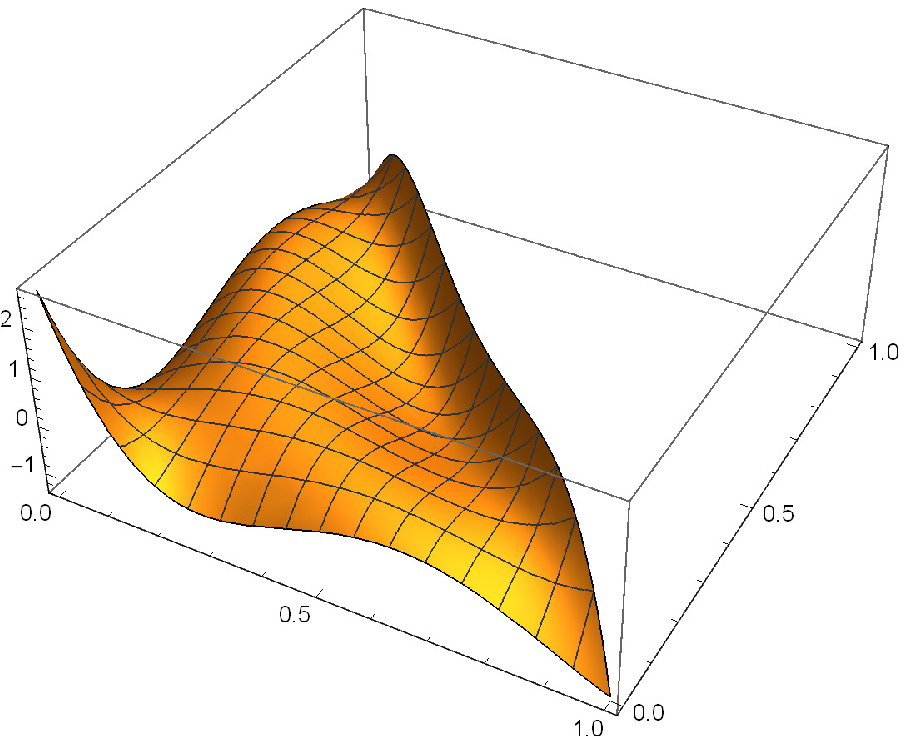}%
}%
&
{\includegraphics[
height=1.2168in,
width=1.5273in
]%
{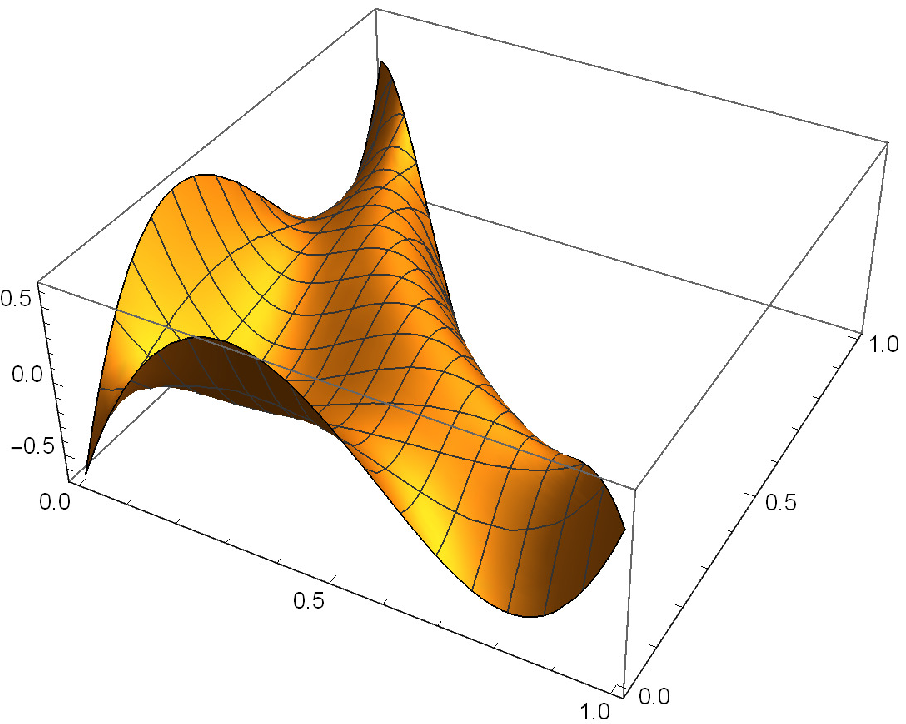}%
}%
\\%
{\includegraphics[
height=1.2721in,
width=1.51in
]%
{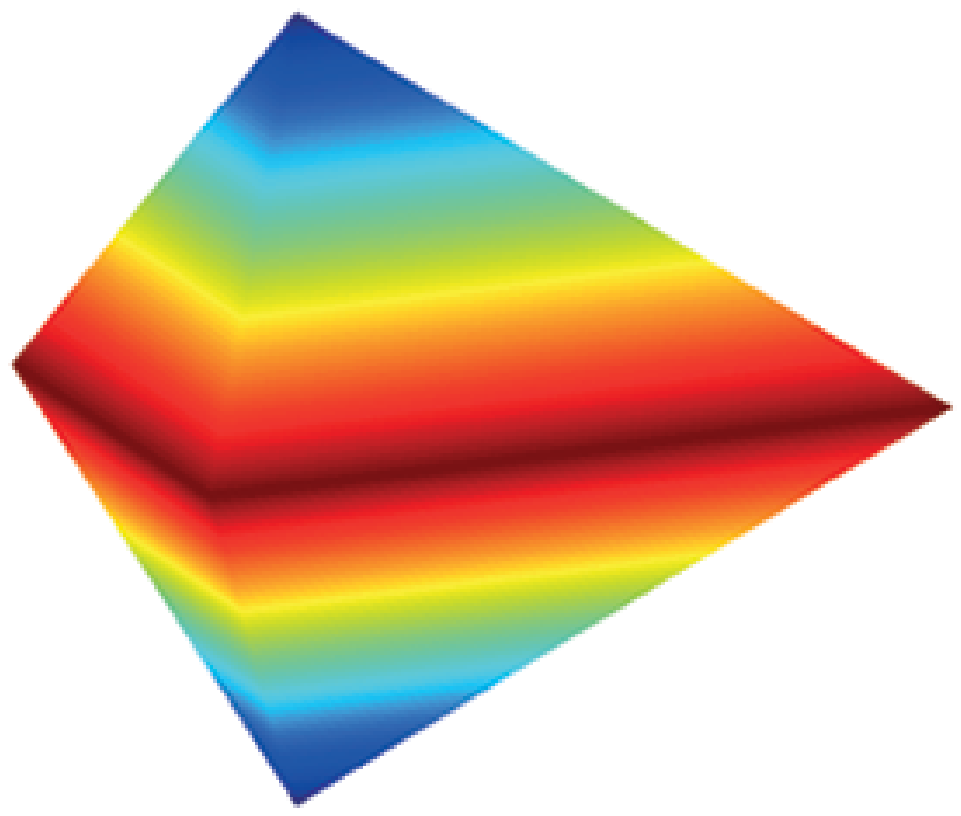}%
}%
&
{\includegraphics[
height=1.2695in,
width=1.5866in
]%
{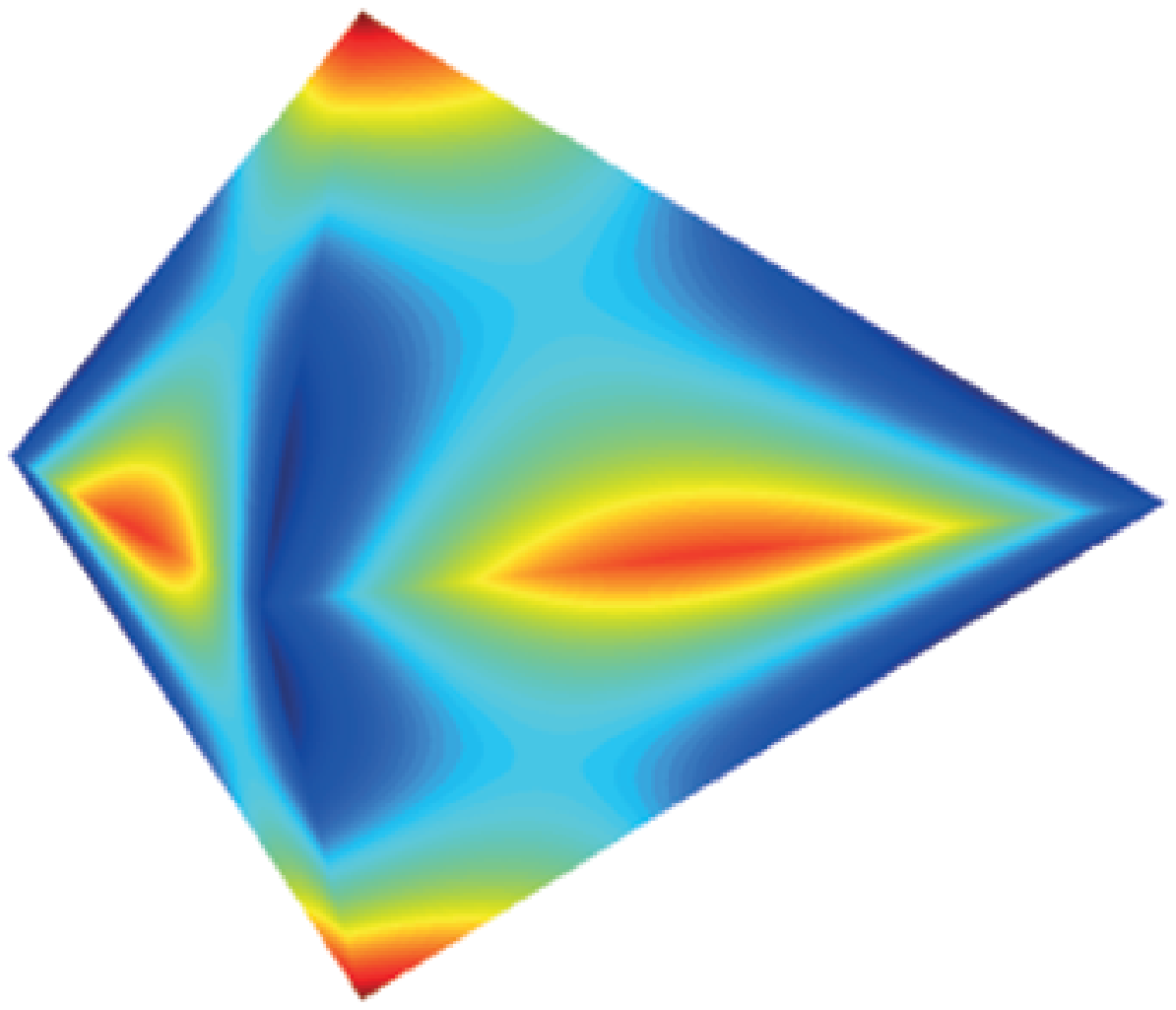}%
}%
&
{\includegraphics[
height=1.2721in,
width=1.5004in
]%
{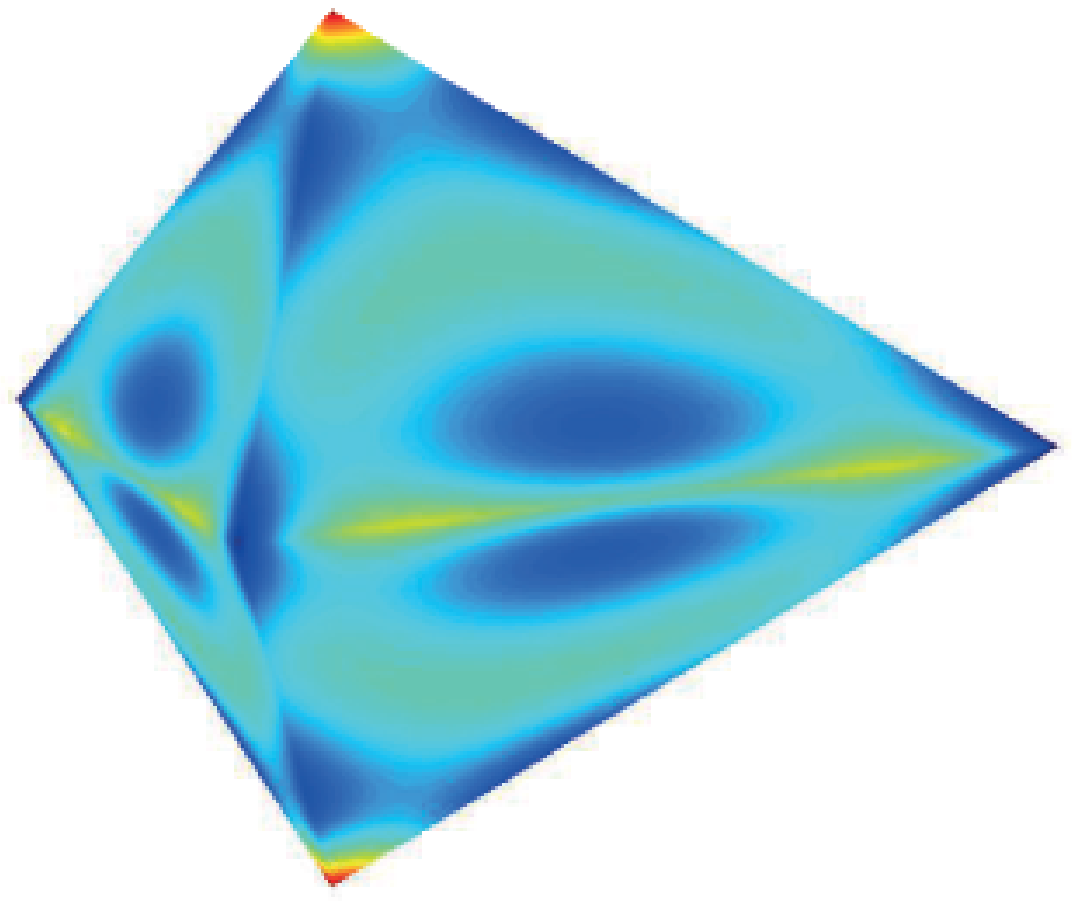}%
}%
&
{\includegraphics[
height=1.2687in,
width=1.4183in
]%
{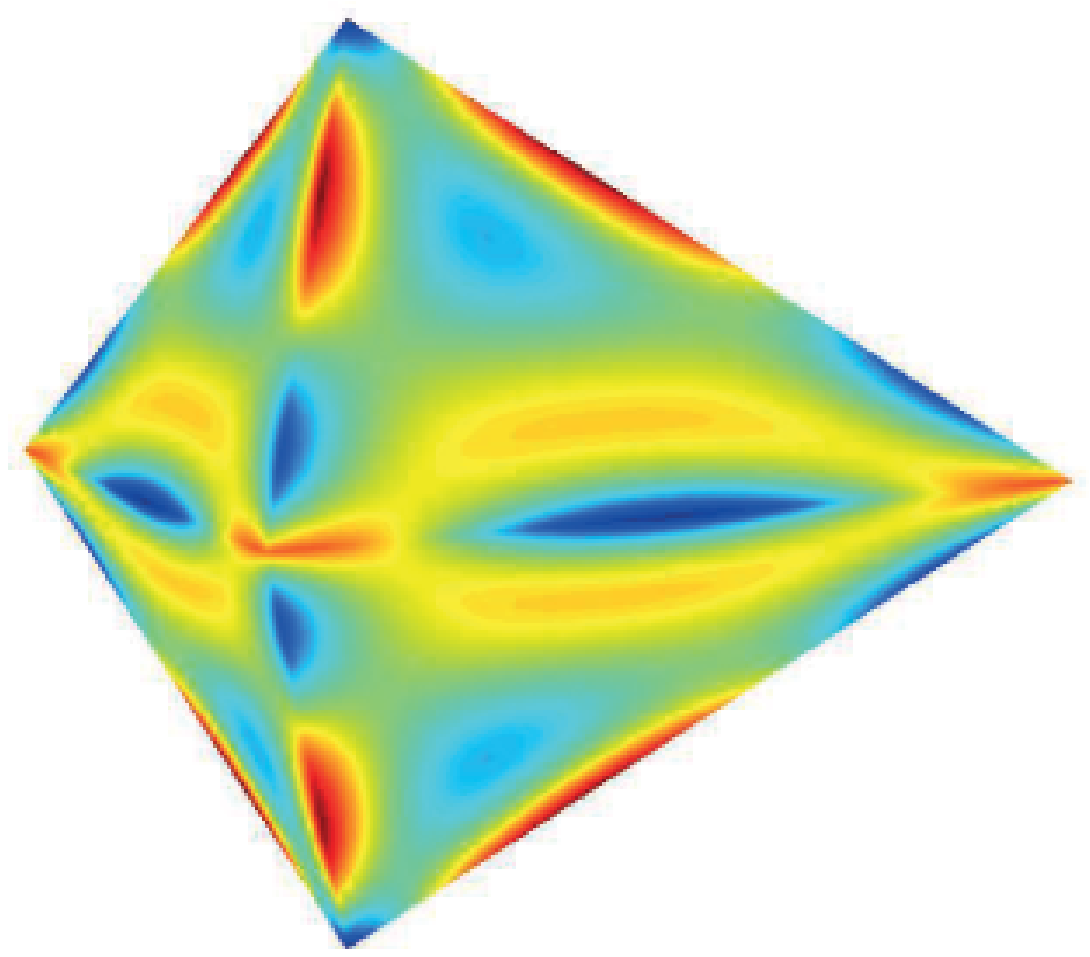}%
}%
\\
\multicolumn{1}{|c|}{$b_{1,0}^{\operatorname*{refl}},B_{1,0}%
^{T,\operatorname*{nc}}$} & \multicolumn{1}{|c|}{$b_{2,0}%
^{\operatorname*{refl}},B_{2,0}^{T,\operatorname*{nc}}$} &
\multicolumn{1}{|c|}{$b_{4,0}^{\operatorname*{refl}},B_{4,0}%
^{T,\operatorname*{nc}}$} & \multicolumn{1}{|c|}{$b_{4,1}%
^{\operatorname*{refl}},B_{4,1}^{T,\operatorname*{nc}}$}\\\hline
\end{tabular}
\caption{Orthogonal polynomials of reflection type and corresponding non-conforming basis functions which are supported on two adjacent tetrahedrons.
The common facet is horizontal and the two tetrahedrons are on top of each other.\label{FigRefl}}%
\end{table}%

\end{example}

\subsection{Error Analysis}

In this subsection we present the error analysis for the Galerkin
discretization (\ref{ncGalDisc}) with the non-conforming finite element space
$S_{\mathcal{G}}^{p}$ and subspaces thereof. The analysis is based on the
second Strang lemma and has been presented for an intrinsic version of
$S_{\mathcal{G}}^{p}$ in \cite{ccss_2012}.

For any inner facet $T\in\mathcal{F}$ and any $v\in S_{\mathcal{G}}^{p}$,
condition (\ref{poinprep}) implies $\int_{T}\left[  v\right]  _{T}=0$ : hence,
the jump $\left[  v\right]  _{T}$ is always zero-mean valued. Let $h_{T}$
denote the diameter of $T$. The combination of a Poincar\'{e} inequality with
a trace inequality then yields%
\begin{equation}
\left\Vert \left[  u\right]  _{T}\right\Vert _{L^{2}\left(  T\right)  }\leq
Ch_{T}\left\vert \left[  u\right]  _{T}\right\vert _{H^{1}\left(  T\right)
}\leq\tilde{C}h_{T}^{1/2}\left\vert u\right\vert _{H_{\operatorname*{pw}}%
^{1}\left(  \omega_{T}\right)  }, \label{poincarecondspa}%
\end{equation}
where%
\[
\left\vert u\right\vert _{H_{\operatorname*{pw}}^{p}\left(  \omega_{T}\right)
}:=\left(  \sum_{K\subset\omega_{T}}\left\vert u\right\vert _{H^{p}\left(
K\right)  }^{2}\right)  ^{1/2}\text{.}%
\]
In a similar fashion we obtain for all boundary facets $T\in\mathcal{F}%
_{\partial\Omega}$ and all $u\in S_{\mathcal{G}}^{p}$ the estimate%
\begin{equation}
\left\Vert u\right\Vert _{L^{2}\left(  T\right)  }\leq\tilde{C}h_{T}%
^{1/2}\left\vert u\right\vert _{H_{\operatorname*{pw}}^{1}\left(  \omega
_{T}\right)  }. \label{poincarecondspb}%
\end{equation}

{We say that the exact solution }${u}\in H${$_{0}^{1}\left(  \Omega\right)  $
is piecewise smooth over the partition $\mathcal{P=}\left(  \Omega_{j}\right)
_{j=1}^{J}$},{ if there exists some positive integer $s$ such that
\[
u_{|\Omega_{j}}\in H^{1+s}\left(  \Omega_{j}\right)  \quad\text{for
}j=1,2,\ldots,J.
\]
We write }$u\in${$PH^{1+s}(\Omega)$ and refer for further properties and
generalizations to non-integer values of }$s${, e.g., to \cite[Sec.
4.1.9]{SauterSchwab2010}.}

{For the approximation results, the finite element meshes $\mathcal{G}$ are
assumed to be \textit{compatible} with the partition $\mathcal{P}$ in the
following sense: for all $K\in\mathcal{G}$, there exists a single index $j$
such that }$\overset{\circ}{{K}}${$\cap\Omega_{j}\neq\emptyset$. }

The proof that $\left\vert \cdot\right\vert _{H_{\operatorname*{pw}}%
^{1}\left(  \Omega\right)  }$ is a norm on $S_{\mathcal{G}}^{p}$ is similar as
in \cite[Sect. 10.3]{scottbrenner3}: For $w\in H_{0}^{1}\left(  \Omega\right)
$ this follows from $\left\vert w\right\vert _{H_{\operatorname*{pw}}%
^{1}\left(  \Omega\right)  }=\left\Vert \nabla w\right\Vert $ and a Friedrichs
inequality; for $w\in S_{\mathcal{G}}^{p}$ the condition $\left\Vert
\nabla_{\mathcal{G}}w\right\Vert =0$ implies that $\left.  w\right\vert _{K}$
is constant on all simplices $K\in\mathcal{G}$. The combination with $\int
_{T}w=0$ for all $T\in\mathcal{F}_{\partial\Omega}$ leads to $\left.
w\right\vert _{K}=0$ for the outmost simplex layer via a Poincar\'{e}
inequality, i.e., $\left.  w\right\vert _{K}=0$ for all $K\in\mathcal{G}$
having at least one facet on $\partial\Omega$. This argument can be iterated
step by step over simplex layers towards the interior of $\Omega$ to finally
obtain $w=0$.

\begin{theorem}
\label{thm-convergence-piecewise}Let $\Omega\subset\mathbb{R}^{d}$ be a
bounded, polygonal ($d=2$) or polyhedral ($d=3$) Lipschitz domain and let
$\mathcal{G}$ be a regular simplicial finite element mesh for $\Omega$. Let
the diffusion matrix $\mathbb{A}\in L^{\infty}\left(  \Omega,\mathbb{R}%
_{\operatorname*{sym}}^{d\times d}\right)  $ satisfy assumption (\ref{aeps})
and let $f\in L^{2}\left(  \Omega\right)  $. As an additional assumption on
the regularity, we require that the exact solution of (\ref{varform})
satisfies $u\in PH{^{1+s}\left(  \Omega\right)  }$ for some positive integer
$s$ and $\left\Vert \mathbb{A}\right\Vert _{{PW^{{r},\infty}\left(
\Omega\right)  }}<\infty$ holds with $r:=\min\left\{  p,s\right\}  $. Let the
continuous problem (\ref{varform}) be discretized by the non-conforming
Galerkin method (\ref{ncGalDisc}) with a finite dimensional space $S$ which
satisfies $S_{\mathcal{G},\operatorname*{c}}^{p}\subset S\subset
S_{\mathcal{G}}^{p}$ on{ a compatible mesh $\mathcal{G}$}. Then,
(\ref{ncGalDisc}) has a unique solution which satisfies%
\[
\left\vert u-u_{S}\right\vert _{H_{\operatorname*{pw}}^{1}\left(
\Omega\right)  }\leq Ch^{r}\left\Vert u\right\Vert _{PH{^{1+r}\left(
\Omega\right)  }}.
\]
The constant $C$ only depends on $a_{\min}$, $a_{\max}$, $\left\Vert
\mathbb{A}\right\Vert _{{PW^{{r},\infty}\left(  \Omega\right)  }}$, $p$, $r$,
and the shape regularity of the mesh.
\end{theorem}

%

\proof
The second Strang lemma (cf. \cite[Theo. 4.2.2]{Ciarlet}) applied to the
non-conforming Galerkin discretization (\ref{ncGalDisc}) implies the existence
of a unique solution which satisfies the error estimate%
\[
\left\vert u-u_{S}\right\vert _{H_{\operatorname*{pw}}^{1}\left(
\Omega\right)  }\leq\left(  1+\frac{a_{\max}}{a_{\min}}\right)  \inf_{v\in
S}\left\vert u-v\right\vert _{H_{\operatorname*{pw}}^{1}\left(  \Omega\right)
}+\frac{1}{a_{\min}}\sup_{v\in S}\frac{\left\vert \mathcal{L}_{u}\left(
v\right)  \right\vert }{\left\vert v\right\vert _{H_{\operatorname*{pw}}%
^{1}\left(  \Omega\right)  }},
\]
where%
\[
\mathcal{L}_{u}\left(  v\right)  :=a_{\mathcal{G}}\left(  u,v\right)  -\left(
f,v\right)  .
\]

The approximation properties of $S$ are inherited from the approximation
properties of $S_{\mathcal{G},\operatorname*{c}}^{p}$ in the first infimum
because of the inclusion $S_{\mathcal{G},\operatorname*{c}}^{p}\subset S$. For
the second term we obtain%
\begin{equation}
\mathcal{L}_{u}\left(  v\right)  =\left(  \mathbb{A}\nabla u,\nabla
_{\mathcal{G}}v\right)  -\left(  f,v\right)  . \label{defL-piecewise}%
\end{equation}
Note that $f\in L^{2}\left(  \Omega\right)  $ implies that $\operatorname{div}%
\left(  \mathbb{A}\nabla u\right)  \in L^{2}\left(  \Omega\right)  $ and, in
turn, that the normal jump $\left[  \mathbb{A}\nabla u\cdot\mathbf{n}%
_{T}\right]  _{T}$ equals zero and the restriction $\left.  \left(
\mathbb{A}\nabla u\cdot\mathbf{n}_{T}\right)  \right\vert _{T}$ is well
defined for all $T\in\mathcal{F}$. We may apply simplexwise integration by
parts to (\ref{defL-piecewise}) to obtain%
\[
\mathcal{L}_{u}\left(  v\right)  =-\sum_{T\in\mathcal{F}_{\Omega}}\,\int
_{T}\left(  \mathbb{A}\nabla u\cdot\mathbf{n}_{T}\right)  \left[  v\right]
_{T}+\sum_{T\in\mathcal{F}_{\partial\Omega}}\int_{T}\left(  \mathbb{A}\nabla
u\cdot\mathbf{n}_{T}\right)  v.
\]

Let $K_{T}${ be one simplex in $\omega_{T}$. For }$1\leq i\leq d$,{ let}
$q_{i}\in\mathbb{P}_{d}^{p-1}\left(  K_{T}\right)  $ denote the best
approximation of $w_{i}:=\left.  \left(  \sum_{j=1}^{d}A_{i,j}\partial
_{j}u\right)  \right\vert _{K_{T}}$ with respect to the $H^{1}\left(
K_{T}\right)  $ norm. Then, $\left.  q_{i}\right\vert _{T}n_{T,i}\in
\mathbb{P}_{d-1}^{p-1}\left(  T\right)  $ for $1\leq i\leq d$, and the
inclusion $S\subset S_{\mathcal{G}}^{p}$ implies%
\begin{align}
\left\vert \mathcal{L}_{u}\left(  v\right)  \right\vert \leq &  \left\vert
-\sum_{T\in\mathcal{F}_{\Omega}}\int_{T}\left(  \sum_{i=1}^{d}\left(
w_{i}-q_{i}\right)  \cdot n_{T,i}\right)  \left[  v\right]  _{T}\right\vert
\label{Luest1}\\
&  +\left\vert \sum_{T\in\mathcal{F}_{\partial\Omega}}\int_{T}\left(
\sum_{i=1}^{d}\left(  w_{i}-q_{i}\right)  \cdot n_{T,i}\right)  v\right\vert
\nonumber\\
\leq &  \sum_{T\in\mathcal{F}_{\Omega}}\left\Vert \left[  v\right]
_{T}\right\Vert _{L^{2}\left(  T\right)  }\sum_{i=1}^{d}\left\Vert w_{i}%
-q_{i}\right\Vert _{L^{2}\left(  T\right)  }\nonumber\\
&  +\sum_{T\in\mathcal{F}_{\partial\Omega}}\left\Vert v\right\Vert
_{L^{2}\left(  T\right)  }\sum_{i=1}^{d}\left\Vert w_{i}-q_{i}\right\Vert
_{L^{2}\left(  T\right)  }.\nonumber
\end{align}
Standard trace estimates and approximation properties lead to%
\begin{align}
\left\Vert w_{i}-q_{i}\right\Vert _{L^{2}\left(  T\right)  }  &  \leq C\left(
h_{T}^{-1/2}\left\Vert w_{i}-q_{i}\right\Vert _{L^{2}\left(  K_{T}\right)
}+h_{T}^{1/2}\left\vert w_{i}-q_{i}\right\vert _{H^{1}\left(  K_{T}\right)
}\right) \label{Luest2}\\
&  \leq Ch_{T}^{r-1/2}\left\vert w_{i}\right\vert _{H^{r}\left(  K_{T}\right)
}\leq Ch_{T}^{r-1/2}\left\Vert u\right\Vert _{H^{1+r}\left(  K_{T}\right)
},\nonumber
\end{align}
where $C$ depends only on $p$, $r$, $\left\Vert \mathbb{A}\right\Vert
_{{W^{r}\left(  K_{T}\right)  }}$, and the shape regularity of the mesh{.}The
combination of (\ref{Luest1}), (\ref{Luest2}) and (\ref{poincarecondspa}%
),(\ref{poincarecondspb}) along with the shape regularity of the mesh leads to
the consistency estimate%
\begin{align*}
\left\vert \mathcal{L}_{u}\left(  v\right)  \right\vert  &  \leq C\left(
\sum_{T\in\mathcal{F}_{\Omega}}h_{T}^{r}\left\Vert u\right\Vert _{H^{1+r}%
\left(  K_{T}\right)  }\left\vert v\right\vert _{H_{\operatorname*{pw}}%
^{1}\left(  \omega_{T}\right)  }+\sum_{T\in\mathcal{F}_{\partial\Omega}}%
h_{T}^{r}\left\Vert u\right\Vert _{H^{1+r}\left(  K_{T}\right)  }\left\vert
v\right\vert _{H_{\operatorname*{pw}}^{1}\left(  \omega_{T}\right)  }\right)
\\
&  \leq\tilde{C}h^{r}\left\Vert u\right\Vert _{PH{^{1+r}\left(  \Omega\right)
}}\left\vert v\right\vert _{H_{\operatorname*{pw}}^{1}\left(  \Omega\right)
},
\end{align*}
which completes the proof.\hfill%
\endproof

\begin{remark}
If one chooses in (\ref{poinprep}) a degree $p^{\prime}<p$ for the
orthogonality relations in (\ref{poinprep}), then the order of convergence
behaves like $h^{r^{\prime}}\left\Vert e\right\Vert _{H^{1+r^{\prime}}\left(
\Omega\right)  }$, with $r^{\prime}:=\min\left\{  p^{\prime},s\right\}  $,
because the best approximations $q_{i}$ now belong to $P_{d-1}^{p^{\prime}%
-1}\left(  T\right)  $.
\end{remark}

\section{Explicit Construction of Non-Conforming Crouzeix-Raviart Finite
Elements\label{SecExplConstr}}

\subsection{Jacobi Polynomials}

Let $\alpha,\beta>-1$. The \emph{Jacobi polynomial} $P_{n}^{\left(
\alpha,\beta\right)  }$ is a polynomial of degree $n$ such that
\[
\int_{-1}^{1}P_{n}^{\left(  \alpha,\beta\right)  }\left(  x\right)  \,q\left(
x\right)  \left(  1-x\right)  ^{\alpha}\left(  1+x\right)  ^{\beta}\,dx=0
\]
for all polynomials $q$ of degree less than $n$, and (cf. \cite[Table
18.6.1]{NIST:DLMF})%
\begin{equation}
P_{n}^{\left(  \alpha,\beta\right)  }\left(  1\right)  =\frac{\left(
\alpha+1\right)  _{n}}{n!},\qquad P_{n}^{\left(  \alpha,\beta\right)  }\left(
-1\right)  =\left(  -1\right)  ^{n}\frac{\left(  \beta+1\right)  _{n}}{n!}.
\label{Pnormalization}%
\end{equation}
Here the \emph{shifted factorial} is defined by $\left(  a\right)
_{n}:=a\left(  a+1\right)  \ldots\left(  a+n-1\right)  $ for $n>0$ and
$\left(  a\right)  _{0}:=1$. The Jacobi polynomial has an explicit expression
in terms of a \emph{terminating Gauss hypergeometric series} (see (cf.
\cite[18.5.7]{NIST:DLMF}))%
\begin{equation}
\,\mbox{}_{2}F_{1}\!\left(
\genfrac{}{}{0pt}{}{-n,b}{c}%
;z\right)  :=\sum_{k=0}^{n}\frac{\left(  -n\right)  _{k}\left(  b\right)
_{k}}{\left(  c\right)  _{k}\,k!}\,z^{k} \label{F21rep}%
\end{equation}
as follows%
\begin{equation}
P_{n}^{\left(  \alpha,\beta\right)  }\left(  x\right)  =\frac{\left(
\alpha+1\right)  _{n}}{n!}\,\,\mbox{}_{2}F_{1}\!\left(
\genfrac{}{}{0pt}{}{-n,n+\alpha+\beta+1}{\alpha+1}%
;\frac{1-x}{2}\right)  . \label{defJP}%
\end{equation}

\subsection{Orthogonal Polynomials on Triangles\label{SecOrthoPolyT}}

Recall that $\widehat{T}$ is the (closed) unit triangle in $\mathbb{R}^{2}$
with vertices $\widehat{\mathbf{A}}_{0}=\left(  0,0\right)  ^{\intercal}$,
$\widehat{\mathbf{A}}_{1}=\left(  1,0\right)  ^{\intercal}$, and
$\widehat{\mathbf{A}}_{3}=\left(  0,1\right)  ^{\intercal}$. An orthogonal
basis for the space $\mathbb{P}_{n,n-1}^{\perp}\left(  \widehat{T}\right)  $
was introduced in \cite{Proriol_ortho} and is given by the functions $b_{n,k}%
$, $0\leq k\leq n$,
\begin{equation}
b_{n,k}(\mathbf{x}):=\left(  x_{1}+x_{2}\right)  ^{k}\,P_{n-k}^{(0,2k+1)}%
\left(  2\left(  x_{1}+x_{2}\right)  -1\right)  \,P_{k}^{\left(  0,0\right)
}\left(  \frac{x_{1}-x_{2}}{x_{1}+x_{2}}\right)  , \label{defbnk}%
\end{equation}
where $P_{k}^{\left(  0,0\right)  }$ are the Legendre polynomials (see
\cite[18.7.9]{NIST:DLMF})\footnote{The Legendre polynomials with normalization
$P_{k}^{\left(  0,0\right)  }\left(  1\right)  =1$ for all $k=0,1,\ldots$ can
be defined \cite[Table 18.9.1]{NIST:DLMF} via the three-term recursion%
\begin{equation}
P_{0}^{\left(  0,0\right)  }\left(  x\right)  =1;\quad P_{1}^{\left(
0,0\right)  }\left(  x\right)  =x;\quad\text{and\quad}\left(  k+1\right)
P_{k+1}^{\left(  0,0\right)  }\left(  x\right)  =\left(  2k+1\right)
xP_{k}^{\left(  0,0\right)  }\left(  x\right)  -kP_{k-1}^{\left(  0,0\right)
}\left(  x\right)  \quad\text{for }k=1,2,\ldots, \label{rekLk}%
\end{equation}
from which the well-known relation $P_{k}^{\left(  0,0\right)  }\left(
x\right)  =\left(  -1\right)  ^{k}P_{k}^{\left(  0,0\right)  }\left(
x\right)  $ for all $k\in\mathbb{N}_{0}$ follows.}. From (\ref{rekLk})
(footnote) it follows that these polynomials satisfy the following symmetry
relation%
\begin{equation}
b_{n,k}\left(  x_{1},x_{2}\right)  =\left(  -1\right)  ^{k}b_{n,k}\left(
x_{2},x_{1}\right)  \quad\forall n\geq0,\forall\left(  x_{1},x_{2}\right)  .
\label{symmbnk}%
\end{equation}
By combining (\ref{F21rep}) - (\ref{defbnk}), an elementary calculation leads
to\footnote{Further special values are%
\[%
\begin{array}
[c]{ll}%
b_{n,0}\left(  0,0\right)  =P_{n}^{\left(  0,1\right)  }\left(  -1\right)
=\left(  -1\right)  ^{n}\frac{\left(  2\right)  _{n}}{n!}=\left(  -1\right)
^{n}\left(  n+1\right)  , & b_{n,k}\left(  0,0\right)  =0,~1\leq k\leq n,\\
b_{n,k}\left(  1,0\right)  =P_{n-k}^{\left(  0,2k+1\right)  }\left(  1\right)
P_{k}^{\left(  0,0\right)  }\left(  1\right)  =1,0\leq k\leq n, &
b_{n,k}\left(  0,1\right)  =P_{n-k}^{\left(  0,2k+1\right)  }\left(  1\right)
P_{k}^{\left(  0,0\right)  }\left(  -1\right)  =\left(  -1\right)  ^{k},~0\leq
k\leq n.
\end{array}
\]
} $b_{n,0}\left(  0,0\right)  =\left(  -1\right)  ^{n}\left(  n+1\right)  $.

Let
\begin{equation}
E^{\operatorname{I}}:=\overline{\widehat{\mathbf{A}}_{0}\widehat{\mathbf{A}%
}_{1}}\text{,\quad}E^{\operatorname{II}}:=\overline{\widehat{\mathbf{A}}%
_{0}\widehat{\mathbf{A}}_{2}}\text{,\quad and\quad}E^{\operatorname{III}%
}:=\overline{\widehat{\mathbf{A}}_{1}\widehat{\mathbf{A}}_{2}}
\label{defedges}%
\end{equation}
denote the edges of $\widehat{T}$. For $\operatorname{Z}\in\left\{
\operatorname{I},\operatorname{II},\operatorname{III}\right\}  $, we introduce
the linear restriction operator for the edge $E^{\operatorname{Z}}$ by
$\gamma^{\operatorname{Z}}:C^{0}\left(  \widehat{T}\right)  \rightarrow
C^{0}\left(  \left[  0,1\right]  \right)  $ by%
\begin{equation}
\gamma^{\operatorname{I}}u:=u\left(  \cdot,0\right)  ,\quad\gamma
^{\operatorname{II}}u:=u\left(  0,\cdot\right)  ,\quad\gamma
^{\operatorname{III}}u=u\left(  1-\cdot,\cdot\right)  \label{defgamma}%
\end{equation}
which allows to define%
\[
b_{n,k}^{\operatorname{I}}:=\gamma^{\operatorname{I}}b_{n,k},\quad
b_{n,k}^{\operatorname{II}}:=\gamma^{\operatorname{II}}b_{n,k},\quad
b_{n,k}^{\operatorname{III}}:=\gamma^{\operatorname{III}}b_{n,k}%
,\quad\text{for }k=0,1,\ldots,n.
\]

\begin{lemma}
\label{LemRestLinIndep}For any $\operatorname{Z}\in\left\{  \operatorname{I}%
,\operatorname{II},\operatorname{III}\right\}  $, each of the systems $\left(
b_{n,k}^{\operatorname{Z}}\right)  _{k=0}^{n}$, form a basis of $\mathbb{P}%
_{n}\left(  \left[  0,1\right]  \right)  $.
\end{lemma}

%

\proof
First note that $\left\{  x^{j}\left(  x-1\right)  ^{n-j}:0\leq j\leq
n\right\}  $ is a basis for $\mathbb{P}_{n}\left(  \left[  0,1\right]
\right)  $; this follows from expanding the right-hand side of $x^{m}%
=x^{m}\left(  x-\left(  x-1\right)  \right)  ^{n-m}$. Specialize the formula
\cite[18.5.8]{NIST:DLMF}%
\[
P_{m}^{\left(  \alpha,\beta\right)  }\left(  s\right)  =\frac{\left(
\alpha+1\right)  _{m}}{m!}\left(  \frac{1+s}{2}\right)  ^{m}~_{2}F_{1}\left(
\genfrac{}{}{0pt}{}{-m,-m-\beta}{\alpha+1}%
;\frac{s-1}{s+1}\right)
\]
to $m=n-k$, $\alpha=0,\beta=2k+1,s=2x-1$ to obtain%
\begin{align}
b_{n,k}^{\operatorname{I}}\left(  x\right)   &  =x^{n}~_{2}F_{1}\left(
\genfrac{}{}{0pt}{}{k-n,-n-k-1}{1}%
;\frac{x-1}{x}\right) \label{defbtilde}\\
&  \overset{\text{(\ref{F21rep})}}{=}\sum_{i=0}^{n-k}\frac{\left(  k-n\right)
_{i}\left(  -n-k-1\right)  _{i}}{i!i!}x^{n-i}\left(  x-1\right)  ^{i}.
\label{clnk}%
\end{align}
The highest index $i$ of $x^{n-i}\left(  x-1\right)  ^{i}$ in $b_{n,k}%
^{\operatorname{I}}\left(  x\right)  $ is $n-k$ with coefficient
$\dfrac{\left(  2k+2\right)  _{n-k}}{\left(  n-k\right)  !}\neq0$. Thus the
matrix expressing $\left[  b_{n,0}^{\operatorname{I}},\ldots,b_{n,n}%
^{\operatorname{I}}\right]  $ in terms of $\left[  \left(  x-1\right)
^{n},x\left(  x-1\right)  ^{n-1},\ldots,x^{n}\right]  $ is triangular and
nonsingular; hence $\left\{  b_{n,k}^{\operatorname{I}}:0\leq k\leq n\right\}
$ is a basis of $\mathbb{P}_{n}\left(  \left[  0,1\right]  \right)  $. The
symmetry relation $b_{n,k}^{\operatorname{II}}=\left(  -1\right)  ^{k}%
b_{n,k}^{\operatorname{I}}$ for $0\leq k\leq n$ (cf. (\ref{symmbnk})) shows
that $\left\{  b_{n,k}^{\operatorname{II}}:0\leq k\leq n\right\}  $ is also a
basis of $\mathbb{P}_{n}\left(  \left[  0,1\right]  \right)  $. Finally
substituting $x_{1}=1-x,x_{2}=x$ in $b_{n,k}$ results in%
\begin{equation}
b_{n,k}^{\operatorname{III}}\left(  x\right)  =P_{n-k}^{\left(  0,2k+1\right)
}\left(  1\right)  P_{k}^{\left(  0,0\right)  }\left(  1-2x\right)  ,
\label{bnkIII}%
\end{equation}
and $P_{n-k}^{\left(  0,2k+1\right)  }\left(  1\right)  =1$ (from
(\ref{Pnormalization})). Clearly $\left\{  P_{k}^{\left(  0,0\right)  }\left(
1-2x\right)  :0\leq k\leq n\right\}  $ is a basis for $\mathbb{P}_{n}\left(
\left[  0,1\right]  \right)  $.%
\endproof

\begin{lemma}
\label{TraceLemma}Let $v\in\mathbb{P}_{n}\left(  \left[  0,1\right]  \right)
$. Then, there exist unique orthogonal polynomials $u^{\operatorname{Z}}%
\in\mathbb{P}_{n,n-1}^{\perp}\left(  \widehat{T}\right)  $, $\operatorname{Z}%
\in\left\{  \operatorname{I},\operatorname{II},\operatorname{III}\right\}  $
with $v=\gamma^{\operatorname{Z}}u^{\operatorname{Z}}$. Thus, the linear
extension operator $\mathcal{E}^{\operatorname{Z}}:\mathbb{P}_{n}\left(
\left[  0,1\right]  \right)  \rightarrow\mathbb{P}_{n,n-1}^{\perp}\left(
\widehat{T}\right)  $ is well defined by $\mathcal{E}^{\operatorname{Z}%
}v:=u^{\operatorname{Z}}$.
\end{lemma}

%

\proof
From Lemma \ref{LemRestLinIndep} we conclude that $\gamma^{\operatorname{Z}}$
is surjective. Since the polynomial spaces are finite dimensional the
assertion follows from%
\[
\dim\mathbb{P}_{n}\left(  \left[  0,1\right]  \right)  =n+1=\dim
\mathbb{P}_{n,n-1}^{\perp}\left(  \widehat{T}\right)  .
\]%
\endproof

The orthogonal polynomials can be lifted to a general triangle $T$.

\begin{definition}
Let $T$ denote a triangle and $\chi_{T}$ an affine pullback to the reference
triangle $\widehat{T}$. Then, the space of orthogonal polynomials of degree
$n$ on $T$ is%
\[
\mathbb{P}_{n,n-1}^{\perp}\left(  T\right)  :=\left\{  v\circ\chi_{T}%
^{-1}:v\in\mathbb{P}_{n,n-1}^{\perp}\left(  \widehat{T}\right)  \right\}  .
\]

\end{definition}

From the transformation rule for integrals one concludes that for any
$u=v\circ\chi_{T}^{-1}\in\mathbb{P}_{n,n-1}^{\perp}\left(  T\right)  $ and all
$q\in\mathbb{P}_{n-1}\left(  T\right)  $ it holds%
\begin{equation}
\int_{T}uq=\int_{T}\left(  v\circ\chi_{T}^{-1}\right)  q=2\left\vert
T\right\vert \int_{\widehat{T}}v\left(  q\circ\chi_{T}\right)  =0
\label{orthogonalityT}%
\end{equation}
since $q\circ\chi_{T}\in\mathbb{P}_{n-1}\left(  \widehat{T}\right)  $. Here
$\left\vert T\right\vert $ denotes the area of the triangle $T$.

\subsection{Totally Symmetric Orthogonal Polynomials\label{SecFullSymPoly}}

In this section, we will decompose the space of orthogonal polynomials
$\mathbb{P}_{n,n-1}^{\perp}\left(  \widehat{T}\right)  $ into three
irreducible modules (see \S \ref{SubSecIrr}) and thus, obtain a direct sum
decomposition $\mathbb{P}_{n,n-1}^{\perp}\left(  \widehat{T}\right)
=\mathbb{P}_{n,n-1}^{\perp,\operatorname*{sym}}\left(  \widehat{T}\right)
\oplus\mathbb{P}_{n,n-1}^{\perp,\operatorname*{refl}}\left(  \widehat
{T}\right)  \oplus\mathbb{P}_{n,n-1}^{\perp,\operatorname{sign}}\left(
\widehat{T}\right)  $. We will derive an explicit representation for a basis
of the space of totally symmetric polynomials $\mathbb{P}_{n,n-1}%
^{\perp,\operatorname*{sym}}\left(  \widehat{T}\right)  $ in
\S \ref{ConSymmBas} and of the space of reflection symmetric polynomials
$\mathbb{P}_{n,n-1}^{\perp,\operatorname*{refl}}\left(  \widehat{T}\right)  $
in \S \ref{Sectaureflcomp}.

We start by introducing, for functions on triangles, the notation of total
symmetry. For an arbitrary triangle $T$ with vertices $\mathbf{A}_{0}$,
$\mathbf{A}_{1}$, $\mathbf{A}_{2}$, we introduce the set of permutations
$\Pi=\left\{  \left(  i,j,k\right)  :i,j,k\in\left\{  0,1,2\right\}  \text{
pairwise disjoint}\right\}  $. For $\pi=\left(  i,j,k\right)  \in\Pi$, define
the affine mapping $\chi_{\pi}:T\rightarrow T$ by
\begin{equation}
\chi_{\pi}\left(  \mathbf{x}\right)  =\mathbf{A}_{i}+x_{1}\left(
\mathbf{A}_{j}-\mathbf{A}_{i}\right)  +x_{2}\left(  \mathbf{A}_{k}%
-\mathbf{A}_{i}\right)  . \label{defchipi}%
\end{equation}
We say a function $u$, defined on $T$, has \textit{total symmetry} if%
\[
u=u\circ\chi_{\pi}\quad\forall\pi\in\Pi\text{.}%
\]

The space of totally symmetric orthogonal polynomials is%
\begin{equation}
\mathbb{P}_{n,n-1}^{\perp,\operatorname*{sym}}\left(  \widehat{T}\right)
:=\left\{  u\in\mathbb{P}_{n,n-1}^{\perp}\left(  \widehat{T}\right)  :u\text{
has total symmetry}\right\}  . \label{defpsymhat}%
\end{equation}

The construction of a basis of $\mathbb{P}_{n,n-1}^{\perp,\operatorname*{sym}%
}\left(  \widehat{T}\right)  $ requires some algebraic tools which we develop
in the following.

\subsubsection{The decomposition of $\mathbb{P}_{n,n-1}^{\perp}\left(
\widehat{T}\right)  $ or $\mathbb{P}_{n}\left(  \left[  0,1\right]  \right)  $
into irreducible $\mathcal{S}_{3}$ modules\label{SubSecIrr}}

We use the operator $\mathcal{\gamma}^{\operatorname{I}}$ (cf. (\ref{defgamma}%
)) to set up an action of the symmetric group $\mathcal{S}_{3}$ on
$\mathbb{P}_{n}\left(  \left[  0,1\right]  \right)  $ by transferring its
action on $\mathbb{P}_{n,n-1}^{\perp}\left(  \widehat{T}\right)  $ on the
basis $\left\{  b_{n,k}\right\}  $. It suffices to work with two generating
reflections. On the triangle $\chi_{\left\{  0,2,1\right\}  }\left(
x_{1},x_{2}\right)  =\left(  x_{2},x_{1}\right)  $ and thus $b_{n,k}\circ
\chi_{\left\{  0,2,1\right\}  }=\left(  -1\right)  ^{k}b_{n,k}$ (this follows
from (\ref{symmbnk})). The action of $\chi_{\left\{  0,2,1\right\}  }$ is
mapped to $\sum_{k=0}^{n}\alpha_{k}b_{n,k}^{\operatorname{I}}\mapsto\sum
_{k=0}^{n}\left(  -1\right)  ^{k}\alpha_{k}b_{n,k}^{\operatorname{I}}$, and
denoted by $R$. For the other generator we use $\chi_{\left\{  1,0,2\right\}
}\left(  x_{1},x_{2}\right)  =\left(  1-x_{1}-x_{2},x_{2}\right)  $. Under
$\gamma^{\operatorname{I}}$ this corresponds to the map $\sum_{k=0}^{n}%
\alpha_{k}b_{n,k}^{\operatorname{I}}\left(  x\right)  \mapsto\sum_{k=0}%
^{n}\alpha_{k}b_{n,k}^{\operatorname{I}}\left(  1-x\right)  $ which is denoted
by $M$. We will return later to transformation formulae expressing%
\[
b_{n,k}\circ\chi_{\left\{  1,0,2\right\}  }\left(  x_{1},x_{2}\right)
=\left(  1-x_{1}\right)  ^{k}P_{n-k}^{\left(  0,2k+1\right)  }\left(
1-2x_{1}\right)  P_{k}^{\left(  0,0\right)  }\left(  \frac{1-x_{1}-2x_{2}%
}{1-x_{1}}\right)
\]
in the $\left\{  b_{n,k}\right\}  $-basis. Observe that $\left(  MR\right)
^{3}=I$ because $\chi_{\left\{  1,0,2\right\}  }\circ\chi_{\left\{
0,2,1\right\}  }\left(  x_{1},x_{2}\right)  =\left(  1-x_{1}-x_{2}%
,x_{1}\right)  $ and this mapping is of period 3. It follows that each of
$\left\{  M,R\right\}  $ and $\left\{  \chi_{\left\{  1,0,2\right\}  }%
,\chi_{\left\{  0,2,1\right\}  }\right\}  $ generates (an isomorphic copy of)
$\mathcal{S}_{3}$. It is a basic fact that the relations $M^{2}=I,R^{2}=I$ and
$\left(  MR\right)  ^{3}=I$ define $\mathcal{S}_{3}$. The representation
theory of $\mathcal{S}_{3}$ informs us that there are three nonisomorphic
irreducible representations:%
\begin{align*}
\tau_{\operatorname*{triv}}  &  :\chi_{\left\{  0,2,1\right\}  }%
\rightarrow1,\chi_{\left\{  1,0,2\right\}  }\rightarrow1;\\
\tau_{\operatorname{sign}}  &  :\chi_{\left\{  0,2,1\right\}  }\rightarrow
-1,\chi_{\left\{  1,0,2\right\}  }\rightarrow-1;\\
\tau_{\operatorname*{refl}}  &  :\chi_{\left\{  0,2,1\right\}  }%
\rightarrow\sigma_{1}:=%
\begin{bmatrix}
-1 & 0\\
0 & 1
\end{bmatrix}
,\chi_{\left\{  1,0,2\right\}  }\rightarrow\sigma_{2}:=%
\begin{bmatrix}
\frac{1}{2} & 1\\
\frac{3}{4} & -\frac{1}{2}%
\end{bmatrix}
.
\end{align*}
(The subscript \textquotedblleft\textit{$\operatorname*{refl}$}%
\textquotedblright\textit{ }designates the \textit{reflection }%
representation).\textit{ }Then the eigenvectors of $\sigma_{1},\sigma_{2}$
with $-1$ as eigenvalue are $\left(  -1,0\right)  ^{\intercal}$ and $\left(
2,-3\right)  ^{\intercal}$ respectively; these two vectors are a basis for
$\mathbb{R}^{2}$. Similarly the eigenvectors of $\sigma_{1}$ and $\sigma_{2}$
with eigenvalue $+1$, namely $\left(  0,1\right)  ^{\intercal}$, $\left(
2,1\right)  ^{\intercal}$, form a basis. Form a direct sum
\[
\mathbb{P}_{n,n-1}^{\perp}\left(  \widehat{T}\right)  :=\left(
{\displaystyle\bigoplus\limits_{j\geq0}}
E_{j}^{\left(  \operatorname*{triv}\right)  }\right)  \oplus\left(
{\displaystyle\bigoplus\limits_{j\geq0}}
E_{j}^{\left(  \operatorname{sign}\right)  }\right)  \oplus\left(
{\displaystyle\bigoplus\limits_{j\geq0}}
E_{j}^{\left(  \operatorname*{refl}\right)  }\right)  ,
\]
where the $E_{j}^{\left(  \operatorname*{triv}\right)  },E_{j}^{\left(
\operatorname{sign}\right)  },E_{j}^{\left(  \operatorname*{refl}\right)  }$
are $\mathcal{S}_{3}$-irreducible and realizations of the representations
$\tau_{\operatorname*{triv}},\tau_{\operatorname{sign}},\tau
_{\operatorname*{refl}}$ respectively. Let $d_{\operatorname*{triv}}\left(
n\right)  ,d_{\operatorname{sign}}\left(  n\right)  ,d_{\operatorname*{refl}%
}\left(  n\right)  $ denote the respective multiplicities, so that
$d_{\operatorname*{triv}}\left(  n\right)  +d_{\operatorname{sign}}\left(
n\right)  +2d_{\operatorname*{refl}}\left(  n\right)  =n+1$. The case $n$ even
or odd are handled separately. If $n=2m$ is even then the number of
eigenvectors of $R$ having $-1$ as eigenvalue equals $m$ (the cardinality of
$\left\{  1,3,5,\ldots,2m-1\right\}  $). The same property holds for $M$ since
the eigenvectors of $M$ in the basis $\left\{  x^{2m}\left(  x-1\right)
^{2m-j}\right\}  $ are explicitly given by $\left\{  x^{2m-2\ell}\left(
x-1\right)  ^{2\ell}-x^{2\ell}\left(  x-1\right)  ^{2m-2\ell}:0\leq\ell\leq
m\right\}  $. Each $E_{j}^{\left(  \operatorname*{refl}\right)  }$ contains
one $\left(  -1\right)  $-eigenvector of $\chi_{\left\{  1,0,2\right\}  }$ and
one of $\chi_{\left\{  0,2,1\right\}  }$ and each $E_{j}^{\left(
\operatorname{sign}\right)  }$ consists of one $\left(  -1\right)
$-eigenvector of $\chi_{\left\{  0,2,1\right\}  }$. This gives the equation
$d_{\operatorname*{refl}}\left(  n\right)  +d_{\operatorname{sign}}\left(
n\right)  =m$. Each $E_{j}^{\left(  \operatorname*{refl}\right)  }$ contains
one $\left(  +1\right)  $-eigenvector of $\chi_{\left\{  1,0,2\right\}  }$ and
one of $\chi_{\left\{  0,2,1\right\}  }$ and each $E_{j}^{\left(
\operatorname*{triv}\right)  }$ consists of one $\left(  +1\right)
$-eigenvector of $\chi_{\left\{  0,2,1\right\}  }$. There are $m+1$
eigenvectors with eigenvalue $1$ of each of $\chi_{\left\{  1,0,2\right\}  }$
and $\chi_{\left\{  0,2,1\right\}  }$ thus $d_{\operatorname*{refl}}\left(
n\right)  +d_{\operatorname*{triv}}\left(  n\right)  =m+1$.

If $n=2m+1$ is odd then the eigenvector multiplicities are $m+1$ for both
eigenvalues $+1,-1$. By similar arguments we obtain the equations
$d_{\operatorname*{refl}}\left(  n\right)  +d_{\operatorname{sign}}\left(
n\right)  =m+1$, $d_{\operatorname*{refl}}\left(  n\right)
+d_{\operatorname*{triv}}\left(  n\right)  =m+1$. It remains to find one last
relation for both, even and odd cases.

To finish the determination of the multiplicities $d_{\operatorname*{triv}%
}\left(  n\right)  ,d_{\operatorname{sign}}\left(  n\right)
,d_{\operatorname*{refl}}\left(  n\right)  $ it suffices to find
$d_{\operatorname*{triv}}\left(  n\right)  $. This is the dimension of the
space of polynomials in $\mathbb{P}_{n,n-1}^{\perp}\left(  \widehat{T}\right)
$ which are invariant under both $\chi_{\left\{  0,2,1\right\}  }$ and
$\chi_{\left\{  1,0,2\right\}  }.$ Since these two group elements generate
$\mathcal{S}_{3}$ this is equivalent to being invariant under each element of
$\mathcal{S}_{3}$ .This property is called \textit{totally symmetric}. Under
the action of $\gamma^{\operatorname{I}}$ this corresponds to the space of
polynomials in $\mathbb{P}_{n}\left(  \left[  0,1\right]  \right)  $ which are
invariant under both $R$ and $M$. We appeal to the classical theory of
symmetric polynomials: suppose $\mathcal{S}_{3}$ acts on polynomials in
$\left(  y_{1},y_{2},y_{3}\right)  $ by permutation of coordinates then the
space of symmetric (invariant under the group) polynomials is exactly the
space of polynomials in $\left\{  e_{1},e_{2},e_{3}\right\}  $ the elementary
symmetric polynomials, namely $e_{1}=y_{1}+y_{2}+y_{3}$, $e_{2}=y_{1}%
y_{2}+y_{1}y_{3}+y_{2}y_{3}$, $e_{3}=y_{1}y_{2}y_{3}$. To apply this we set up
an affine map from $\widehat{T}$ to the triangle in $\mathbb{R}^{3}$ with
vertices $\left(  2,-1,-1\right)  $, $\left(  -1,2,-1\right)  $, $\left(
-1,-1,2\right)  $. The formula for the map is%
\[
y\left(  x\right)  =\left(  2-3x_{1}-3x_{2},3x_{1}-1,3x_{2}-1\right)  .
\]
The map takes $\left(  0,0\right)  ,\left(  1,0\right)  ,\left(  0,1\right)  $
to the three vertices respectively. The result is%
\begin{align*}
e_{1}\left(  y\left(  x\right)  \right)   &  =0,\\
e_{2}\left(  y\left(  x\right)  \right)   &  =-9\left(  x_{1}^{2}+x_{1}%
x_{2}+x_{2}^{2}-x_{1}-x_{2}\right)  -3,\\
e_{3}\left(  y\left(  x\right)  \right)   &  =\left(  3x_{1}-1\right)  \left(
3x_{2}-1\right)  \left(  2-3x_{1}-3x_{2}\right)  .
\end{align*}
Thus any totally symmetric polynomial on $\widehat{T}$ is a linear combination
of $e_{2}^{a}e_{3}^{b}$ with uniquely determined coefficients. The number of
linearly independent totally symmetric polynomials in $\left(
{\displaystyle\bigoplus\limits_{j=1}^{n}}
\mathbb{P}_{n,n-1}^{\perp}\left(  \widehat{T}\right)  \right)  \oplus
\mathbb{P}_{0}\left(  \widehat{T}\right)  $ equals the number of solutions of
$0\leq2a+3b\leq n$ with $a,b=0,1,2,\ldots$. As a consequence
$d_{\operatorname*{triv}}\left(  n\right)  =\operatorname*{card}\left\{
\left(  a,b\right)  :2a+3b=n\right\}  $. This number is the coefficient of
$t^{n}$ in the power series expansion of
\[
\frac{1}{\left(  1-t^{2}\right)  \left(  1-t^{3}\right)  }=\left(
1+t^{2}+t^{3}+t^{4}+t^{5}+t^{7}\right)  \left(  1+2t^{6}+3t^{12}%
+\ldots\right)  .
\]
From $d_{\operatorname*{triv}}\left(  n\right)  =\operatorname*{card}\left(
\left\{  0,2,4,\ldots\right\}  \cap\left\{  n,n-3,n-6,\ldots\right\}  \right)
$ we deduce the formula (cf. (\ref{dtrivfirst}))%
\[
d_{\operatorname*{triv}}\left(  n\right)  =\left\lfloor \frac{n}%
{2}\right\rfloor -\left\lfloor \frac{n-1}{3}\right\rfloor .
\]

As a consequence: if $n=2m$ then $d_{\operatorname{sign}}\left(  n\right)
=d_{\operatorname*{triv}}\left(  n\right)  -1$ and $d_{\operatorname*{refl}%
}\left(  n\right)  =m+1-d_{\operatorname*{triv}}\left(  n\right)  $; if
$n=2m+1$ then $d_{\operatorname{sign}}\left(  n\right)
=d_{\operatorname*{triv}}\left(  n\right)  $ and $d_{\operatorname*{refl}%
}\left(  n\right)  =m+1-d_{\operatorname*{triv}}\left(  n\right)  $. From this
the following can be derived: $d_{\operatorname{sign}}\left(  n\right)
=\left\lfloor \frac{n-1}{2}\right\rfloor -\left\lfloor \frac{n-1}%
{3}\right\rfloor $ and $d_{\operatorname*{refl}}\left(  n\right)
=\left\lfloor \frac{n+2}{3}\right\rfloor $. Here is a table of values in terms
of $n\operatorname{mod}6$:%
\[%
\begin{vmatrix}
n & d_{\operatorname*{triv}}\left(  n\right)  & d_{\operatorname{sign}}\left(
n\right)  & d_{\operatorname*{refl}}\left(  n\right) \\
6m & m+1 & m & 2m\\
6m+1 & m & m & 2m+1\\
6m+2 & m+1 & m & 2m+1\\
6m+3 & m+1 & m+1 & 2m+1\\
6m+4 & m+1 & m & 2m+2\\
6m+5 & m+1 & m+1 & 2m+2
\end{vmatrix}
.
\]

\subsubsection{Construction of totally symmetric polynomials\label{ConSymmBas}%
}

Let $M$ and $R$ denote the linear maps $Mp\left(  x_{1},x_{2}\right)
:=p\left(  1-x_{1}-x_{2},x_{2}\right)  $ and $Rp\left(  x_{1},x_{2}\right)
:=p\left(  x_{2},x_{1}\right)  $ respectively. Both are automorphisms of
$\mathbb{P}_{n,n-1}^{\perp}\left(  \widehat{T}\right)  $. Note $Mp=p\circ
\chi_{\left\{  1,0,2\right\}  }$ and $Rp=p\circ\chi_{\left\{  0,2,1\right\}
}$ (cf.\ Section \ref{SubSecIrr}).

\begin{proposition}
Suppose $0\leq k\leq n$ then%
\begin{align}
Rb_{n,k}  &  =\left(  -1\right)  ^{k}b_{n,k};\label{Rbnk}\\
Mb_{n,k}  &  =\left(  -1\right)  ^{n}\sum_{j=0}^{n}~_{4}F_{3}\left(
\genfrac{}{}{0pt}{}{-j,j+1,-k,k+1}{-n,n+2,1}%
;1\right)  \frac{2j+1}{n+1}b_{n,j}. \label{Mbnkeq}%
\end{align}

\end{proposition}%

\proof
The $_{4}F_{3}$-sum is understood to terminate at $k$ to avoid the $0/0$
ambiguities in the formal $_{4}F_{3}$-series. The first formula was shown in
Section \ref{SubSecIrr}. The second formula is a specialization of
transformations in \cite[Theorem 1.7(iii)]{Dunkl_3symm}: this paper used the
shifted Jacobi polynomial $R_{m}^{\left(  \alpha,\beta\right)  }\left(
s\right)  =\frac{m!}{\left(  \alpha+1\right)  _{m}}P_{m}^{\left(  \alpha
,\beta\right)  }\left(  1-2s\right)  $. Setting $\alpha=\beta=\gamma=0$ in the
formulas in \cite[Theorem 1.7(iii)]{Dunkl_3symm} results in $b_{n,k}=\left(
-1\right)  ^{k}\dfrac{\theta_{n,k}}{k!\left(  n-k\right)  !}$ and
$Mb_{n,k}=\dfrac{\phi_{n,k}}{k!\left(  n-k\right)  !}$, where $\theta_{n,k}$,
$\phi_{n,k}$ are the polynomials introduced in \cite[p.690]{Dunkl_3symm}. More
precisely, the arguments $v_{1},v_{2},v_{3}$ in $\theta_{n,k}$ and $\phi
_{n,k}$ are specialized to $v_{1}=x_{1},v_{2}=x_{2}$ and $v_{3}=1-x_{1}-x_{2}%
$.%
\endproof

\begin{proposition}
The range of $I+RM+MR$ is exactly the subspace $\left\{  p\in\mathbb{P}%
_{n,n-1}^{\perp}\left(  \widehat{T}\right)  :RMp=p\right\}  $.
\end{proposition}

%

\proof
By direct computation $\left(  MR\right)  ^{3}=I$ (cf. Section \ref{SubSecIrr}%
). This implies $\left(  RM\right)  ^{2}=MR$. If $p$ satisfies $RMp=p$ then
$Mp=Rp$ and $p=MRp$. Now suppose $RMp=p$ then $\left(  I+RM+MR\right)
\frac{1}{3}p=p$; hence $p$ is in the range of $I+RM+MR$. Conversely suppose
$p=\left(  I+RM+MR\right)  p^{\prime}$ for some polynomial $p^{\prime}$, then,
$RM\left(  I+RM+MR\right)  p^{\prime}=\left(  RM+\left(  RM\right)
^{2}+I\right)  p^{\prime}=p$.%
\endproof

Let $M_{i,j}^{\left(  n\right)  },R_{i,j}^{\left(  n\right)  }$ denote the
matrix entries of $M,R$ with respect to the basis $\left\{  b_{n,k}:0\leq
k\leq n\right\}  $, respectively (that is $Mb_{n,k}=\sum_{j=0}^{n}%
b_{n,j}M_{j,k}^{\left(  n\right)  }$) . Let $S_{i,j}^{\left(  n\right)  }$
denote the matrix entries of $MR+RM+I$. Then%
\begin{align*}
R_{i,j}^{\left(  n\right)  }  &  =\left(  -1\right)  ^{i}\delta_{i,j}%
;~M_{i,j}^{\left(  n\right)  }=\left(  -1\right)  ^{n}~_{4}F_{3}\left(
\genfrac{}{}{0pt}{}{-i,i+1,-j,j+1}{-n,n+2,1}%
;1\right)  \frac{2i+1}{n+1};\\
S_{i,j}^{\left(  n\right)  }  &  =\left(  \left(  -1\right)  ^{j}+\left(
-1\right)  ^{i}\right)  M_{i,j}^{\left(  n\right)  }+\delta_{i,j}.
\end{align*}
Thus $S_{i,j}^{\left(  n\right)  }=2M_{i,j}^{\left(  n\right)  }+\delta_{i,j}$
if both $i,j$ are even, $S_{i,j}^{\left(  n\right)  }=-2M_{i,j}^{\left(
n\right)  }+\delta_{i,j}$ if both $i,j$ are odd , and $S_{i,j}^{\left(
n\right)  }=0$ if $i-j\equiv1\operatorname{mod}2$.%
\endproof

\begin{corollary}
For $0\leq k\leq\frac{n}{2}$ each polynomial $r_{n,2k}:=2\sum\limits_{0\leq
j\leq n/2}M_{2j,2k}^{\left(  n\right)  }b_{n,2j}+b_{n,2k}$ is totally
symmetric and for $0\leq k\leq\frac{n-1}{2}$ each polynomial $r_{n,2k+1}%
=-2\sum\limits_{0\leq j\leq\left(  n-1\right)  /2}M_{2j+1,2k+1}^{\left(
n\right)  }b_{n,2j+1}+b_{n,2k+1}$ satisfies $Mp=-p=Rp$ (the sign representation).
\end{corollary}

%

\proof
The pattern of zeroes in $\left[  M_{i,j}^{\left(  n\right)  }\right]  $ shows
that $r_{n,2k}=\left(  MR+RM+I\right)  b_{n,2k}\in\operatorname*{span}\left\{
b_{n,2j}\right\}  $ and thus satisfies $Rr_{n,2k}=r_{n,2k}$; combined with
$RMr_{n,2k}=r_{n,2k}$ this shows $r_{n,2k}$ is totally symmetric. A similar
argument applies to $\left(  MR+RM+I\right)  b_{n,2k+1}$.%
\endproof

\begin{theorem}
\label{TheoCD14}The functions $b_{n,k}^{\operatorname*{sym}}$, $0\leq k\leq
d_{\operatorname*{triv}}\left(  n\right)  -1$, as in (\ref{defbpksym}) form a
basis for the totally symmetric polynomials in $\mathbb{P}_{n,n-1}^{\perp
}\left(  \widehat{T}\right)  $.
\end{theorem}

%

\proof
We use the homogeneous form of the $b_{n,m}$ as in \cite{Dunkl_3symm}, that
is, set%
\begin{gather*}
b_{n,2m}^{\prime}\left(  v\right)  =\left(  v_{1}+v_{2}+v_{3}\right)
^{n}b_{n,2m}\left(  \frac{v_{1}}{v_{1}+v_{2}+v_{3}},\frac{v_{2}}{v_{1}%
+v_{2}+v_{3}}\right) \\
=\left(  v_{1}+v_{2}+v_{3}\right)  ^{n-2m}P_{n-2m}^{\left(  0,4m+1\right)
}\left(  \frac{v_{1}+v_{2}-v_{3}}{v_{1}+v_{2}+v_{3}}\right)  \left(
v_{1}+v_{2}\right)  ^{2m}P_{2m}^{\left(  0,0\right)  }\left(  \frac
{v_{1}-v_{2}}{v_{1}+v_{2}}\right)  .
\end{gather*}
Formally $b_{n,j}^{\prime}\left(  v\right)  =\left(  -1\right)  ^{j}\left(
j!\left(  n-j\right)  !\right)  ^{-1}\theta_{n,j}\left(  v\right)  $ with
$\theta_{n,j}$ as in \cite[p.690]{Dunkl_3symm}. The expansion of such a
polynomial is a sum of monomials $v_{1}^{n_{1}}v_{2}^{n_{2}}v_{3}^{n_{3}}$
with $\sum_{i=1}^{3}n_{i}=n$. Symmetrizing the monomial results in the sum of
$v_{1}^{m_{1}}v_{2}^{m_{2}}v_{3}^{m_{3}}$ where $\left(  m_{1},m_{2}%
,m_{3}\right)  $ ranges over all permutations of $\left(  n_{1},n_{2}%
,n_{3}\right)  $. The argument is based on the occurrence of certain indices
in $b_{n,m}$. For a more straightforward approach to the coefficients we use
the following expansions (with $\ell=n-2k,\beta=2k+1$):%
\begin{gather}
\left(  v_{1}+v_{2}+v_{3}\right)  ^{\ell}P_{\ell}^{\left(  0,\beta\right)
}\left(  \frac{v_{1}+v_{2}-v_{3}}{v_{1}+v_{2}+v_{3}}\right)  =\left(
-1\right)  ^{\ell}\left(  v_{1}+v_{2}+v_{3}\right)  ^{\ell}P_{\ell}^{\left(
\beta,0\right)  }\left(  \frac{-v_{1}-v_{2}+v_{3}}{v_{1}+v_{2}+v_{3}}\right)
\label{Pn(0,beta)}\\
=\left(  -1\right)  ^{\ell}\frac{\left(  \beta+1\right)  _{\ell}}{\ell!}%
\sum_{i=0}^{\ell}\frac{\left(  -\ell\right)  _{i}\left(  \ell+\beta+1\right)
_{i}}{i!\left(  \beta+1\right)  _{i}}\left(  v_{1}+v_{2}\right)  ^{i}\left(
v_{1}+v_{2}+v_{3}\right)  ^{\ell-i};\nonumber
\end{gather}
and%
\[
\left(  v_{1}+v_{2}\right)  ^{2k}P_{2k}^{\left(  0,0\right)  }\left(
\frac{v_{1}-v_{2}}{v_{1}+v_{2}}\right)  =\frac{1}{\left(  2k\right)  !}%
\sum_{j=0}^{2k}\frac{\left(  -2k\right)  _{j}\left(  -2k\right)  _{j}\left(
-2k\right)  _{2k-j}}{j!}v_{2}^{j}v_{1}^{2k-j}.
\]
First let $n=2m$. The highest power of $v_{3}$ that can occur in
$b_{2m,2m-2k}^{\prime}$ is $2k$, with corresponding coefficient $\frac{\left(
4m-4k+1\right)  _{2k}}{\left(  2k\right)  !}\sum_{j=0}^{2m-2k}c_{j}v_{2}%
^{j}v_{1}^{2m-j}$ for certain coefficients $\left\{  c_{j}\right\}  $. Recall
that $d_{\operatorname*{triv}}\left(  n\right)  $ is the number of solutions
$\left(  i,j\right)  $ of the equation $3j+2i=2m$ (with $i,j=0,1,2,\ldots$).
The solutions can be listed as $\left(  m,0\right)  ,\left(  m-3,2\right)
,\left(  m-6,4\right)  \ldots\left(  m-3\ell,2\ell\right)  $ where
$\ell=d_{\operatorname*{triv}}\left(  n\right)  -1$. By hypothesis $\left(
m-3k,2k\right)  $ occurs in the list and thus $m-3k\geq0$ and $m-k\geq2k$.
There is only one possible permutation of $v_{1}^{m-k}v_{2}^{m-k}v_{3}^{2k}$
that occurs in $b_{2m,2m-2k}^{\prime}$ and its coefficient is $\frac{\left(
2k-2m\right)  _{m-k}^{3}}{\left(  2m-2k\right)  !}\neq0$. Hence there is a
triangular pattern for the occurrence of $v_{1}^{m}v_{2}^{m}$, $v_{1}%
^{m-1}v_{2}^{m-1}v_{3}^{2}$, $v_{1}^{m-2}v_{2}^{m-2}v_{3}^{4},\ldots$in the
symmetrizations of $b_{2m,2m}^{\prime}$, $b_{2m,2m-2}^{\prime}$, \ldots\ with
nonzero numbers on the diagonal and this proves the basis property when $n=2m$.

Now let $n=2m+1$. The highest power of $v_{3}$ that can occur in
$b_{2m+1,2m-2k}^{\prime}$ is $2k+1$, with coefficient $\frac{\left(
4m-4k+1\right)  _{2k+1}}{\left(  2k+1\right)  !}\sum_{j=0}^{2m-2k}c_{j}%
v_{2}^{j}v_{1}^{2m-j}$ for certain coefficients $\left\{  c_{j}\right\}  $.
The solutions of $3j+2i=2m+1$ can be listed as $\left(  m-1,1\right)  ,\left(
m-4,3\right)  ,\left(  m-7,5\right)  \ldots\left(  m-1-3\ell,2\ell+1\right)  $
where $\ell=d_{\operatorname*{triv}}\left(  n\right)  -1$. By hypothesis
$\left(  m-1-3k,2k+1\right)  $ occurs in this list, thus $m-k\geq2k+1$. There
is only one possible permutation of $v_{1}^{m-k}v_{2}^{m-k}v_{3}^{2k+1}$ that
occurs in $b_{2m+1,2m-2k}^{\prime}$ and its coefficient is $\frac{\left(
2k-2m\right)  _{m-k}^{3}}{\left(  2m-2k\right)  !}\neq0$. As above, there is a
triangular pattern for the occurrence of $v_{1}^{m}v_{2}^{m}v_{3}$,
$v_{1}^{m-1}v_{2}^{m-1}v_{3}^{3}$, $v_{1}^{m-2}v_{2}^{m-2}v_{3}^{5},\ldots$ in
the symmetrizations of $b_{2m+1,2m}^{\prime}$, $b_{2m+1,2m-2}^{\prime}$,
\ldots\ with nonzero numbers on the diagonal and this proves the basis
property when $n=2m+1$.%
\endproof

The totally symmetric orthogonal polynomials can be lifted to a general
triangle $T$.

\begin{definition}
\label{DefBasissymT}Let $T$ denote a triangle. The space of totally symmetric,
orthogonal polynomials of degree $n$ is%
\begin{align}
\mathbb{P}_{n,n-1}^{\perp,\operatorname*{sym}}\left(  T\right)   &  :=\left\{
u\in\mathbb{P}_{n,n-1}^{\perp}\left(  T\right)  :u\text{ has total
symmetry}\right\} \label{defbnkt}\\
&  =\operatorname*{span}\left\{  b_{n,m}^{T,\operatorname*{sym}}:0\leq m\leq
d_{\operatorname*{triv}}\left(  n\right)  -1\right\}  , \label{defBnmT}%
\end{align}
where the lifted symmetric basis functions are given by $b_{n,m}%
^{T,\operatorname*{sym}}:=b_{n,m}^{\operatorname*{sym}}\circ\chi_{T}^{-1}$ for
$b_{n,m}^{\operatorname*{sym}}$ as in Theorem \ref{TheoCD14} and an affine
pullback $\chi_{T}:\widehat{T}\rightarrow T$.
\end{definition}

\subsubsection{A Basis for the $\tau_{\operatorname*{refl}}$ component of
$\mathbb{P}_{n,n-1}^{\perp}\left(  T\right)  $\label{Sectaureflcomp}}

As explained in Section \ref{SubSecIrr} the space $\mathbb{P}_{n,n-1}^{\perp
}\left(  \widehat{T}\right)  $ can be decomposed into the $\tau
_{\operatorname*{triv}}$-, the $\tau_{\operatorname{sign}}$- and the
$\tau_{\operatorname*{refl}}$-component. A basis for the $\tau
_{\operatorname*{triv}}$ component are the fully symmetric basis functions
(cf. Section \ref{ConSymmBas}).

Next, we will construct a basis for all of $\mathbb{P}_{n,n-1}^{\perp}\left(
\widehat{T}\right)  $ by extending the totally symmetric one. It is
straightforward to adjoin the $d_{\operatorname{sign}}\left(  n\right)  $
basis, using the same technique as for the fully symmetric ones: the monomials
which appear in $p$ with $Rp=-p=Mp$ must be permutations of $v_{1}^{n_{1}%
}v_{2}^{n_{2}}v_{3}^{n_{3}}$ with $n_{1}>n_{2}>n_{3}$. As in Theorem
\ref{TheoCD14} for $n=2m$ argue on monomials $v_{1}^{m-k}v_{2}^{m-1-k}%
v_{3}^{2k+1}$ and the polynomials $b_{2m,2m-2k-1}^{\prime}$ with $0\leq k\leq
d_{\operatorname{sign}}\left(  n\right)  -1=d_{\operatorname*{triv}}\left(
n\right)  -2$, and for $n=2m+1$ use the monomials $v_{1}^{m+1-k}v_{2}%
^{m-k}v_{3}^{2k}$ and $b_{2m+1,2m-2k}$ with $0\leq k\leq
d_{\operatorname*{triv}}\left(  n\right)  -1=d_{\operatorname{sign}}\left(
n\right)  -1.$

As we will see when constructing a basis for the non-conforming finite element
space, the $\tau_{\operatorname{sign}}$ component of $\mathbb{P}%
_{n,n-1}^{\perp}\left(  \widehat{T}\right)  $ is not relevant, in contrast to
the $\tau_{\operatorname*{refl}}$ component. In this section, we will
construct a basis for the $\tau_{\operatorname*{refl}}$ polynomials in
$\mathbb{P}_{n,n-1}^{\perp}\left(  \widehat{T}\right)  $. Each such polynomial
is an eigenvector of $RM+MR$ with eigenvalue $-1$. We will show that the
polynomials%
\begin{equation}
b_{n,k}^{\operatorname*{refl}}=\frac{1}{3}\left(  2I-RM-MR\right)
b_{n,2k},~0\leq k\leq\frac{n-1}{3}, \label{defqnk}%
\end{equation}
are linearly independent (and the same as introduced in (\ref{brefl1stdef}))
and, subsequently, that the set%
\begin{equation}
\left\{  RMb_{n,k}^{\operatorname*{refl}},MRb_{n,k}^{\operatorname*{refl}%
}:0\leq k\leq\frac{n-1}{3}\right\}  \label{basisreflrmmr}%
\end{equation}
is a basis for the $\tau_{\operatorname*{refl}}$ subspace of $\mathbb{P}%
_{n,n-1}^{\perp}\left(  \widehat{T}\right)  $. (The upper limit of $k$ is as
in (\ref{basisreflrmmr}) $d_{\operatorname*{refl}}\left(  n\right)  -1$ (cf.
(\ref{defdreflp})).) Note that%
\begin{equation}%
\begin{array}
[c]{c}%
RMb_{n,k}^{\operatorname*{refl}}=\frac{1}{3}\left(  2RM-MR-I\right)
b_{n,2k},\\
MRb_{n,k}^{\operatorname*{refl}}=\frac{1}{3}\left(  2MR-I-RM\right)  b_{n,2k},
\end{array}
\label{defqnkplus}%
\end{equation}
because $\left(  RM\right)  ^{2}=MR$. Thus the calculation of these
polynomials follows directly from the formulae for $\left[  M_{ij}\right]  $
and $\left[  R_{ij}\right]  $. The method of proof relies on complex
coordinates for the triangle.

\begin{lemma}
\label{P00(v1v2)}For $k=0,1,2,\ldots$%
\begin{align*}
P_{2k}^{\left(  0,0\right)  }\left(  s\right)   &  =\left(  -1\right)
^{k}\frac{\left(  k+\frac{1}{2}\right)  _{k}}{k!}\sum_{j=0}^{k}\frac{\left(
-k\right)  _{j}^{2}}{j!\left(  \frac{1}{2}-2k\right)  _{j}}\left(
1-s^{2}\right)  ^{k-j},\\
\left(  v_{1}+v_{2}\right)  ^{2k}P_{2k}^{\left(  0,0\right)  }\left(
\frac{v_{1}-v_{2}}{v_{1}+v_{2}}\right)   &  =\left(  -1\right)  ^{k}%
\frac{\left(  k+\frac{1}{2}\right)  _{k}}{k!}\sum_{j=0}^{k}\frac{\left(
-k\right)  _{j}^{2}}{j!\left(  \frac{1}{2}-2k\right)  _{j}}4^{k-j}\left(
v_{1}v_{2}\right)  ^{k-j}\left(  v_{1}+v_{2}\right)  ^{2j}.
\end{align*}

\end{lemma}%

\proof
Start with the formula (specialized from a formula for Gegenbauer polynomials
\cite[18.5.10]{NIST:DLMF})%
\[
P_{2k}^{\left(  0,0\right)  }\left(  s\right)  =\left(  2s\right)  ^{2k}%
\frac{\left(  \frac{1}{2}\right)  _{2k}}{\left(  2k\right)  !}~_{2}%
F_{1}\left(
\genfrac{}{}{0pt}{}{-k,\frac{1}{2}-k}{\frac{1}{2}-2k}%
;\frac{1}{s^{2}}\right)  .
\]
Apply the transformation (cf. \cite[15.8.1]{NIST:DLMF})%
\[
_{2}F_{1}\left(
\genfrac{}{}{0pt}{}{-k,b}{c}%
;t\right)  =\left(  1-t\right)  ^{k}~_{2}F_{1}\left(
\genfrac{}{}{0pt}{}{-k,c-b}{c}%
;\frac{t}{t-1}\right)
\]
with $t=1/s^{2}$; then $\dfrac{t}{t-1}=\dfrac{1}{1-s^{2}}$ and $s^{2k}\left(
1-\frac{1}{s^{2}}\right)  ^{k}=\left(  -1\right)  ^{k}\left(  1-s^{2}\right)
^{k}$. Also $2^{2k}\frac{\left(  \frac{1}{2}\right)  _{2k}}{\left(  2k\right)
!}=\frac{\left(  \frac{1}{2}\right)  _{2k}}{k!\left(  \frac{1}{2}\right)
_{k}}=\frac{\left(  k+\frac{1}{2}\right)  _{k}}{k!}$. This proves the first
formula. Set $s=\dfrac{v_{1}-v_{2}}{v_{1}+v_{2}}$ then $1-s^{2}=\dfrac
{4v_{1}v_{2}}{\left(  v_{1}+v_{2}\right)  ^{2}}$ to obtain the second one.%
\endproof

Introduce complex homogeneous coordinates:%
\begin{align*}
z  &  =\omega v_{1}+\omega^{2}v_{2}+v_{3}\\
\overline{z}  &  =\omega^{2}v_{1}+\omega v_{2}+v_{3}\\
t  &  =v_{1}+v_{2}+v_{3}.
\end{align*}
Recall $\omega=e^{2\pi\mathrm{\operatorname*{i}}/3}=-\frac{1}{2}%
+\frac{\mathrm{\operatorname*{i}}}{2}\sqrt{3}$ and $\omega^{2}=\overline
{\omega}$. The inverse relations are%
\begin{align*}
v_{1}  &  =\frac{1}{3}\left(  -\left(  \omega+1\right)  z+\omega\overline
{z}+t\right) \\
v_{2}  &  =\frac{1}{3}\left(  \omega z-\left(  \omega+1\right)  \overline
{z}+t\right) \\
v_{3}  &  =\frac{1}{3}\left(  z+\overline{z}+t\right)  .
\end{align*}
Suppose $f\left(  z,\overline{z},t\right)  $ is a polynomial in $z$ and
$\bar{z}$ then $Rf\left(  z,\overline{z},t\right)  =f\left(  \overline
{z},z,t\right)  $ and $Mf\left(  z,\overline{z},t\right)  =f\left(
\omega\overline{z},\omega^{2}z,t\right)  $. Thus $RMf\left(  z,\overline
{z},t\right)  =f\left(  \omega^{2}z,\omega\overline{z},t\right)  $ and
$MRf\left(  z,\overline{z},t\right)  =f\left(  \omega z,\omega^{2}\overline
{z},t\right)  $. The idea is to write $b_{n,2k}$ in terms of $z,\overline
{z},t$ and apply the projection $\Pi:=\frac{1}{3}\left(  2I-MR-RM\right)  $.
To determine linear independence it suffices to consider the terms of highest
degree in $z,\overline{z}$ thus we set $t=v_{1}+v_{2}+v_{3}=0$ in the formula
for $b_{n,2k}$ (previously denoted $b_{n,2k}^{\prime}$ using the homogeneous
coordinates, see proof of Theorem \ref{TheoCD14}). From formula
(\ref{Pn(0,beta)}) and Lemma \ref{P00(v1v2)}%
\begin{align*}
b_{n,2k}^{\prime}\left(  v_{1},v_{2},0\right)   &  =\left(  n-2k+2\right)
_{n-2k}\left(  v_{1}+v_{2}\right)  ^{n-2k}\left(  -1\right)  ^{k}\frac{\left(
k+\frac{1}{2}\right)  _{k}}{k!}\\
&  \times\sum_{j=0}^{k}\frac{\left(  -k\right)  _{j}^{2}}{j!\left(  \frac
{1}{2}-2k\right)  _{j}}4^{k-j}\left(  v_{1}v_{2}\right)  ^{k-j}\left(
v_{1}+v_{2}\right)  ^{2j}.
\end{align*}
The coefficient of $\left(  v_{1}v_{2}\right)  ^{k}\left(  v_{1}+v_{2}\right)
^{n-2k}$ in $b_{n,2k}^{\prime}\left(  v_{1},v_{2},0\right)  $ is nonzero, and
this is the term with highest power of $v_{1}v_{2}$. Thus $\left\{
b_{n,2k}^{\prime}\left(  v_{1},v_{2},0\right)  :0\leq k\leq\frac{n-2}%
{3}\right\}  $ is a basis for $\mathrm{\operatorname*{span}}\left\{  \left(
v_{1}v_{2}\right)  ^{k}\left(  v_{1}+v_{2}\right)  ^{n-2k}:0\leq k\leq
\frac{n-2}{3}\right\}  $. The next step is to show that the projection $\Pi$
has trivial kernel. In the complex coordinates $v_{1}+v_{2}=-\frac{1}%
{3}\left(  z+z-t\right)  =-\frac{1}{3}\left(  z+z\right)  $ and $v_{1}%
v_{2}=\frac{1}{9}\left(  z^{2}-z\overline{z}+\overline{z}^{2}\right)  $
(discarding terms of lower order in $z,\overline{z}$, that is, set $t=0$).

\begin{proposition}
\label{Prop17}If $\Pi\sum_{k=0}^{\left\lfloor \left(  n-1\right)
/3\right\rfloor }c_{k}\left(  z+\overline{z}\right)  ^{n-2k}\left(
z^{2}-z\overline{z}+\overline{z}^{2}\right)  ^{k}=0$ then $c_{k}=0$ for all
$k$.
\end{proposition}

%

\proof
For any polynomial $f\left(  z,\overline{z}\right)  $ we have $\Pi f\left(
z,\overline{z}\right)  =\frac{1}{3}\left(  2f\left(  z,\overline{z}\right)
-f\left(  \omega^{2}z,\omega\overline{z}\right)  -f\left(  \omega z,\omega
^{2}\overline{z}\right)  \right)  $. In particular%
\begin{gather*}
\Pi\left(  z+\overline{z}\right)  ^{n-2k}\left(  z^{2}-z\overline{z}%
+\overline{z}^{2}\right)  ^{k}=\Pi\left(  z+\overline{z}\right)
^{n-3k}\left(  z^{3}+\overline{z}^{3}\right)  ^{k}\\
=\frac{1}{3}\left\{  2\left(  z+\overline{z}\right)  ^{n-3k}-\left(
\omega^{2}z+\omega\overline{z}\right)  ^{n-3k}-\left(  \omega z+\omega
^{2}\overline{z}\right)  ^{n-3k}\right\}  \left(  z^{3}+\overline{z}%
^{3}\right)  ^{k}.
\end{gather*}
By hypothesis $n-3k\geq1$. Evaluate the expression at $z=e^{\pi
\mathrm{\operatorname*{i}/}6}+\varepsilon$ where $\varepsilon$ is real and
near $0$. Note $e^{\pi\mathrm{\operatorname*{i}/}6}=\frac{1}{2}\left(
\sqrt{3}+\mathrm{\operatorname*{i}}\right)  $. Then%
\begin{align*}
z+\overline{z}  &  =\sqrt{3}+2\varepsilon,\\
\omega^{2}z+\omega\overline{z}  &  =-\varepsilon,\\
\omega z+\omega^{2}\overline{z}  &  =-\sqrt{3}-\varepsilon,\\
z^{3}+\overline{z}^{3}  &  =3\varepsilon+3\sqrt{3}\varepsilon^{2}%
+2\varepsilon^{3},
\end{align*}
and%
\begin{align*}
&  \frac{1}{3}\left\{  2\left(  z+\overline{z}\right)  ^{n-3k}-\left(
\omega^{2}z+\omega\overline{z}\right)  ^{n-3k}-\left(  \omega z+\omega
^{2}\overline{z}\right)  ^{n-3k}\right\}  \left(  z^{3}+\overline{z}%
^{3}\right)  ^{k}\\
&  =\frac{1}{3}\left\{  \left(  2-\left(  -1\right)  ^{n-3k}\right)
\times3^{\left(  n-3k\right)  /2}-\left(  -\varepsilon\right)  ^{n-3k}%
+C\varepsilon+O\left(  \varepsilon^{2}\right)  \right\}  \varepsilon
^{k}\left(  3+3\sqrt{3}\varepsilon+2\varepsilon^{2}\right)  ^{k},
\end{align*}
where $C=3^{\left(  n--3k-1\right)  /2}\left(  n-3k\right)  \left(  4-2\left(
-1\right)  ^{n-3k}\right)  $ (binomial theorem). The dominant term in the
right-hand side is $\left(  2-\left(  -1\right)  ^{n-3k}\right)  3^{\left(
n-k\right)  /2-1}\varepsilon^{k}$. Now suppose $\Pi\sum\limits_{k=0}%
^{\left\lfloor \left(  n-1\right)  /3\right\rfloor }c_{k}\left(
z+\overline{z}\right)  ^{n-2k}\left(  z^{2}-z\overline{z}+\overline{z}%
^{2}\right)  ^{k}=0$. Evaluate the polynomial at $z=e^{\pi
\mathrm{\operatorname*{i}/}6}+\varepsilon$. Let $\varepsilon\rightarrow0$
implying $c_{0}=0$. Indeed write the expression as%
\[
\sum_{k=0}^{\left\lfloor \left(  n-1\right)  /3\right\rfloor }c_{k}\left(
2-\left(  -1\right)  ^{n-3k}\right)  3^{\left(  n-k\right)  /2-1}%
\varepsilon^{k}\left(  1+O\left(  \varepsilon\right)  \right)  =0.
\]
Since $2-\left(  -1\right)  ^{n-3k}\geq1$ this shows $c_{k}=0$ for all $k$.%
\endproof

We have shown:

\begin{proposition}
Suppose $\Pi\sum\limits_{k=0}^{\left\lfloor \left(  n-1\right)
/3\right\rfloor }c_{k}b_{n,2k}=0$ then $c_{k}=0$ for all $k$; the cardinality
of the set (\ref{basisreflrmmr}) is $d_{\operatorname*{refl}}\left(  n\right)
$.
\end{proposition}

\begin{theorem}
\label{TheoCDBasis}\

\begin{enumerate}
\item[a.] The polynomials $\left\{  \Pi b_{n,2k}:0\leq k\leq\frac{n-1}%
{3}\right\}  $ are linearly independent.

\item[b.] The set $\left\{  RM\Pi b_{n,2k},MR\Pi b_{n,2k}:0\leq k\leq
\frac{n-1}{3}\right\}  $ is linearly independent and defines a basis for the
$\tau_{\operatorname*{refl}}$ component of $\mathbb{P}_{n,n-1}^{\perp}\left(
\widehat{T}\right)  $.
\end{enumerate}
\end{theorem}

%

\proof
In general $\Pi z^{a}\overline{z}^{b}=z^{a}\overline{z}^{b}$ if $a-b\equiv
1,2\operatorname{mod}3$ and $\Pi z^{a}\overline{z}^{b}=0$ if $a-b\equiv
0\operatorname{mod}3$. Expand the polynomials $w_{k}\left(  z,\overline
{z}\right)  :=\Pi\left(  z+\overline{z}\right)  ^{n-3k}\left(  z^{3}%
+\overline{z}^{3}\right)  ^{k}$ by the binomial theorem to obtain%
\[
\Pi\left(  z+\overline{z}\right)  ^{n-3k}\left(  z^{3}+\overline{z}%
^{3}\right)  ^{k}=\sum_{\substack{j=0\\n-2j\equiv1,2\operatorname{mod}%
3}}^{n-3k}\binom{n-3k}{j}z^{n-3k-j}\overline{z}^{j}\left(  z^{3}+\overline
{z}^{3}\right)  ^{k}.
\]
Then%
\begin{align*}
RMw_{k}\left(  z,\overline{z}\right)   &  =\sum_{j=0,n-2j\equiv
1,2\operatorname{mod}3}^{n-3k}\binom{n-3k}{j}\omega^{2j-n}z^{n-3k-j}%
\overline{z}^{j}\left(  z^{3}+\overline{z}^{3}\right)  ^{k},\\
MRw_{k}\left(  z,\overline{z}\right)   &  =\sum_{j=0,n-2j\equiv
1,2\operatorname{mod}3}^{n-3k}\binom{n-3k}{j}\omega^{n-2j}z^{n-3k-j}%
\overline{z}^{j}\left(  z^{3}+\overline{z}^{3}\right)  ^{k}.
\end{align*}
Firstly we show that $\left\{  RMw_{k},MRw_{k}\right\}  $ is linearly
independent for $0\leq k\leq\frac{n-1}{3}$. For each value of
$n\operatorname{mod}3$ we select the highest degree terms from $RMw_{k}$ and
$MRw_{k}$: (i) $n=3m+1$, $\omega^{2}z^{3m+1}+\omega\overline{z}^{3m+1}$ and
$\omega z^{3m+1}+\omega^{2}\overline{z}^{3m+1}$, (ii) $n=3m+2,$ $\omega
z^{3m+2}+\omega^{2}\overline{z}^{3m+2}$ and $\omega^{2}z^{3m+2}+\omega
\overline{z}^{3m+2}$, (iii) $n=3m$, $\left(  n-3k\right)  \left(  \omega
^{2}z^{3m}\overline{z}+\omega z\overline{z}^{3m}\right)  $ and $\left(
n-3k\right)  \left(  \omega z^{3m}\overline{z}+\omega^{2}z\overline{z}%
^{3m}\right)  $ (by hypothesis $n-3k\geq1$). In each case the two terms are
linearly independent (the determinant of the coefficients is $\pm\left(
\omega-\omega^{2}\right)  =\mp\mathrm{\operatorname*{i}}\sqrt{3}$). Secondly
the same argument as in the previous theorem shows that $\sum_{k=0}%
^{\left\lfloor \left(  n-1\right)  /3\right\rfloor }\left\{  c_{k}%
RMw_{k}+d_{k}MRw_{k}\right\}  =0$ implies $c_{k}RMw_{k}+d_{k}MRw_{k}=0$ for
all $k$. By the first part it follows that $c_{k}=0=d_{k}$. This completes the
proof.%
\endproof

\begin{remark}
The basis $b_{n,k}$ for $\mathbb{P}_{n,n-1}^{\perp}\left(  \widehat{T}\right)
$ in (\ref{defbnk}) is mirror symmetric with respect to the angular bisector
in $\widehat{T}$ through the origin for even $k$ and is mirror skew-symmetric
for odd $k$. This fact makes the point $\mathbf{0}$ in $\widehat{T}$ special
compared to the other vertices. As a consequence the functions defined in
Theorem \ref{TheoCDBasis}.a reflects the special role of $\mathbf{0}$. Part b
shows that it is possible to define a basis with functions which are either
symmetric with respect to the angle bisector in $\widehat{T}$ through $\left(
1,0\right)  ^{\intercal}$ or through $\left(  0,1\right)  ^{\intercal}$ by
\textquotedblleft rotating\textquotedblright\ the functions $\Pi b_{n,2k}$ to
these vertices:%
\[
RM\left(  \Pi b_{n,2k}\right)  \left(  x_{1},x_{2}\right)  =\left(  \Pi
b_{n,2k}\right)  \left(  x_{2},1-x_{1}-x_{2}\right)  \quad\text{and\quad
}MR\left(  \Pi b_{n,2k}\right)  \left(  x_{1},x_{2}\right)  =\left(  \Pi
b_{n,2k}\right)  \left(  1-x_{1}-x_{2},x_{1}\right)  .
\]
Since the dimension of $E^{\left(  \operatorname*{refl}\right)  }$ is
$2d_{\operatorname*{refl}}\left(  n\right)  =2\left\lfloor \frac{n+2}%
{3}\right\rfloor $ is not (always) a multiple of $3$, it is, in general, not
possible to define a basis where all three vertices of the triangle are
treated in a symmetric way.
\end{remark}

\begin{definition}
\label{RemMirror}Let%
\begin{equation}
\mathbb{P}_{n,n-1}^{\perp,\operatorname*{refl}}\left(  \widehat{T}\right)
:=\operatorname*{span}\left\{  RM\Pi b_{n,2k},MR\Pi b_{n,2k}:0\leq k\leq
\frac{n-1}{3}\right\}  . \label{defPnn-1refl}%
\end{equation}
This space is lifted to a general triangle $T$ by fixing a vertex $\mathbf{P}$
of $T$ and setting%
\begin{equation}
\mathbb{P}_{n,n-1}^{\perp,\operatorname*{refl}}\left(  T\right)  :=\left\{
u\circ\chi_{\mathbf{P},T}^{-1}:u\in\mathbb{P}_{n,n-1}^{\perp
,\operatorname*{refl}}\left(  \widehat{T}\right)  \right\}  ,
\label{defReflspace}%
\end{equation}
where the lifting $\chi_{\mathbf{P},T}$ is an affine pullback $\chi
_{\mathbf{P},T}:\widehat{T}\rightarrow T$ which maps $\mathbf{0}$ to
$\mathbf{P}$.

The basis $b_{n,k}^{\operatorname*{refl}}$ to describe the restrictions of
facet-oriented, non-conforming finite element functions to the facets is
related to a reduced space and defined as in (\ref{defqnk}) with lifted
versions%
\begin{equation}
b_{n,k}^{\mathbf{P},T}:=b_{n,k}^{\operatorname*{refl}}\circ\chi_{\mathbf{P}%
,T}^{-1},\quad0\leq k\leq\frac{n-1}{3}. \label{defbmitttriangle}%
\end{equation}

\end{definition}

\begin{remark}
The construction of the spaces $\mathbb{P}_{p,p-1}^{\perp,\operatorname*{sym}%
}\left(  T\right)  \ $and $\mathbb{P}_{p,p-1}^{\perp,\operatorname*{refl}%
}\left(  T\right)  $ (cf. Definitions \ref{DefBasissymT} and \ref{RemMirror})
implies the direct sum decomposition%
\begin{equation}
\operatorname*{span}\left\{  b_{p,2k}\circ\chi_{\mathbf{P},T}^{-1}:0\leq
k\leq\left\lfloor p/2\right\rfloor \right\}  =\mathbb{P}_{p,p-1}%
^{\perp,\operatorname*{sym}}\left(  T\right)  \oplus\mathbb{P}_{p,p-1}%
^{\perp,\operatorname*{refl}}\left(  T\right)  . \label{defmirrorP}%
\end{equation}
It is easy to verify that the basis functions $b_{p,k}^{\mathbf{P},T}$ are
mirror symmetric with respect to the angle bisector in $T$ through
$\mathbf{P}$. However, the space $\mathbb{P}_{n,n-1}^{\perp
,\operatorname*{refl}}\left(  T\right)  $ is independent of the choice of the
vertex $\mathbf{P}$.

In Appendix \ref{AA} we will define further sets of basis functions for the
$\tau_{\operatorname*{refl}}$ component of $\mathbb{P}_{n,n-1}^{\perp}\left(
\widehat{T}\right)  $ -- different choices might be preferable for different
kinds of applications.
\end{remark}

\subsection{Simplex-Supported and Facet-Oriented Non-Conforming Basis
Functions}

In this section, we will define non-conforming Crouzeix-Raviart type functions
which are supported either on one single tetrahedron or on two tetrahedrons
which share a common facet. As a prerequisite, we study in
\S \ref{SecOrthPolyTriStar} piecewise orthogonal polynomials on triangle
stars, i.e., on a collection of triangles which share a common vertex and
cover a neighborhood of this vertex (see Notation \ref{Notstar}). We will
derive conditions such that these functions are continuous across common edges
and determine the dimension of the resulting space. This allows us to
determine the non-conforming Courzeix-Raviart basis functions which are either
supported on a single tetrahedron (see \S \ref{SecBasSKnc}) or on two adjacent
tetrahedrons (see \S \ref{SecBasSTnc}) by \textquotedblleft
closing\textquotedblright\ triangle stars either by a single triangle or
another triangle star.

\subsubsection{Orthogonal Polynomials on Triangle
Stars\label{SecOrthPolyTriStar}}

The construction of the functions $B_{p,k}^{K,\operatorname*{nc}}$ and
$B_{p,k}^{T,\operatorname*{nc}}$ as in (\ref{DefSkpsym}) and
(\ref{defedgesupp}) requires some results of continuous, piecewise orthogonal
polynomials on \textit{triangle stars} which we provide in this section.

\begin{notation}
\label{Notstar}A subset $C\subset\Omega$ is a \emph{triangle star} if $C$ is
the union of some, say $m_{C}\geq3$, triangles $T\in\mathcal{F}_{C}%
\subset\mathcal{F}$, i.e., $C=%
{\displaystyle\bigcup\limits_{T\in\mathcal{F}_{C}}}
T$ and there exists some vertex $\mathbf{V}_{C}\in\mathcal{V}$ such that%
\begin{equation}%
\begin{array}
[c]{ll}%
\mathbf{V}_{C}\text{ is a vertex of }T & \forall T\in\mathcal{F}_{C},\\
\exists\text{ a continuous, piecewise affine mapping }\chi:D_{m_{C}%
}\rightarrow C & \text{such that }\chi\left(  0\right)  =\mathbf{V}_{C}.
\end{array}
\label{deftrianglestar}%
\end{equation}
Here, $D_{k}$ denotes the regular closed $k$-gon (in $\mathbb{R}^{2}$).
\end{notation}

For a triangle star $C$, we define%
\[
\mathbb{P}_{p,p-1}^{\perp}\left(  C\right)  :=\left\{  u\in C^{0}\left(
C\right)  \mid\forall T\in\mathcal{F}_{C}:\left.  u\right\vert _{T}%
\in\mathbb{P}_{p,p-1}^{\perp}\left(  T\right)  \right\}  .
\]
In the next step, we will explicitly characterize the space $\mathbb{P}%
_{p,p-1}^{\perp}\left(  C\right)  $ by defining a set of basis functions. Set
$\mathbf{A}:=\mathbf{V}_{C}$ (cf. (\ref{deftrianglestar})) and pick an outer
vertex in $\mathcal{F}_{C}$, denote it by $\mathbf{A}_{1}$, and number the
remaining vertices $\mathbf{A}_{2},\ldots,\mathbf{A}_{m_{C}}$ in
$\mathcal{F}_{C}$ counterclockwise. We use the cyclic numbering convention
$\mathbf{A}_{m_{C}+1}:=\mathbf{A}_{1}$ and also for similar quantities.

For $1\leq\ell\leq m_{C}$, let $\mathbf{e}_{\ell}:=\left[  \mathbf{A}%
,\mathbf{A}_{\ell}\right]  $ be the straight line (convex hull) between and
including $\mathbf{A}$, $\mathbf{A}_{\ell}$. Let $T_{\ell}\in\mathcal{F}_{C}$
be the triangle with vertices $\mathbf{A}$, $\mathbf{A}_{\ell}$,
$\mathbf{A}_{\ell+1}$. Then we choose the affine pullbacks to the reference
element $\widehat{T}$ by%
\[
\chi_{\ell}\left(  x_{1},x_{2}\right)  :=\left\{
\begin{array}
[c]{ll}%
\mathbf{A}+x_{1}\left(  \mathbf{A}_{\ell}-\mathbf{A}\right)  +x_{2}\left(
\mathbf{A}_{\ell+1}-\mathbf{A}\right)  & \text{if }\ell\text{ is odd,}\\
\mathbf{A}+x_{1}\left(  \mathbf{A}_{\ell+1}-\mathbf{A}\right)  +x_{2}\left(
\mathbf{A}_{\ell}-\mathbf{A}\right)  & \text{if }\ell\text{ is even.}%
\end{array}
\right.
\]
In this way, the common edges $\mathbf{e}_{\ell}$ are parametrized by
$\chi_{\ell-1}\left(  t,0\right)  =\chi_{\ell}\left(  t,0\right)  $ if
$3\leq\ell\leq m_{C}$ is odd and by $\chi_{\ell-1}\left(  0,t\right)
=\chi_{\ell}\left(  0,t\right)  $ if $2\leq\ell\leq m_{C}$ is even. The final
edge $\mathbf{e}_{1}$ is parametrized by $\chi_{1}\left(  t,0\right)
=\chi_{m_{C}}\left(  t,0\right)  $ if $m_{C}$ is even and by $\chi_{1}\left(
t,0\right)  =\chi_{m_{C}}\left(  0,t\right)  $ (with interchanged arguments!)
otherwise. We introduce the set%
\[
R_{p,C}:=\left\{
\begin{array}
[c]{ll}%
\left\{  0,\ldots,p\right\}  & \text{if }m_{C}\text{ is even,}\\
\left\{  2\ell:0\leq\ell\leq\left\lfloor \frac{p}{2}\right\rfloor \right\}  &
\text{if }m_{C}\text{ is odd}%
\end{array}
\right.
\]
and define the functions (cf. (\ref{defbnkt}), (\ref{defReflspace}),
(\ref{defmirrorP}))%
\begin{equation}
\left.  b_{p,k}^{C}\right\vert _{T_{\ell}}:=b_{p,k}\circ\chi_{\ell}%
^{-1},\qquad\forall k\in R_{p,C}. \label{basispatch}%
\end{equation}

\begin{lemma}
\label{Lemcont}For a triangle star $C$, a basis for $\mathbb{P}_{p,p-1}%
^{\perp}\left(  C\right)  $ is given by $b_{p,k}^{C}$, $k\in R_{p,C}$. Further%
\begin{equation}
\dim\mathbb{P}_{p,p-1}^{\perp}\left(  C\right)  =\left\{
\begin{array}
[c]{cc}%
p+1 & \text{if }m_{C}\text{ is even,}\\
\left\lfloor \frac{p}{2}\right\rfloor +1 & \text{if }m_{C}\text{ is odd.}%
\end{array}
\right.  \label{dimformula}%
\end{equation}

\end{lemma}%

\proof
We show that $\left(  b_{p,k}^{C}\right)  _{k\in R_{p,C}}$ is a basis of
$\mathbb{P}_{p,p-1}^{\perp}\left(  C\right)  $ and the dimension formula.

\textbf{Continuity across }$\mathbf{e}_{\ell}$ \textbf{for odd }$3\leq\ell\leq
m_{C}$.

The definition of the lifted orthogonal polynomials (see (\ref{defbnkt}),
(\ref{defReflspace}), (\ref{defmirrorP})) implies that the continuity across
$\mathbf{e}_{\ell}$ for odd $3\leq\ell\leq m_{C}$ is equivalent to%
\[
\sum_{k=0}^{p}\alpha_{p,k}^{\left(  \ell-1\right)  }b_{p,k}^{\operatorname{I}%
}=\sum_{k=0}^{p}\alpha_{p,k}^{\left(  \ell\right)  }b_{p,k}^{\operatorname{I}%
}.
\]
From Lemma \ref{LemRestLinIndep} we conclude that the continuity across such
edges is equivalent to
\begin{equation}
\alpha_{p,k}^{\left(  \ell-1\right)  }=\alpha_{p,k}^{\left(  \ell\right)
}\qquad\forall0\leq k\leq p.
\end{equation}

\textbf{Continuity across }$\mathbf{e}_{\ell}$ \textbf{for even }$2\leq
\ell\leq m_{C}$.

Note that $\chi_{2}\left(  0,t\right)  =\chi_{3}\left(  0,t\right)  $. Taking
into account (\ref{defbnkt}), (\ref{defReflspace}), (\ref{defmirrorP}) we see
that the continuity across $\mathbf{e}_{\ell}$ is equivalent to%
\[
\sum_{k=0}^{p}\alpha_{p,k}^{\left(  2\right)  }b_{p,k}^{\operatorname{II}%
}=\sum_{k=0}^{p}\alpha_{p,k}^{\left(  3\right)  }b_{p,k}^{\operatorname{II}}.
\]

From Lemma \ref{LemRestLinIndep} we conclude that the continuity across
$\mathbf{e}_{\ell}$ for even $2\leq\ell\leq m_{C}$ is again equivalent to
\begin{equation}
\alpha_{p,k}^{\left(  \ell-1\right)  }=\alpha_{p,k}^{\left(  \ell\right)
}\qquad\forall0\leq k\leq p.
\end{equation}

\textbf{Continuity across }$\mathbf{e}_{1}$

For even $m_{C}$ the previous argument also applies for the edge
$\mathbf{e}_{1}$ and the functions $b_{p,k}^{C}$, $0\leq k\leq p$, are
continuous across $\mathbf{e}_{1}$. For odd $m_{C}$, note that $\chi
_{1}\left(  t,0\right)  =\chi_{m_{C}}\left(  0,t\right)  $. Taking into
account (\ref{defbnkt}), (\ref{defReflspace}), (\ref{defmirrorP}) we see that
the continuity across $\mathbf{e}_{1}$ is equivalent to
\[
\sum_{k=0}^{p}\alpha_{p,k}^{\left(  1\right)  }b_{p,k}^{\operatorname{I}}%
=\sum_{k=0}^{p}\alpha_{p,k}^{\left(  m_{C}\right)  }b_{p,k}^{\operatorname{II}%
}.
\]
Using the symmetry relation (\ref{symmbnk}) we conclude that this is
equivalent to%
\[
\sum_{k=0}^{p}\alpha_{p,k}^{\left(  1\right)  }b_{p,k}^{\operatorname{I}}%
=\sum_{k=0}^{p}\alpha_{p,k}^{\left(  m_{C}\right)  }\left(  -1\right)
^{k}b_{p,k}^{\operatorname{I}}.
\]
From Lemma \ref{LemRestLinIndep} we conclude that this, in turn, is equivalent
to%
\begin{equation}%
\begin{array}
[c]{cl}%
\alpha_{p,k}^{\left(  1\right)  }=\alpha_{p,k}^{\left(  m_{C}\right)  } &
k\text{ is even,}\\
\alpha_{p,k}^{\left(  1\right)  }=-\alpha_{p,k}^{\left(  m_{C}\right)  } &
k\text{ is odd.}%
\end{array}
\end{equation}
From the above reasoning, the continuity of $b_{p,k}^{C}$ across
$\mathbf{e}_{1}$ follows if $\alpha_{n,k}^{\left(  \ell\right)  }=0$ for odd
$k$ and all $1\leq\ell\leq m_{C}$.

The proof of the dimension formula (\ref{dimformula}) is trivial.%
\endproof

\subsubsection{A Basis for the Symmetric Non-Conforming Space
$S_{K,\operatorname*{nc}}^{p}$\label{SecBasSKnc}}

In this section, we will prove that $S_{K,\operatorname*{nc}}^{p}$ (cf.
(\ref{DefSkpsym})) satisfies%
\[
S_{K,\operatorname*{nc}}^{p}\oplus S_{K,\operatorname*{c}}^{p}=S_{K}%
^{p}:=\left\{  u\in S_{\mathcal{G}}^{p}:\operatorname*{supp}u\subset
K\right\}  ,
\]
where $S_{\mathcal{G}}^{p}$ is defined in (\ref{hpfinele}) and, moreover, that
the functions $B_{p,k}^{K,\operatorname*{nc}}$, $k=0,1,\ldots
,d_{\operatorname*{triv}}\left(  p\right)  -1$, as in (\ref{DefBpkKhut}),
(\ref{DefSkpsym}) form a basis of $S_{K,\operatorname*{nc}}^{p}$.

Let $T$ denote one facet of $K$ and let $C:=\partial K\backslash\overset
{\circ}{T}$. Since $C$ is a triangle star with $m_{C}=3$, we can apply Lemma
\ref{Lemcont} to obtain that
\[
\left.  S_{K}^{p}\right\vert _{C}:=\left\{  \left.  u\right\vert _{C}:u\in
S_{K}^{p}\right\}  \subset\operatorname*{span}\left\{  b_{p,2k}^{C}:0\leq
k\leq\left\lfloor \frac{p}{2}\right\rfloor \right\}  .
\]
The continuity of $b_{p,2k}^{C}$ implies that the restriction $b_{p,2k}%
^{\partial T}:=\left.  b_{p,2k}^{C}\right\vert _{\partial T}$ is continuous.
From (\ref{bnkIII}) we conclude that
\begin{equation}
\left.  b_{p,2k}^{\partial T}\right\vert _{E}=P_{2k}^{E}\qquad\forall
E\subset\partial T, \label{defbn2kdT}%
\end{equation}
where $P_{2k}^{E}$ is the Legendre polynomial of even degree $2k$ scaled to
the edge $E$ with endpoint values $+1$ and symmetry with respect to the
midpoint of $E$. Hence, we are looking for orthogonal polynomials
$\mathbb{P}_{p,p-1}^{\perp}\left(  T\right)  $ whose traces on $\partial T$
are linear combination of $b_{p,2k}^{\partial T}$, $0\leq k\leq\left\lfloor
\frac{p}{2}\right\rfloor $. From (\ref{symmbnk}) we deduce that they have
total symmetry, i.e., belong to the space $\mathbb{P}_{p,p-1}^{\perp
,\operatorname*{sym}}\left(  T\right)  $ (cf. Definition \ref{DefBasissymT}).
For $0\leq m\leq d_{\operatorname*{triv}}\left(  p\right)  -1$, let
$b_{p,m}^{\partial K,\operatorname*{sym}}:\partial K\rightarrow\mathbb{R}$ be
defined facet-wise for any $T\subset\partial K$ by%
\begin{equation}
\left.  b_{p,m}^{\partial K,\operatorname*{sym}}\right\vert _{T}%
:=b_{p,m}^{T,\operatorname*{sym}}\qquad0\leq m\leq d_{\operatorname*{triv}%
}\left(  p\right)  -1. \label{defsurfacesym}%
\end{equation}
Finally, we extend the function $b_{p,m}^{\partial K,\operatorname*{sym}}$ to
the total simplex $K$ by polynomial extension (cf. (\ref{DefBpkKhut}),
(\ref{18b}))%
\begin{equation}
B_{p,m}^{K,\operatorname*{nc}}=\sum_{\mathbf{N}\in\mathcal{N}^{p}\cap\partial
K}b_{p,m}^{\partial K,\operatorname*{sym}}\left(  \mathbf{N}\right)  \left.
B_{p,\mathbf{N}}^{\mathcal{G}}\right\vert _{K}\qquad0\leq m\leq
d_{\operatorname*{triv}}\left(  p\right)  -1. \label{defsimplsupp}%
\end{equation}
These functions are the same as those introduced in Definition
\ref{DefSymSpaceK}. The above reasoning leads to the following Proposition.

\begin{proposition}
For a simplex $K$, the space of non-conforming, simplex-supported
Crouzeix-Raviart finite elements can be chosen as in (\ref{DefSkpsym}) and the
functions $B_{p,k}^{K,\operatorname*{nc}}$, $0\leq k\leq
d_{\operatorname*{triv}}\left(  p\right)  -1$ are linearly independent.
\end{proposition}

\subsubsection{A Basis for $S_{T,\operatorname*{nc}}^{p}$\label{SecBasSTnc}}

Let $T\in\mathcal{F}_{\Omega}$ be an inner facet and $K_{1},K_{2}%
\in\mathcal{G}$ such that $T=K_{1}\cap K_{2}$ and $\omega_{T}=K_{1}\cup K_{2}$
(cf. (\ref{defoftrianglesubsets})) with the convention that the unit normal
$\mathbf{n}_{T}$ points into $K_{2}$. In this section, we will prove that a
space $\tilde{S}_{T,\operatorname*{nc}}^{p}$ which satisfies%
\begin{equation}
\tilde{S}_{T,\operatorname*{nc}}^{p}\oplus\left(
{\displaystyle\bigoplus\limits_{i=1}^{2}}
S_{K_{i},\operatorname*{nc}}^{p}\right)  \oplus\left(
{\displaystyle\bigoplus\limits_{i=1}^{2}}
S_{K_{i},\operatorname*{c}}^{p}\right)  \oplus S_{T,\operatorname*{c}}%
^{p}=S_{T}^{p}:=\left\{  u\in S_{\mathcal{G}}^{p}:\operatorname*{supp}%
u\subset\omega_{T}\right\}  \label{directsplitfacetspace}%
\end{equation}
can be chosen as $\tilde{S}_{T,\operatorname*{nc}}^{p}%
:=S_{T,\operatorname*{nc}}^{p}$ (cf. (\ref{STpncdef})) and, moreover, that the
functions $B_{p,k}^{T,\operatorname*{nc}}$, $k=0,1,\ldots
,d_{\operatorname*{refl}}\left(  p\right)  -1$, as in (\ref{defedgesupp}) form
a basis of $S_{T,\operatorname*{nc}}^{p}$.

Let $C_{i}:=\left(  \partial K_{i}\right)  \backslash\overset{\circ}{T}$,
$i=1,2$, denote the triangle star (cf. Notation \ref{Notstar}) formed by the
three remaining triangles of $\partial K_{i}$. We conclude from Lemma
\ref{Lemcont} that a basis for $\mathbb{P}_{p,p-1}^{\perp}\left(
C_{i}\right)  $ is given by $b_{p,2k}^{C_{i}}$, $0\leq k\leq\left\lfloor
\frac{p}{2}\right\rfloor $ (cf. (\ref{basispatch})). Any function $u$ in
$S_{T}^{p}$ satisfies%
\begin{equation}%
\begin{array}
[c]{cl}%
\gamma_{K_{i}}u\in\mathbb{P}_{p}\left(  K_{i}\right)  & i=1,2,\\
\left.  \left(  \gamma_{K_{i}}u\right)  \right\vert _{T^{\prime}}\in
\mathbb{P}_{p,p-1}^{\perp}\left(  T^{\prime}\right)  & \forall T^{\prime
}\subset C_{i},\quad i=1,2,\\
\left[  u\right]  _{T}\in\mathbb{P}_{p,p-1}^{\perp}\left(  T\right)  . &
\end{array}
\label{jumpfacet}%
\end{equation}
Since any function in $S_{T}^{p}$ is continuous on $C_{i}$, we conclude from
Lemma \ref{Lemcont} (with $m_{C_{i}}=3$) that%
\begin{equation}
\left.  u\right\vert _{C_{i}}\in\mathbb{P}_{p,p-1}^{\perp}\left(
C_{i}\right)  \text{\quad and\quad}\left.  \gamma_{K_{i}}u\right\vert
_{\partial T}\in\operatorname*{span}\left\{  b_{p,2k}^{\partial T}:0\leq
k\leq\left\lfloor \frac{p}{2}\right\rfloor \right\}  \qquad i=1,2
\label{twosidedlinearcombis}%
\end{equation}
with $b_{p,2k}^{\partial T}$ as in (\ref{defbn2kdT}).

To identify a space $\tilde{S}_{T,\operatorname*{nc}}^{p}$ which satisfies
(\ref{directsplitfacetspace}) we consider the jump condition in
(\ref{jumpfacet}) restricted to the boundary $\partial T$. The symmetry of the
functions $b_{p,2k}^{\partial T}$ implies that $\left[  u\right]  _{T}%
\in\mathbb{P}_{p,p-1}^{\perp,\operatorname*{sym}}\left(  T\right)  $, i.e.,
there is a function $q_{1}\in S_{K_{1},\operatorname*{nc}}^{p}$(see
(\ref{DefSkpsym})) such that $\left[  u\right]  _{T}=\left.  q_{1}\right\vert
_{T}$ and $\tilde{u}$, defined by $\left.  \tilde{u}\right\vert _{K_{1}}%
=u_{1}+q_{1}$ and $\left.  \tilde{u}\right\vert _{K_{2}}=u_{2}$, is continuous
across $T$. On the other hand, all functions $u\in S_{T}^{p}$ whose
restrictions $\left.  u\right\vert _{\omega_{T}}$ are discontinuous can be
found in $S_{K_{1},\operatorname*{nc}}^{p}\oplus S_{K_{2},\operatorname*{nc}%
}^{p}$. In view of the direct sum in (\ref{directsplitfacetspace}) we may thus
assume that the functions in $\tilde{S}_{T,\operatorname*{nc}}^{p}$ are
continuous in $\omega_{T}$.

To finally arrive at a direct decomposition of the space in the right-hand
side of (\ref{directsplitfacetspace}) we have to split the spaces
$\mathbb{P}_{p,p-1}^{\perp}\left(  C_{i}\right)  $ into a direct sum of the
spaces of totally symmetric orthogonal polynomials and the spaces introduced
in Definition \ref{RemMirror} and glue them together in a continuous way. We
introduce the functions $b_{p,k}^{C_{i},\operatorname*{sym}}:=\left.
b_{p,k}^{\partial K_{i},\operatorname*{sym}}\right\vert _{C_{i}}$, $0\leq
k\leq d_{\operatorname*{triv}}\left(  p\right)  -1$, with $b_{p,k}^{\partial
K_{i},\operatorname*{sym}}$ as in (\ref{defsurfacesym}) and define
$b_{p,k}^{C_{i},\operatorname*{refl}}$, $0\leq k\leq d_{\operatorname*{refl}%
}\left(  p\right)  -1$, piecewise by $\left.  b_{p,k}^{C_{i}%
,\operatorname*{refl}}\right\vert _{T^{\prime}}:=b_{p,k}^{\mathbf{A}%
_{i},T^{\prime}}$ for $T^{\prime}\subset C_{i}$ with $b_{p,k}^{\mathbf{A}%
_{i},T^{\prime}}$ as in (\ref{defbmitttriangle}). The mirror symmetry of
$b_{p,k}^{\mathbf{A}_{i},T^{\prime}}$ with respect to the angular bisector in
$T^{\prime}$ through $\mathbf{A}_{i}$ implies the continuity of $b_{p,k}%
^{C_{i},\operatorname*{refl}}$. Hence,%
\begin{equation}
\mathbb{P}_{p,p-1}^{\perp}\left(  C_{i}\right)  =\operatorname*{span}\left\{
\left.  b_{p,k}^{C_{i},\operatorname*{sym}}\right\vert _{C_{i}}:0\leq k\leq
d_{\operatorname*{triv}}\left(  p\right)  -1\right\}  \oplus
\operatorname*{span}\left\{  b_{p,k}^{C_{i},\operatorname*{refl}}:0\leq k\leq
d_{\operatorname*{refl}}\left(  p\right)  -1\right\}  . \label{paux}%
\end{equation}
Since the traces of $b_{p,k}^{C_{i},\operatorname*{sym}}$ and $b_{p,k}%
^{C_{i},\operatorname*{refl}}$ at $\partial T$ are continuous and are, from
both sides, the same linear combinations of edge-wise Legendre polynomials of
even degree, the gluing $\left.  b_{p,k}^{\partial\omega_{T}%
,\operatorname*{sym}}\right\vert _{C_{i}}:=b_{p,k}^{C_{i},\operatorname*{sym}%
}$ and $\left.  b_{p,k}^{\partial\omega_{T},\operatorname*{refl}}\right\vert
_{\dot{C}_{i}}:=b_{p,k}^{C_{i},\operatorname*{refl}}$, $i=1,2$, defines
continuous functions on $\partial\omega_{T}$. Since the space
$S_{T,\operatorname*{nc}}^{p}$ must satisfy a direct sum decomposition (cf.
(\ref{directsplitfacetspace})), it suffices to consider the functions
$b_{p,k}^{\partial\omega_{T},\operatorname*{refl}}$ for the definition of
$S_{T,\operatorname*{nc}}^{p}$. The resulting non-conforming facet-oriented
space $S_{T,\operatorname*{nc}}^{p}$ was introduced in Definition
\ref{PropBpmTncbasis} and $\tilde{S}_{T,\operatorname*{nc}}^{p}$ can be chosen
to be $S_{T,\operatorname*{nc}}^{p}$.

\begin{proposition}
\label{CornonorthoT}For any $u\in S_{T,\operatorname*{nc}}^{p}$, the following
implication holds%
\[
\left.  u\right\vert _{T}\in\left.  S_{T,\operatorname*{nc}}^{p}\right\vert
_{T}\cap\mathbb{P}_{p,p-1}^{\perp}\left(  T\right)  \implies u=0.
\]

\end{proposition}%

\proof
Assume there exists $u\in S_{T,\operatorname*{nc}}^{p}$ with $\left.
u\right\vert _{T}\in\left.  S_{T,\operatorname*{nc}}^{p}\right\vert _{T}%
\cap\mathbb{P}_{p,p-1}^{\perp}\left(  T\right)  $. Let $K$ be a simplex
adjacent to $T$. Then $u_{K}=\left.  u\right\vert _{K}$ satisfies $\left.
u_{K}\right\vert _{T^{\prime}}\in\mathbb{P}_{p,p-1}^{\perp}\left(  T^{\prime
}\right)  $ for all $T^{\prime}\subset\partial K$ and, thus, $u_{K}\in
S_{K,\operatorname*{nc}}^{p}$. Since $\left.  S_{K,\operatorname*{nc}}%
^{p}\right\vert _{T^{\prime}}\cap\left.  S_{T,\operatorname*{nc}}%
^{p}\right\vert _{T^{\prime}}=\left\{  0\right\}  $ for $T^{\prime}\in\partial
K\backslash\overset{\circ}{T}$ we conclude that $u_{K}=0$.%
\endproof

Note that Definition \ref{PropBpmTncbasis} and Proposition \ref{CornonorthoT}
neither imply a priori that the functions
\[
\left.  B_{p,k}^{T,\operatorname*{nc}}\right\vert _{K},\quad\forall
T\subset\partial K,\quad k=0,\ldots,d_{\operatorname*{refl}}\left(  p\right)
-1
\]
are linearly independent nor that%
\begin{equation}
\forall T\subset\partial K\quad\text{it holds\quad}\sum_{T^{\prime}\subset
C}\left.  B_{p,m}^{T^{\prime},\operatorname*{nc}}\right\vert _{T}%
=\mathbb{P}_{p,p-1}^{\perp,\operatorname*{refl}}\left(  T\right)
\quad\text{for the triangle star }C=\partial K\backslash\overset{\circ}{T}
\label{reflspacecontained}%
\end{equation}
holds. These properties will be proved next. Recall the projection $\Pi
=\frac{1}{3}\left(  2I-MR-RM\right)  $ from Proposition \ref{Prop17}. We
showed (Theorem \ref{TheoCDBasis}.a) that$\left\{  b_{p,k}%
^{\operatorname*{refl}}:0\leq k\leq\frac{p-1}{3}\right\}  $ is linearly
independent, where $b_{p,k}^{\operatorname*{refl}}:=\Pi b_{p,2k}$.
Additionally $Rb_{p,k}^{\operatorname*{refl}}=b_{p,k}^{\operatorname*{refl}}$
which implies $b_{p,k}^{\operatorname*{refl}}\left(  0,x_{1}\right)
=b_{p,k}^{\operatorname*{refl}}\left(  x_{1},0\right)  $, and the restriction
$x_{1}\longmapsto b_{p,k}^{\operatorname*{refl}}\left(  x_{1},1-x_{1}\right)
$ is invariant under $x_{1}\mapsto1-x_{1}$. For four non-coplanar points
$A_{0},A_{1},A_{2},A_{3}$ let $K$ denote the tetrahedron with these vertices.
For any $k$ such that $0\leq k\leq\frac{p-1}{3}$ define a piecewise polynomial
on the faces of $K$ as follows: choose a local $\left(  x_{1},x_{2}\right)
$-coordinate system for $A_{0}A_{1}A_{2}$ so that the respective coordinates
are $\left(  0,0\right)  ,\left(  1,0\right)  ,\left(  0,1\right)  $, and
define $Q_{k}^{\left(  0\right)  }$ on the facet equal to $b_{p,k}%
^{\operatorname*{refl}}$. Similarly define $Q_{k}^{\left(  0\right)  }$ on
$A_{0}A_{2}A_{3}$ and $A_{0}A_{3}A_{1}$ (with analogously chosen local
$\left(  x_{1},x_{2}\right)  $-coordinate systems), by the property
$b_{p,k}^{\operatorname*{refl}}\left(  0,x_{1}\right)  =b_{p,k}%
^{\operatorname*{refl}}\left(  x_{1},0\right)  $. $Q_{k}^{\left(  0\right)  }$
is continuous at the edges $A_{0}A_{1}$, $A_{0}A_{2}$, and $A_{0}A_{3}$. The
values at the boundary of the triangle star equal $b_{p,k}%
^{\operatorname*{refl}}\left(  x_{1},1-x_{1}\right)  $; note the symmetry and
thus the orientation of the coordinates on the edges $A_{1}A_{2}$, $A_{2}%
A_{3}$, $A_{3}A_{1}$ is immaterial. The value of $Q_{k}^{\left(  0\right)  }$
on the triangle $A_{1}A_{2}A_{3}$ is taken to be a degree $p$ polynomial,
totally symmetric, with values agreeing with $b_{p,k}^{\operatorname*{refl}%
}\left(  x_{1},1-x_{1}\right)  $ on each edge.

Similarly $Q_{k}^{\left(  1\right)  },Q_{k}^{\left(  2\right)  }%
,Q_{k}^{\left(  3\right)  }$ are defined by taking $A_{1},A_{2},A_{3}$ as the
center of the construction, respectively.

\begin{theorem}
\label{THeoQkilinindep}a) The functions $Q_{k}^{\left(  i\right)  }$, $0\leq
k\leq d_{\operatorname*{refl}}\left(  p\right)  -1$, $i=0,1,2,3$ are linearly independent.

b) Property (\ref{reflspacecontained}) holds.
\end{theorem}

The proof involves a series of steps. The argument will depend on the values
of the functions on the three rays $A_{0}A_{1}$, $A_{0}A_{2}$, $A_{0}A_{3}$,
each one of them is given coordinates $t$ so that $t=0$ at $A_{0}$ and $t=1$
at the other end-point. For a fixed $k$ let $q\left(  t\right)  =b_{p,k}%
^{\operatorname*{refl}}\left(  t,0\right)  $, $\widehat{q}\left(  t\right)
=b_{p,k}^{\operatorname*{refl}}\left(  1-t,0\right)  $ and $\widetilde
{q}\left(  t\right)  =b_{p,k}^{\operatorname*{refl}}\left(  t,1-t\right)  $.

\begin{lemma}
\label{Lemsump=0}Suppose $0\leq k\leq\frac{p-1}{3}$ and $0\leq t\leq1$ then
$q\left(  t\right)  +\widehat{q}\left(  t\right)  +\widetilde{q}\left(
t\right)  =0$.
\end{lemma}

%

\proof
The actions of $RM$ and $MR$ on polynomials $f\left(  x_{1},x_{2}\right)  $
are given by $MRf\left(  x_{1},x_{2}\right)  =f\left(  1-x_{1}-x_{2}%
,x_{1}\right)  $ and $RMf\left(  x_{1},x_{2}\right)  =f\left(  x_{2}%
,1-x_{1}-x_{2}\right)  $. Polynomials of $\tau_{\operatorname*{refl}}$-type
satisfy $f+RMf+MRf=0$. Apply this relation to $b_{p,k}^{\operatorname*{refl}}$
with $x_{1}=t$ and $x_{2}=0$ with the result%
\[
b_{p,k}^{\operatorname*{refl}}\left(  t,0\right)  +b_{p,k}%
^{\operatorname*{refl}}\left(  1-t,t\right)  +b_{p,k}^{\operatorname*{refl}%
}\left(  0,1-t\right)  =0.
\]
The fact that $b_{p,k}^{\operatorname*{refl}}\left(  x_{1},x_{2}\right)
=b_{p,k}^{\operatorname*{refl}}\left(  x_{2},x_{1}\right)  $ finishes the
proof.%
\endproof

\textbf{Proof of Theorem \ref{THeoQkilinindep}. }Consider the contribution of
$Q_{k}^{\left(  1\right)  }$ to the values on the ray $A_{0}A_{1}$: because
$Q_{k}^{\left(  1\right)  }$ is constructed taking the origin at $A_{1}$ and
because of the reverse orientation of the ray we see that the value of
$Q_{k}^{\left(  1\right)  }$ is given by $\widehat{q}$. The value of
$Q_{k}^{\left(  1\right)  }$ on the ray $A_{0}A_{2}$ is $\widetilde{q}$ (by
the symmetry of $\widetilde{q}$ the orientation of the ray does not matter).
The other functions are handled similarly, and the contributions to the three
rays are given in this table:%
\[%
\begin{vmatrix}
& Q_{k}^{\left(  0\right)  } & Q_{k}^{\left(  1\right)  } & Q_{k}^{\left(
2\right)  } & Q_{k}^{\left(  3\right)  }\\
A_{0}A_{1} & q & \widehat{q} & \widetilde{q} & \widetilde{q}\\
A_{0}A_{2} & q & \widetilde{q} & \widehat{q} & \widetilde{q}\\
A_{0}A_{3} & q & \widetilde{q} & \widetilde{q} & \widehat{q}%
\end{vmatrix}
.
\]

We use $q_{k},\widetilde{q}_{k},\widehat{q}_{k}$ to denote the polynomials
corresponding to $b_{p,k}^{\operatorname*{refl}}$. Suppose that the linear
combination $\sum_{k=0}^{\left\lfloor \left(  p-1\right)  /3\right\rfloor
}\sum_{i=0}^{3}c_{k,i}Q_{k}^{\left(  i\right)  }=0$. Evaluate the sum on the
three rays to obtain the equations:%
\begin{align*}
0  &  =\sum_{k=0}^{\left\lfloor \left(  p-1\right)  /3\right\rfloor }\left\{
c_{k,0}q_{k}+c_{k,1}\widehat{q}_{k}+\left(  c_{k,2}+c_{k,3}\right)
\widetilde{q}_{k}\right\}  =\sum_{k=0}^{\left\lfloor \left(  p-1\right)
/3\right\rfloor }\left\{  \left(  c_{k,1}-c_{k,0}\right)  \widehat{q}%
_{k}+\left(  c_{k,2}+c_{k,3}-c_{k,0}\right)  \widetilde{q}_{k}\right\} \\
0  &  =\sum_{k=0}^{\left\lfloor \left(  p-1\right)  /3\right\rfloor }\left\{
c_{k,0}q_{k}+c_{k,2}\widehat{q}_{k}+\left(  c_{k,1}+c_{k,3}\right)
\widetilde{q}_{k}\right\}  =\sum_{k=0}^{\left\lfloor \left(  p-1\right)
/3\right\rfloor }\left\{  \left(  c_{k,2}-c_{k,0}\right)  \widehat{q}%
_{k}+\left(  c_{k,1}+c_{k,3}-c_{k,0}\right)  \widetilde{q}_{k}\right\}  ,\\
0  &  =\sum_{k=0}^{\left\lfloor \left(  p-1\right)  /3\right\rfloor }\left\{
c_{k,0}q_{k}+c_{k,3}\widehat{q}_{k}+\left(  c_{k,1}+c_{k,2}\right)
\widetilde{q}_{k}\right\}  =\sum_{k=0}^{\left\lfloor \left(  p-1\right)
/3\right\rfloor }\left\{  \left(  c_{k,3}-c_{k,0}\right)  \widehat{q}%
_{k}+\left(  c_{k,1}+c_{k,2}-c_{k,0}\right)  \widetilde{q}_{k}\right\}  .
\end{align*}
We used Lemma \ref{Lemsump=0} to eliminate $q_{k}$ from the equations. In
Theorem \ref{TheoCDBasis}.b we showed the linear independence of $\left\{
RMb_{p,k}^{\operatorname*{refl}},MRb_{p,k}^{\operatorname*{refl}}:0\leq
k\leq\frac{p-1}{3}\right\}  $, and in Lemma \ref{LemRestLinIndep} that the
restriction map $f\mapsto f\left(  x_{1},0\right)  $ is an isomorphism from
the orthogonal polynomials $\mathbb{P}_{p,p-1}^{\bot}$ to $\mathbb{P}%
_{p}\left(  \left[  0,1\right]  \right)  $. Thus the projection of the set is
also linearly independent, that is, $\left\{  \widetilde{q}_{k},\widehat
{q}_{k}:0\leq k\leq\frac{p-1}{3}\right\}  $ is a linearly independent set of
polynomials on $0\leq t\leq1$. This implies all the coefficients in the above
equations vanish: the $\widehat{q}_{k}$ terms show $c_{k,0}=c_{k,1}%
=c_{k,2}=c_{k,3}$ and then the $\widetilde{q}_{k}$-terms show $2c_{k,0}%
-c_{k,0}=c_{k,0}=0$.

To prove (\ref{reflspacecontained}) it suffices to transfer the statement to
the reference element $\widehat{T}$. The pullbacks of the restrictions
$\left.  B_{p,m}^{T^{\prime},\operatorname*{nc}}\right\vert _{T}$, $T^{\prime
}\subset C$, are given by%
\begin{equation}
b_{n,k}^{\operatorname*{refl}}=\Pi b_{n,2k}\text{,}\quad\widetilde{b}%
_{n,k}^{\operatorname*{refl}}:=RM\Pi b_{n,2k}\text{,\quad}\widehat{b}%
_{n,k}^{\operatorname*{refl}}:=MR\Pi b_{n,2k}\text{, }k=0,\ldots
d_{\operatorname*{refl}}\left(  n\right)  -1\text{.} \label{bbb}%
\end{equation}
Since $b_{n,k}^{\operatorname*{refl}}\in\mathbb{P}_{n,n-1}^{\perp
,\operatorname*{refl}}\left(  \widehat{T}\right)  $ (cf. (\ref{T1+T2})) it
follows%
\[
\mathbb{P}_{n,n-1}^{\perp,\operatorname*{refl}}\left(  \widehat{T}\right)
\overset{\text{(\ref{defPnn-1refl})}}{=}\operatorname*{span}\left\{
\widetilde{b}_{n,k}^{\operatorname*{refl}},\widehat{b}_{n,k}%
^{\operatorname*{refl}}:0\leq k\leq d_{\operatorname*{refl}}\left(  n\right)
-1\right\}  =\operatorname*{span}\left\{  b_{n,k}^{\operatorname*{refl}%
},\widetilde{b}_{n,k}^{\operatorname*{refl}},\widehat{b}_{n,k}%
^{\operatorname*{refl}}:0\leq k\leq d_{\operatorname*{refl}}\left(  n\right)
-1\right\}  .
\]%
\endproof

\section{Properties of Non-Conforming Crouzeix-Raviart Finite
Elements\label{SecPropNC}}

The non-conforming Crouzeix-Raviart finite element space $S_{\mathcal{G}%
,\operatorname*{nc}}^{p}$ satisfies $S_{\mathcal{G},\operatorname*{c}}%
^{p}\subsetneq S_{\mathcal{G},\operatorname*{nc}}^{p}\subset S_{\mathcal{G}%
}^{p}$ (cf. Section \ref{SecNCFECR}). In this section, we will present a basis
for $S_{\mathcal{G},\operatorname*{nc}}^{p}$ and discuss whether the inclusion
$S_{\mathcal{G},\operatorname*{nc}}^{p}\subset S_{\mathcal{G}}^{p}$, in fact,
is an equality.

\subsection{A Basis for Non-Conforming Crouzeix-Raviart Finite Elements}

We have defined conforming and non-conforming sets of functions which are
spanned by functions with local support. In this section, we will investigate
the linear independence of these functions. We introduce the following spaces%
\[
S_{\operatorname*{sym},\operatorname*{nc}}^{p}:=%
{\displaystyle\bigoplus\limits_{K\in\mathcal{G}}}
S_{K,\operatorname*{nc}}^{p},\quad S_{\operatorname*{refl},\operatorname*{nc}%
}^{p}:=%
{\displaystyle\bigoplus\limits_{T\in\mathcal{F}_{\Omega}}}
S_{T,\operatorname*{nc}}^{p},
\]
where $S_{K,\operatorname*{nc}}^{p}$ and $S_{T,\operatorname*{nc}}^{p}$ are as
in Definitions \ref{DefSymSpaceK} and \ref{PropBpmTncbasis}. For some $0\leq
k\leq d_{\operatorname*{refl}}\left(  p\right)  -1$, we introduce the subspace
$S_{\operatorname*{refl},\operatorname*{nc}}^{p,k}\subset
S_{\operatorname*{refl},\operatorname*{nc}}^{p}$ by%
\[
S_{\operatorname*{refl},\operatorname*{nc}}^{p,k}:=%
{\displaystyle\bigoplus\limits_{T\in\mathcal{F}_{\Omega}}}
\left\{  B_{p,m}^{T,\operatorname*{nc}}:0\leq m\leq k\right\}  .
\]
Further we will need the conforming finite element space (cf. (\ref{hpfinele}%
), Def. \ref{Deflocbasis}), where the vertex-oriented functions are removed,
i.e.,%
\[
\tilde{S}_{\mathcal{G},\operatorname*{c}}^{p}:=\left(
{\displaystyle\bigoplus\limits_{E\in\mathcal{E}_{\Omega}}}
S_{E,\operatorname*{c}}^{p}\right)  \oplus\left(
{\displaystyle\bigoplus\limits_{T\in\mathcal{F}_{\Omega}}}
S_{T,\operatorname*{c}}^{p}\right)  \oplus\left(
{\displaystyle\bigoplus\limits_{K\in\mathcal{G}}}
S_{K,\operatorname*{c}}^{p}\right)  .
\]

\begin{theorem}
\label{Theorem33}The sums%
\begin{equation}
S_{\mathcal{G},\operatorname*{c}}^{p}\oplus S_{\operatorname*{sym}%
,\operatorname*{nc}}^{p},\qquad S_{\operatorname*{sym},\operatorname*{nc}}%
^{p}\oplus S_{\operatorname*{refl},\operatorname*{nc}}^{p} \label{directsums}%
\end{equation}
are direct. The sum%
\begin{equation}
S_{\mathcal{G},\operatorname*{c}}^{p}+S_{\operatorname*{refl}%
,\operatorname*{nc}}^{p} \label{nondirectsum}%
\end{equation}
is \textbf{not} direct. The sum%
\begin{equation}
\tilde{S}_{\mathcal{G},\operatorname*{c}}^{p}\oplus S_{\operatorname*{sym}%
,\operatorname*{nc}}^{p}\oplus S_{\operatorname*{refl},\operatorname*{nc}%
}^{p,0} \label{finalspace}%
\end{equation}
is direct.
\end{theorem}

%

\proof
\textbf{Part 1.} We prove that the sum $S_{\operatorname*{sym}%
,\operatorname*{nc}}^{p}\oplus S_{\operatorname*{refl},\operatorname*{nc}}%
^{p}$ is direct.

From Proposition \ref{CornonorthoT} we know that the sum $\left.
S_{T,\operatorname*{nc}}^{p}\right\vert _{T}\oplus\mathbb{P}_{p,p-1}^{\perp
}\left(  T\right)  $ is direct. Let $\Pi_{T}:L^{2}\left(  T\right)
\rightarrow\mathbb{P}_{p-1}\left(  T\right)  $ denote the $L^{2}\left(
T\right)  $ orthogonal projection. Since $\mathbb{P}_{p-1}\left(  T\right)  $
is the orthogonal complement of $\mathbb{P}_{p,p-1}^{\perp}\left(  T\right)  $
in $\mathbb{P}_{p}\left(  T\right)  $ and since $\mathbb{P}_{p,p-1}^{\perp
}\left(  T\right)  \cap\left.  S_{T,\operatorname*{nc}}^{p}\right\vert
_{T}=\left\{  0\right\}  $, the restricted mapping $\Pi_{T}:\left.
S_{T,\operatorname*{nc}}^{p}\right\vert _{T}\rightarrow\mathbb{P}_{p-1}\left(
T\right)  $ is injective and the functions $q_{p,k}^{T}:=\Pi_{T}\left(
\left.  B_{p,k}^{T,\operatorname*{nc}}\right\vert _{T}\right)  $, $0\leq k\leq
d_{\operatorname*{refl}}\left(  p\right)  -1$, are linearly independent and
belong to $\mathbb{P}_{p-1}\left(  T\right)  $. We define the functionals%
\[
J_{p,k}^{T}\left(  w\right)  :=\int_{T}wq_{p,k}^{T}\quad0\leq k\leq
d_{\operatorname*{refl}}\left(  p\right)  -1\text{.}%
\]

Next we consider a general linear combination and show that the condition%
\begin{equation}
\sum_{K\subset\mathcal{G}}\sum_{i=0}^{d_{\operatorname*{triv}}\left(
p\right)  -1}\alpha_{i}^{K}B_{p,i}^{K,\operatorname*{nc}}+\sum_{K\subset
\mathcal{G}}\sum_{T^{\prime}\subset\partial K}\sum_{j=0}%
^{d_{\operatorname*{refl}}\left(  p\right)  -1}\beta_{j}^{T^{\prime}}%
B_{p,j}^{T^{\prime},\operatorname*{nc}}\overset{!}{=}0 \label{ureplinindep1}%
\end{equation}
implies that all coefficients are zero.

We apply the functionals $J_{p,k}^{T}$ to (\ref{ureplinindep1}) and use the
orthogonality between $\mathbb{P}_{p,p-1}^{\perp}\left(  T\right)  $ and
$q_{p,k}^{T}$ to obtain%
\begin{equation}
\sum_{K\subset\mathcal{G}}\sum_{T^{\prime}\subset\partial K}\sum
_{j=0}^{d_{\operatorname*{refl}}\left(  p\right)  -1}\beta_{j}^{T^{\prime}%
}J_{p,k}^{T}\left(  B_{p,j}^{T^{\prime},\operatorname*{nc}}\right)
\overset{!}{=}0. \label{ureplinindep2}%
\end{equation}
For $T^{\prime}\neq T$ it holds $J_{p,k}^{T}\left(  B_{p,i}^{T^{\prime
},\operatorname*{nc}}\right)  =0$ since $\left.  \left.  B_{p,i}^{T^{\prime
},\operatorname*{nc}}\right\vert _{K}\right\vert _{T}$ is an orthogonal
polynomial. Thus, equation (\ref{ureplinindep2}) is equivalent to%
\begin{equation}
\sum_{j=0}^{d_{\operatorname*{refl}}\left(  p\right)  -1}\beta_{j}^{T}%
J_{p,k}^{T}\left(  B_{p,j}^{T,\operatorname*{nc}}\right)  \overset{!}{=}0.
\label{ureplinindep3}%
\end{equation}
The matrix $\left(  J_{p,k}^{T}\left(  B_{p,j}^{T,\operatorname*{nc}}\right)
\right)  _{k,j=0}^{d_{\operatorname*{refl}}\left(  p\right)  -1}$ is regular
because%
\[
J_{p,k}^{T}\left(  B_{p,j}^{T,\operatorname*{nc}}\right)  =\int_{T}%
B_{p,j}^{T,\operatorname*{nc}}q_{p,k}^{T}=\int_{T}B_{p,j}%
^{T,\operatorname*{nc}}\Pi_{T}\left(  \left.  B_{p,k}^{T,\operatorname*{nc}%
}\right\vert _{T}\right)  =\int_{T}B_{p,j}^{T,\operatorname*{nc}}%
B_{p,k}^{T,\operatorname*{nc}}%
\]
and $\left(  \left.  B_{p,k}^{T,\operatorname*{nc}}\right\vert _{T}\right)
_{k}$ are linearly independent. Hence we conclude from (\ref{ureplinindep3})
that all coefficients $\beta_{j}^{T}$ are zero and the condition
(\ref{ureplinindep1}) reduces to%
\[
\sum_{K\subset\mathcal{G}}\sum_{i=0}^{d_{\operatorname*{triv}}\left(
p\right)  -1}\alpha_{i}^{K}B_{p,i}^{K,\operatorname*{nc}}\overset{!}{=}0.
\]
The left-hand side is a piecewise continuous function so that the condition is
equivalent to $\sum_{i=0}^{d_{\operatorname*{triv}}\left(  p\right)  -1}%
\alpha_{i}^{K}B_{p,i}^{K,\operatorname*{nc}}\overset{!}{=}0$ for all
$K\in\mathcal{G}$. Since $B_{p,i}^{K,\operatorname*{nc}}$ is a basis for
$\left.  S_{K,\operatorname*{nc}}^{p}\right\vert _{K}$ we conclude that all
$\alpha_{i}^{K}$ are zero.

\textbf{Part 2. }Next we prove that $\left(  S_{\operatorname*{sym}%
,\operatorname*{nc}}^{p}\oplus S_{\operatorname*{refl},\operatorname*{nc}%
}^{p,0}\right)  \cap\tilde{S}_{\mathcal{G},\operatorname*{c}}^{p}=\left\{
0\right\}  $ and we show this by contradiction. Let $u\in\left(
S_{\operatorname*{sym},\operatorname*{nc}}^{p}\oplus S_{\operatorname*{refl}%
,\operatorname*{nc}}^{p,0}\right)  \cap\tilde{S}_{\mathcal{G}%
,\operatorname*{c}}^{p}$ which satisfies $u\neq0$. We decompose
$u=u_{\operatorname*{sym}}+u_{\operatorname*{refl}}$ with
$u_{\operatorname*{sym}}\in S_{\operatorname*{sym},\operatorname*{nc}}^{p}$
and $u_{\operatorname*{refl}}\in S_{\operatorname*{refl},\operatorname*{nc}%
}^{p}$. We prove by contradiction that $u_{\operatorname*{sym}}\in
C^{0}\left(  \Omega\right)  $. Assume that $u_{\operatorname*{sym}}\notin
C^{0}\left(  \Omega\right)  $. Then, there exists a facet $T\subset
\mathcal{F}_{\Omega}$ such that $\left[  u_{\operatorname*{sym}}\right]
_{T}\neq0$. Then, $\left[  u_{\operatorname*{refl}}\right]  _{T}=-\left[
u_{\operatorname*{sym}}\right]  _{T}$ is a necessary condition for the
continuity of $u$. However, $\left[  u_{\operatorname*{sym}}\right]  _{T}%
\in\mathbb{P}_{p,p-1}^{\perp,\operatorname*{sym}}\left(  T\right)  $ while
$\left[  u_{\operatorname*{refl}}\right]  _{T}\in\mathbb{P}_{p,p-1}%
^{\perp,\operatorname*{refl}}\left(  T\right)  $ and there is a contradiction
because $\mathbb{P}_{p,p-1}^{\perp,\operatorname*{sym}}\left(  T\right)
\cap\mathbb{P}_{p,p-1}^{\perp,\operatorname*{refl}}\left(  T\right)  =\left\{
0\right\}  $. Hence, $u_{\operatorname*{sym}}\in C^{0}\left(  \Omega\right)  $
and, in turn, $u_{\operatorname*{refl}}\in C^{0}\left(  \Omega\right)  $.
\medskip

Since $u\neq0$, at least, one of the functions $u_{\operatorname*{sym}}$ and
$u_{\operatorname*{refl}}$ must be different from the zero function.

\textbf{Case a. }We show $u_{\operatorname*{sym}}=0$ by contradiction: Assume
$u_{\operatorname*{sym}}\neq0$. Then, $\left.  u_{\operatorname*{sym}%
}\right\vert _{T}\neq0$ for all facets $T\in\mathcal{F}$. (Proof by
contradiction: If $\left.  u_{\operatorname*{sym}}\right\vert _{T}=0$ for some
$T\in\mathcal{F}$, we pick some $K\in\mathcal{F}$ which has $T$ as a facet.
Since $\left.  u_{\operatorname*{sym}}\right\vert _{K}\in\left.
S_{K,\operatorname*{nc}}^{p}\right\vert _{K}$ we have $\left.
u_{\operatorname*{sym}}\right\vert _{T^{\prime}}=0$ for all facets $T^{\prime
}$ of $K$ and $\left.  u_{\operatorname*{sym}}\right\vert _{K}=0$. Since
$u_{\operatorname*{sym}}$ is continuous in $\Omega$, the restriction $\left.
u_{\operatorname*{sym}}\right\vert _{K^{\prime}}$ is zero for any $K^{\prime
}\in\mathcal{G}$ which shares a facet with $K$. This argument can be applied
inductively to show that $u_{\operatorname*{sym}}=0$ in $\Omega$. This is a
contradiction.) We pick a boundary facet $T\in\mathcal{F}_{\partial\Omega}$.
The condition $u\in\tilde{S}_{\mathcal{G},\operatorname*{c}}^{p}$ implies
$u=0$ on $\partial\Omega$ and, in particular, $\left.  u\right\vert
_{T}=\left.  u_{\operatorname*{sym}}\right\vert _{T}+\left.
u_{\operatorname*{refl}}\right\vert _{T}=0$. We use again the argument
$\mathbb{P}_{p,p-1}^{\perp,\operatorname*{sym}}\left(  T\right)
\cap\mathbb{P}_{p,p-1}^{\perp,\operatorname*{refl}}\left(  T\right)  =\left\{
0\right\}  $ which implies $u_{\operatorname*{sym}}=0$ and this is a
contradiction to the assumption $u_{\operatorname*{sym}}\neq0$.

\textbf{Case b. }From Case a we know that $u_{\operatorname*{sym}}=0$, i.e.,
$u_{\operatorname*{refl}}=u$, and it remains to show $u_{\operatorname*{refl}%
}=0$. The condition $u_{\operatorname*{refl}}\in\tilde{S}_{\mathcal{G}%
,\operatorname*{c}}^{p}$ implies $\left.  u_{\operatorname*{refl}}\right\vert
_{\partial\Omega}=0$ and $u_{\operatorname*{refl}}\left(  \mathbf{V}\right)
=0$ for all vertices $\mathbf{V}\in\mathcal{V}$.

The proof of Case b is similar than the proof of Case a and we start by
showing for a tetrahedron, say $K$, with a facet on the boundary that $\left.
u_{\operatorname*{refl}}\right\vert _{K}=0$ and employ an induction over
adjacent tetrahedrons to prove that $u_{\operatorname*{refl}}=0$ on every
tetrahedron in $\mathcal{G}$.

We consider a boundary facet $T_{0}\in\mathcal{F}_{\partial\Omega}$ with
adjacent tetrahedron $K\subset\mathcal{G}$. We denote the three other facets
of $K$ by $T_{i}$, $1\leq i\leq3$, and for $0\leq i\leq3$, the vertex of $K$
which is opposite to $T_{i}$ by $\mathbf{A}_{i}$.

\textbf{Case b.1. }First we consider the case that there is one and only one
\textit{other} facet, say, $T_{1}$ which lies in $\partial\Omega$. Then
$\left.  u_{\operatorname*{refl}}\right\vert _{T}=\left.  u_{2}\right\vert
_{T}+\left.  u_{3}\right\vert _{T}$ for some $u_{i}\in S_{T_{i}%
,\operatorname*{nc}}^{p,0}:=\operatorname*{span}\left\{  B_{p,0}%
^{T_{i},\operatorname*{nc}}\right\}  $, $i=2,3$. From Theorem
\ref{TheoCDBasis}.b we conclude that the sum $\left.  S_{T_{2}%
,\operatorname*{nc}}^{p,0}\right\vert _{T}\oplus\left.  S_{T_{3}%
,\operatorname*{nc}}^{p,0}\right\vert _{T}$ is direct. The condition $\left.
u_{\operatorname*{refl}}\right\vert _{T}=0$ then implies $u_{2}=u_{3}=0$.
Thus, we have proved $\left.  u_{\operatorname*{refl}}\right\vert _{K}=0$.

\textbf{Case b.2. }The case that there are exactly two \textit{other} facets
which are lying in $\partial\Omega$ can be treated in a similar way.

\textbf{Case b.3. }Next, we consider the case that $T_{i}\in\mathcal{F}%
_{\Omega}$ for $i=1,2,3$. Note that $\left.  u_{\operatorname*{refl}%
}\right\vert _{T}=\sum_{i=1}^{3}\left.  u_{i}\right\vert _{T}$ for some
$u_{i}\in S_{T_{i},\operatorname*{nc}}^{p,0}$. On $T$ we choose a local
$\left(  x_{1},x_{2}\right)  $-coordinate system such that $\mathbf{A}%
_{1}=\mathbf{0}$, $\mathbf{A}_{2}=\left(  1,0\right)  ^{\intercal}$,
$\mathbf{A}_{3}=\left(  0,1\right)  ^{\intercal}$. From (\ref{defqnk}) and
(\ref{defqnkplus}) we conclude that
\[
b_{n,k}^{\operatorname*{refl}}+RMb_{n,k}^{\operatorname*{refl}}+MRb_{n,k}%
^{\operatorname*{refl}}=0.
\]
This implies $\left.  u_{2}\right\vert _{T}=RM\left(  \left.  u_{1}\right\vert
_{T}\right)  =\left.  u_{1}\right\vert _{T}\circ\chi_{\left\{  3,2,1\right\}
}$ and $\left.  u_{3}\right\vert _{T}=MR\left(  \left.  u_{1}\right\vert
_{T}\right)  =\left.  u_{1}\right\vert _{T}\circ\chi_{\left\{  2,3,1\right\}
}$ (cf. (\ref{defchipi})) and, in turn, that the restrictions $u_{i}^{E}$ of
$u_{i}$ to the edge $E_{i}=T_{i}\cap T_{0}$, $1\leq i\leq3$, are the
\textquotedblleft same\textquotedblright, more precisely, the affine pullbacks
of $u_{i}^{E}$ to the interval $\left[  0,1\right]  $ are the same. From Lemma
\ref{TraceLemma}, we obtain that%
\begin{equation}
\left.  u_{1}\right\vert _{T_{1}}\circ\chi_{1}=\left.  u_{2}\right\vert
_{T_{2}}\circ\chi_{2}=\left.  u_{3}\right\vert _{T_{3}}\circ\chi_{3},
\label{affequi}%
\end{equation}
where $\chi_{i}:\widehat{T}\rightarrow T_{i}$ are affine pullbacks to the
reference triangle such that $\chi_{i}\left(  \mathbf{0}\right)
=\mathbf{A}_{0}$.

This implies that the functions $u_{i}$ at $\mathbf{A}_{0}$ have the same
value (say $w_{0}$) and, from the condition $u_{\operatorname*{refl}}\left(
\mathbf{A}_{0}\right)  =3w_{0}=0$, we conclude that $u_{i}\left(
\mathbf{A}_{0}\right)  =0$. The values of $u_{i}$ at the vertex $\mathbf{A}%
_{i}$ of $K$ (which is opposite to $T_{i}$) also coincide and we denote this
value by $v_{0}$. Since $\left.  u_{\operatorname*{refl}}\right\vert _{T}=0$
it holds $u_{\operatorname*{refl}}\left(  \mathbf{A}_{i}\right)  =2w_{0}%
+v_{0}=0$. From $w_{0}=0$ we conclude that also $v_{0}=0$. Let $\chi_{i,T_{0}%
}:\widehat{T}\rightarrow T_{0}$ denote an affine pullback with the property
$\chi_{i,T_{0}}\left(  \mathbf{0}\right)  =\mathbf{A}_{i}$. Hence,
\begin{equation}
\widehat{u_{i}}:=\left.  u_{i}\right\vert _{T_{0}}\circ\chi_{i,T_{0}}^{-1}%
\in\operatorname*{span}\left\{  b_{p,0}^{\operatorname*{refl}}\right\}
\label{defuihut}%
\end{equation}
with values zero at the vertices of $\hat{T}$. Note that%
\begin{equation}
b_{p,0}\left(  0,0\right)  =\left(  -1\right)  ^{p}\left(  p+1\right)
\quad\text{and\quad}b_{p,0}\left(  1,0\right)  =b_{p,0}\left(  0,1\right)  =1.
\label{bp000}%
\end{equation}
The vertex properties (\ref{bp000}) along the definition of $b_{p,k}%
^{\operatorname*{refl}}$ (cf. (\ref{defqnk})) imply that%
\begin{align}
b_{p,0}^{\operatorname*{refl}}\left(  1,0\right)   &  =b_{p,0}%
^{\operatorname*{refl}}\left(  0,1\right)  =\frac{1}{3}\left(  1-\left(
-1\right)  ^{p}\left(  p+1\right)  \right)  =c_{p},\label{defcp}\\
b_{p,0}^{\operatorname*{refl}}\left(  0,0\right)   &  =-2b_{p,0}%
^{\operatorname*{refl}}\left(  1,0\right)  .\nonumber
\end{align}
Since $c_{p}\neq0$ for $p\geq1$ we conclude that $\widehat{u}_{i}=0$ holds.
Relation (\ref{defuihut}) implies $\left.  u_{i}\right\vert _{T_{0}}=0$ and
thus $u_{i}=0$. From $\left.  u_{\operatorname*{refl}}\right\vert _{T}%
=\sum_{i=1}^{3}\left.  u_{i}\right\vert _{T}$ we deduce that $\left.
u_{\operatorname*{refl}}\right\vert _{K}=0$.

The Cases b.1-.3 allow to proceed with the same induction argument as for Case
a and $u_{\operatorname*{refl}}=0$ follows by induction.

\textbf{Part 3. }An inspection of Part 2 shows that, for the proof of Case a,
it was never used that the vertex-oriented basis functions have been removed
from $S_{\mathcal{G},\operatorname*{c}}^{p}$ and Case a holds verbatim for
$S_{\mathcal{G},\operatorname*{c}}^{p}$. This implies that the first sum in
(\ref{directsums}) is direct.

\textbf{Part 4.} The fact that the sum $S_{\mathcal{G},\operatorname*{c}}%
^{p}+S_{\operatorname*{refl},\operatorname*{nc}}^{p}$ is \textit{not} direct
is postponed to Proposition \ref{PropVertremoval}.%
\endproof

\begin{proposition}
\label{PropVertremoval}For any vertex $\mathbf{V\in}\mathcal{V}_{\Omega}$ it
holds $B_{p,\mathbf{V}}^{\mathcal{G}}\in S_{\operatorname*{sym}%
,\operatorname*{nc}}^{p}\oplus S_{\operatorname*{refl},\operatorname*{nc}%
}^{p,0}\oplus\tilde{S}_{\mathcal{G},\operatorname*{c}}^{p}$.
\end{proposition}

%

\proof
We will show the stronger statement $B_{p,\mathbf{V}}^{\mathcal{G}}\in
S_{\operatorname*{refl},\operatorname*{nc}}^{p,0}\oplus\tilde{S}%
_{\mathcal{G},\operatorname*{c}}^{p}$. It suffices to construct a continuous
function $u_{\mathbf{V}}\in S_{\operatorname*{refl},\operatorname*{nc}}^{p}$
which coincides with $B_{p,\mathbf{V}}^{\mathcal{G}}$ at all vertices
$\mathbf{V}^{\prime}\in\mathcal{V}$ and vanishes at $\partial\Omega$; then,
$B_{p,\mathbf{V}}^{\mathcal{G}}-u_{\mathbf{V}}\in\tilde{S}_{\mathcal{G}%
,\operatorname*{c}}^{p}$ and the assertion follows. Recall the known values of
$b_{p,0}^{\operatorname*{refl}}$ at the vertices of the reference triangle and
the definition of $c_{p}$ as in (\ref{defcp}). Let $K\in\mathcal{G}$ be a
tetrahedron with $\mathbf{V}$ as a vertex. The facets of $K$ are denoted by
$T_{i}$, $0\leq i\leq3$, and the vertex which is opposite to $T_{i}$ is
denoted by $\mathbf{A}_{i}$. As a convention we assume that $\mathbf{A}%
_{0}=\mathbf{V}$. For every $T_{i}$, $1\leq i\leq3$, we define the function
$u_{T_{i}}\in S_{T_{i},\operatorname*{nc}}^{p}$ by setting (cf.
(\ref{defbmitttriangle}))%
\[
\left.  u_{T_{i}}\right\vert _{T_{0}}=b_{p,0}^{\operatorname*{refl}}\circ
\chi_{\mathbf{A}_{i},T_{0}}^{-1},
\]
where $\chi_{\mathbf{A}_{i},T_{0}}:\widehat{T}\rightarrow T_{0}$ is an affine
pullback which satisfies $\chi_{\mathbf{A}_{i},T_{0}}\left(  \mathbf{0}%
\right)  =\mathbf{A}_{i}$. (It is easy to see that the definition of
$u_{T_{i}}$ is independent of the side of $T_{i}$, where the tetrahedron $K$
is located.) From (\ref{defqnk}) and (\ref{defqnkplus}) we conclude that
$\left.  \sum_{i=1}^{3}u_{T_{i}}\right\vert _{T_{0}}=0$ holds. We proceed in
the same way for all tetrahedrons $K\in\mathcal{G}_{\mathbf{V}}$ (cf.
(\ref{defoftrianglesubsets})). This implies that
\begin{equation}
\tilde{u}_{\mathbf{V}}:=\sum_{\substack{T\in\mathcal{F}_{\Omega}%
\\\mathbf{V}\in T}}u_{T} \label{sumutildev}%
\end{equation}
vanishes at $\Omega\backslash\overset{\circ}{\omega_{\mathbf{V}}}$ (cf.
(\ref{defoftrianglesubsets})). By construction the function $\tilde
{u}_{\mathbf{V}}$ is continuous. At $\mathbf{V}$, the function $u_{T_{i}}$ has
the value (cf. (\ref{defcp}))%
\[
u_{T_{i}}\left(  \mathbf{V}\right)  =c_{p}%
\]
so that $\tilde{u}_{\mathbf{V}}\left(  \mathbf{V}\right)  =Cc_{p}$, where $C$
is the number of terms in the sum (\ref{sumutildev}). Since $c_{p}>0$ for all
$p\geq1$, the function $u_{\mathbf{V}}:=\frac{1}{Cc_{p}}\tilde{u}_{V}$ is well
defined and has the desired properties.%
\endproof

\begin{remark}
We have seen that the extension of the basis functions of $S_{\mathcal{G}%
,\operatorname*{c}}^{p}$ by the basis functions of $S_{\operatorname*{refl}%
,\operatorname*{nc}}^{p}$ leads to linearly depending functions. On the other
hand, if the basis functions of the subspace $S_{\operatorname*{refl}%
,\operatorname*{nc}}^{p,0}$ are added and the vertex-oriented basis functions
in $S_{\mathcal{G},\operatorname*{c}}^{p}$ are simply removed, one arrives at
a set a linear independent functions which span a larger space than
$S_{\mathcal{G},\operatorname*{c}}^{p}$. Note that $S_{\operatorname*{refl}%
,\operatorname*{nc}}^{p,0}=S_{\operatorname*{refl},\operatorname*{nc}}^{p}$
for $p=1,2,3$.

One could add more basis functions from $S_{\operatorname*{refl}%
,\operatorname*{nc}}^{p}$ but then has to remove further basis functions from
$\tilde{S}_{\mathcal{G},\operatorname*{c}}^{p}$ or formulate side constraints
in order to obtain a set of linearly independent functions.
\end{remark}

We finish this section by an example which shows that there exist meshes with
fairly special topology, where the inclusion%
\begin{equation}
S_{\mathcal{G},\operatorname*{c}}^{p}+S_{\operatorname*{sym}%
,\operatorname*{nc}}^{p}+S_{\operatorname*{refl},\operatorname*{nc}}%
^{p}\subset S_{\mathcal{G}}^{p} \label{maxspace}%
\end{equation}
is strict. We emphasize that the left-hand side in (\ref{maxspace}), for
$p\geq4$, defines a larger space than the space in (\ref{finalspace}) since it
contains \textit{all} non-conforming functions of reflection type.

\begin{example}
Let us consider the octahedron $\Omega$ with vertices $\mathbf{A}^{\pm
}:=\left(  0,0,\pm1\right)  ^{\intercal}$ and $\mathbf{A}_{1}:=\left(
1,0,0\right)  ^{\intercal}$, $\mathbf{A}_{2}:=\left(  0,1,0\right)
^{\intercal}$, $\mathbf{A}_{3}:=\left(  -1,0,0\right)  ^{\intercal}$,
$\mathbf{A}_{4}:=\left(  0,-1,0\right)  ^{\intercal}$. $\Omega$ is subdivided
into a mesh $\mathcal{G}:=\left\{  K_{i}:1\leq i\leq8\right\}  $ consisting of
eight congruent tetrahedrons sharing the origin $\mathbf{0}$ as a common
vertex. The six vertices at $\partial\Omega$ have the special topological
property that each one belongs to exactly four surface facets.

Note that the space defined by the left-hand side of (\ref{maxspace}) does not
contain functions whose restriction to a surface facet, say $T$, belongs to
the $\tau_{\operatorname{sign}}$ component of $\mathbb{P}_{n,n-1}^{\perp
}\left(  T\right)  $. Hence, the inclusion in (\ref{maxspace}) is strict if we
identify a function in $S_{\mathcal{G}}^{p}$ whose restriction to some surface
facet is an orthogonal polynomial of \textquotedblleft$\operatorname{sign}$
type\textquotedblright. Let $\widehat{q}\neq0$ be a polynomial which belongs
to the $\tau_{\operatorname{sign}}$ component of $\mathbb{P}_{n,n-1}^{\perp
}\left(  T\right)  $ on the reference element. Denote the (eight) facet on
$\partial\Omega$ with the vertices $\mathbf{A}^{\pm}$, $\mathbf{A}_{i}$,
$\mathbf{A}_{i+1}$ by $T_{i}^{\pm}$ for $1\leq i\leq4$ (with cyclic numbering
convention) and choose affine pullbacks $\chi_{\pm,i}:\widehat{T}\rightarrow
T_{i}^{\pm}$ as $\chi_{\pm,i}\left(  \mathbf{x}\right)  :=\mathbf{A}^{\pm
}+x_{1}\left(  \mathbf{A}_{i}-\mathbf{A}^{\pm}\right)  +x_{2}\left(
\mathbf{A}_{i+1}-\mathbf{A}^{\pm}\right)  $. Then, it is easy to verify (use
Lemma \ref{Lemcont} with even $m_{C}$) that the function $q:\partial
\Omega\rightarrow\mathbb{R}$, defined by $\left.  q\right\vert _{T_{i}^{\pm}%
}:=\widehat{q}\circ\chi_{\pm,i}^{-1}$ is continuous on $\partial\Omega$. Hence
the \textquotedblleft finite element extension\textquotedblright\ to the
interior of $\Omega$ via
\[
Q:=\sum_{\mathbf{N}\in\mathcal{N}^{p}\cap\partial\Omega}q\left(
\mathbf{N}\right)  B_{p,\mathbf{N}}^{\mathcal{G}}%
\]
defines a function in $S_{\mathcal{G}}^{p}$ which is not in the space defined
by the left-hand side of (\ref{maxspace}).

We state in passing that the space $S_{\mathcal{G}}^{p}$ does not contain any
function whose restriction to a boundary facet, say $T$, belongs to the
$\tau_{\operatorname{sign}}$ component of $\mathbb{P}_{p,p-1}^{\perp}\left(
T\right)  $ if there exists at least one surface vertex which belongs to an
\emph{odd} number of surface facets. In this sense, the topological situation
considered in this example is fairly special.
\end{example}

\section{Conclusion\label{SecConclusion}}

In this article we developed explicit representation of a local basis for
non-conforming finite elements of the Crouzeix-Raviart type. As a model
problem we have considered Poisson-type equations in three-dimensional
domains; however, this approach is by no means limited to this model problem.
Using theoretical conditions in the spirit of the second Strang lemma, we have
derived conforming and non-conforming finite element spaces of arbitrary
order. For these spaces, we also derived sets of local basis functions. To the
best of our knowledge, such explicit representation for general polynomial
order $p$ are not available in the existing literature. The derivation
requires some deeper tools from orthogonal polynomials of triangles, in
particular, the splitting of these polynomials into three irreducible
irreducible $\mathcal{S}_{3}$ modules.

Based on these orthogonal polynomials, simplex- and facet-oriented
non-conforming basis functions are defined. There are two types of
non-conforming basis functions: those whose supports consist of one
tetrahedron and those whose supports consist of two adjacent tetrahedrons. The
first type can be simply added to the conforming $hp$ basis functions. It is
important to note that the span of the functions of the second type contains
also conforming functions and one has to remove some conforming functions in
order to obtain a linearly independent set of functions. We have proposed a
non-conforming space which consists of a) all basis functions of the first
type and b) a reduced set of basis functions of the second type and c) of the
conforming basis functions without the vertex-oriented ones. This leads to a
set of linearly independent functions and is in analogy to the well known
lowest order Crouzeix-Raviart element.

It is interesting to compare these results with high-order Crouzeix-Raviart
finite elements for the two-dimensional case which have been presented in
\cite{ccss_2012}. Facets $T$ of tetrahedrons in 3D correspond to edges $E$ of
triangles in 2D. As a consequence the dimension of the space of orthogonal
polynomials $\mathbb{P}_{p,p-1}^{\perp}\left(  E\right)  $ equals one. For
even degree $p$, one has only non-conforming basis functions of
\textquotedblleft symmetric\textquotedblright\ type (which are supported on a
single triangle) and for odd degree $p$, one has only non-conforming basis
functions of \textquotedblleft reflection\textquotedblright\ type (which are
supported on two adjacent triangles). It turns out that adding the non
conforming symmetric basis function to the conforming $hp$ finite element
space leads to a set of linearly independent functions which is the analogue
of the first sum in (\ref{directsums}). If the non-conforming basis functions
of reflection type are added, the set of vertex-oriented conforming basis
functions have to be removed from the conforming space. This is in analogy to
the properties (\ref{nondirectsum}) and (\ref{finalspace}).

Future research is devoted on numerical experiments and the application of
these functions to system of equations as, e.g., Stokes equation and the
Lam\'{e} system.\medskip

\noindent\textbf{Acknowledgement} This work was supported in part by ENSTA,
Paris, through a visit of S.A. Sauter during his sabbatical. This support is
gratefully acknowledged.

\appendix

\section{Alternative Sets of \textquotedblleft
Reflection-type\textquotedblright\ Basis Functions\label{AA}}

In this Appendix we define further sets of basis functions for the
$\tau_{\operatorname*{refl}}$ component of $\mathbb{P}_{n,n-1}^{\perp}\left(
\widehat{T}\right)  $ -- different choices might be preferable for different
kinds of applications. All these sets have in common that two vertices of
$\widehat{T}$ are special -- any basis function is symmetric/skew symmetric
with respect to the angular bisector of one of these two vertices.

\begin{remark}
The functions $b_{n,2k}$ can be characterized as the range of $I+R$. We
project these functions onto $\tau_{\operatorname*{refl}}$, that is, the space
$E^{\left(  \operatorname*{refl}\right)  }:=\left\{  p:RMp+MRp=-p\right\}  $.
Let
\[
T_{1}:=I-MR\quad\text{and\quad}T_{2}:=I-RM.
\]
The range of both is $E^{\left(  \operatorname*{refl}\right)  }$. We will show
that $\left\{  T_{1}b_{n,2k},T_{2}b_{n,2k},0\leq k\leq\left(  n-2\right)
/3\right\}  $ is a basis for $E^{\left(  \operatorname*{refl}\right)  }$.
Previously we showed $\left\{  RMq_{k},MRq_{k}\right\}  $ is a basis, where
$q_{k}=\left(  2I-MR-RM\right)  b_{n,2k}=\left(  T_{1}+T_{2}\right)  b_{n,2k}$
(cf. (\ref{defqnk}). Observe that%
\begin{align*}
RM\left(  2I-MR-RM\right)   &  =2RM-I-MR=T_{1}-2T_{2},\\
MR\left(  2I-MR-RM\right)   &  =2MR-RM-I=-2T_{1}+T_{2}%
\end{align*}
holds, so the basis is made up out of linear combinations of $\left\{
T_{1}b_{n,2k},T_{2}b_{n,2k},0\leq k\leq\left(  n-1\right)  /3\right\}  $.
These can be written as elements of the range of $T_{1}\left(  I+R\right)  $
and $T_{2}\left(  I+R\right)  $. Different linear combinations will behave
differently under the reflections $R,M,RMR$ (that is $\left(  x,y\right)
\rightarrow(y,x),(1-x-y,y),(x,1-x-y)$ respectively). After some computations
we find%
\begin{align}
R(T_{1}+T_{2})\left(  I+R\right)   &  =\left(  T_{1}+T_{2}\right)
(I+R),\label{T1+T2}\\
R(T_{1}-T_{2})\left(  I+R\right)   &  =-\left(  T_{1}-T_{2}\right)  \left(
I+R\right)  ,\nonumber\\
M(T_{1}-2T_{2})\left(  I+R\right)   &  =\left(  T_{1}-2T_{2}\right)  \left(
I+R\right)  ,\nonumber\\
MT_{1}\left(  I+R\right)   &  =-T_{1}\left(  I+R\right)  ,\nonumber\\
RMR(2T_{1}-T_{2})\left(  I+R\right)   &  =\left(  2T_{1}-T_{2}\right)  \left(
I+R\right)  ,\nonumber\\
RMRT_{2}\left(  I+R\right)   &  =-T_{2}\left(  I+R\right)  .\nonumber
\end{align}
Any two of these types can be used in producing bases from the $b_{n,2k}.$
Also each pair (first two, second two, third two) are orthogonal to each
other. Note $R$ fixes $\left(  0,0\right)  $ and reflects in the line $x=y$,
$M$ fixes $\left(  0,1\right)  $, reflects in $2x+y=1$, and $RMR$ fixes
$\left(  1,0\right)  $, reflects in $x+2y=1$.
\end{remark}

If we allow for a complex valued basis, the three vertices of $\widehat{T}$
can be treated more equally as can be seen from the following remark.

\begin{remark}
The basis functions can be complexified: set $\omega=\operatorname*{e}%
^{2\pi\operatorname*{i}/3}$; any polynomial in $E^{\left(
\operatorname*{refl}\right)  }$ can be expressed as $p=p_{1}+p_{2}$ such that
$MRp=\omega p_{1}+\omega^{2}p_{2}$ (consequently $RMp=\omega^{2}p_{1}+\omega
p_{2}$ ), then%
\begin{align*}
-\frac{1}{3}\left(  \omega T_{1}+\omega^{2}T_{2}\right)  p  &  =p_{1},\\
-\frac{1}{3}\left(  \omega^{2}T_{1}+\omega T_{2}\right)  p  &  =p_{2}.
\end{align*}
These lead to another basis built up from the $b_{n,2k}.$ Let%
\begin{align*}
S_{1}  &  =-\frac{1}{3}\left(  \omega T_{1}+\omega^{2}T_{2}\right)  \left(
I+R\right)  ,\\
S_{2}  &  =-\frac{1}{3}\left(  \omega^{2}T_{1}+\omega T_{2}\right)  \left(
I+R\right)  .
\end{align*}
Applying these operators to $b_{n,2k}$ produces a basis $\left\{
S_{1}b_{n,2k},S_{2}b_{n,2k}:0\leq k\leq\left(  n-1\right)  /3\right\}  $
satisfying%
\begin{align*}
RS_{1}b_{n,2k}  &  =S_{2}b_{n,2k},~RS_{2}b_{n,2k}=S_{1}b_{n,2k},\\
MS_{1}b_{n,2k}  &  =\omega S_{2}b_{n,2k},~MS_{2}b_{n,2k}=\omega^{2}%
S_{1}b_{n,2k},\\
RMRS_{1}b_{n,2k}  &  =\omega^{2}S_{2}b_{n,2k},~RMRS_{2}b_{n,2k}=\omega
S_{1}b_{n,2k}.
\end{align*}
This is a basis which behaves similarly at each vertex.
\end{remark}

\bibliographystyle{abbrv}
\bibliography{ciarlets_sauter_simian}

\end{document}